# ON SYSTEMS OF VERY WEAK UNARY DYADIC ARITHMETIC

Zlatan Damnjanovic
University of Southern California

Formal syntactic objects – terms, formulae, formal proofs, etc. – are naturally conceived as linearly concatenated strings of symbols, i.e. as words in a (usually) finite alphabet. Following Tarski, concatenation theory, the theory specifically intended to serve as a general theoretical framework for characterization and systematic study of such objects, is standardly formulated as the theory of the structure of a free semi-group determined by finitely many generators and an associative binary operation. Tarski originally presented the theory axiomatically with an induction axiom in second-order form (see [7]).

The words obtained from finitely many letters by repeated binary juxtaposition can also be thought of as constructed (less efficiently) by juxtapositions of a single letter $a_i$, $1 \leq i \leq n$, either as a prefix or a suffix to a given word, making use of two distinct unary juxtaposition operations ("front" and "end") associated with each letter $a_i$. Single letters $a_i$ may then be identified with the result of applying either one of the corresponding unary operations to the empty string Ø. That was roughly the approach of Hermes, (see [4]), who independently and about the same time as Tarski, introduced an axiomatic theory of concatenated strings, also in second-order form, allowing only the letter-prefixing operations. He termed these 'successor operations'. Considering the simple case of the words 1, 2, 11, 12, 21, 22, …. in the dyadic alphabet {1,2}, e.g. the binarily constructed string $(1\,\hat{}\,2)\,\hat{}\,1 = 1\,\hat{}\,(2\,\hat{}\,1)$ could obviously also be obtained in any one of the following ways:

$$1_\vee(2_\vee(1_\vee\emptyset)), \quad ((\emptyset\wedge 1)\wedge 2)\wedge 1, \quad (1_\vee(\emptyset\wedge 2))_\vee 1, \quad (1_\vee(2_\vee\emptyset))_\vee 1, \quad \ldots.$$

where the words $a_{i\vee}w$ and $w\wedge a_i$, resp., are obtained by appending the letter $a_i$ as a prefix and suffix, resp., to the word w.

Is there any essential difference between these two ways of formulating concatenation theory, one, Tarski's, as a theory of a free semi-group, and that of Hermes, based on the empty string and finitely many distinct "successor operations"? For $n = 1$, Hermes's theory is equivalent to second-order arithmetic because in the presence of the induction axiom $+$ and $\cdot$ can be defined from successor. Corcoran, Frank and Maloney [1] established that all of the theories considered by Tarski and Hermes are in fact formally mutually interpretable in one another. However, their proof essentially depends on the second-order format in which the theories are presented.

Recent studies of concatenation theory have focused on first-order, and, especially, on weak and very weak formulations of the theory, investigating relationships with analogously formulated theories of other mathematical structures such as integers, sets, sequences and finite trees of various arities, revealing in the process a great deal of information about expressive resources formally contained in those theories despite the apparent weakness of their deductive apparatus (for a brief survey see [8], [2]. Here we take up the question whether the fundamental equivalence between binary and unary formulations of concatenation theory obtains at this level of weak and very weak first-order theories in contexts characterized by complete absence of any axiom or schema of induction. More generally, is binary concatenation at all first-order definable in terms of unary concatenation? We shall prove that the answer to this question is positive modulo the presence of the initial and end segment relations between strings, thus solving a problem posed by Karlov in [5]. Broadly speaking, we study several very weak first-order theories of unary concatenation and investigate their relationships with other well-known very weak first-order theories. On account of the 1-1 correspondence between dyadic strings and positive integers the theories considered may be called systems of "unary dyadic arithmetic". We now describe these theories in detail.

Let $\Sigma^*$ be the set of (finite) words in the dyadic alphabet $\{\underline{1}, \underline{2}\}$, including the empty word e, and let $\sqsubseteq$ be the initial segment relation on dyadic words, i.e. $w_1 \sqsubseteq w_2 \equiv \exists w\, w_1 \hat{}\, w = w_2$ for $w_1, w_2 \in \Sigma^*$. We consider theories whose signatures include a single binary relation symbol $\sqsubseteq$ along with, possibly, infinitely many individual constants $\underline{c}_0, \underline{c}_1, \underline{c}_2,\ldots$, some from among the unary function symbols $S_1, S_2, S_3, S_4$ and $f$, and some of the additional binary relation symbols $\sqsupseteq$ and $|\leq|$. Terms and formulas are constructed in the usual way. Unless indicated otherwise, the vocabularies of these theories will be interpreted in the corresponding reducts of the structure

$$\mathfrak{S}^\dagger = \langle \Sigma^*, s_0, s_1, s_2,\ldots, f_1, f_2, f_3, f_4, f, \sqsubseteq, \sqsupseteq, |\leq| \rangle.$$

Here $s_0, s_1, s_2, \ldots$ stand for distinct words in $\Sigma^*$ relative to some given enumeration of $\Sigma^*$, with $s_0, s_1, s_2$ in particular standing for e, 1 and 2, resp., the function symbols $S_j$ signify the operations $f_j : \Sigma^* \rightarrow \Sigma^*$, for $j = 1,2$, of juxtaposing a single digit as the end-digit with $f_1(w) = w \hat{} \underline{1}$ and $f_2(w) = w \hat{} \underline{2}$, and for $j = 3,4$, of juxtaposing a single digit as the initial digit with $f_3(w) = \underline{1} \hat{} w$ and $f_4(w) = \underline{2} \hat{} w$, for any $w \in \Sigma^*$, the function symbol $f$ signifies the inverse operation $f : \Sigma^* \rightarrow \Sigma^*$ on dyadic words with $f(w) = w^{-1}$ for any $w \in \Sigma^*$, the binary relation symbols $\sqsubseteq$ and $\sqsupseteq$ express the initial and end segment relations $\boldsymbol{\sqsubseteq}$ and $\boldsymbol{\sqsupseteq}$ on $\Sigma^*$, resp., and $|\leq|$ expresses the relation $\boldsymbol{|\leq|}$ where $w_1 \boldsymbol{|\leq|} w_2$ holds just in case $w_1$ is any word in $\Sigma^*$ and $w_2$ a 2-tally (a string of consecutive 2's) such that $w_1$ does not exceed in length $w_2$, the length $lh(w)$ of a word $w \in \Sigma^*$ being the number of digits in w with $lh(e) = 0$. Unless indicated otherwise, we let $D = \{\underline{1},\underline{2}\}$.

Let $\mathcal{L}^+ = \{\underline{e}, S_1, S_2, f, \sqsubseteq, |\leq|\}$. We call a variable free $\mathcal{L}^+$ term t <u>canonical</u> if $S_1$ and $S_2$ are the only function symbols occurring in t. A straightforward induction on the length of words shows that for each $w \in \Sigma^*$ there is a (unique) canonical term t such that $w = val_{\mathfrak{S}+}\, t$, where $\mathfrak{S}^+$ is the reduct of $\mathfrak{S}^\dagger$ interpreting $\mathcal{L}^+$. We write $\underline{w}$ for the canonical term for $w \in \Sigma^*$. The first-order theory $VW^+$ has as its axioms, along with the usual axioms for identity, all instances of the schemas

(VW1) $\neg \underline{s} = \underline{t}$ for any distinct $s, t \in \Sigma^*$,

(VW2) $S_d(\underline{t}) = \underline{S_d(t)}$ for each $t \in \Sigma^*$, $d = 1,2$,

(VW3) $\underline{e} \sqsubseteq \underline{t}$ for each $t \in \Sigma^*$,

(VW4) $\forall x\,(x \sqsubseteq \underline{e} \rightarrow x = \underline{e})$

(VW5) $\forall x\,(x \sqsubseteq S_d(\underline{t}) \leftrightarrow x \sqsubseteq \underline{t} \vee x = S_d(\underline{t}))$ for each $t \in \Sigma^*$, $d = 1,2$,

(VW6) $\forall x\,(x \,|{\leq}|\, \underline{2}^m \leftrightarrow \bigvee_{lh(s)\leq m} x = \underline{s})$ for each $m \geq 0$,

(VW7) $\exists y\, y \,|{\leq}|\, \underline{t} \rightarrow \forall x\,(x \sqsubseteq \underline{t}\ \&\ x \neq \underline{e} \rightarrow \underline{2} \sqsubseteq f(x))$ for each $t \in \Sigma^*$.

(VW8) $f(\underline{t}) = \underline{t^{-1}}$ for each $t \in \Sigma^*$,

(VW9) $\forall x,y\,(f(x) = f(y) \rightarrow x = y)$

(VW10) $\forall x,y\,(x \sqsubseteq y\ \&\ y \sqsubseteq x) \rightarrow x = y)$

In (VW6), the notation '$\underline{2}^m$' stands for the result '$S_2(\ldots S_2(\underline{e})\ldots)$' of m iterated applications of $S_2$ to $\underline{e}$. It will be convenient to let '$\underline{2}^m$' stand for e when m = 0.

In §6 we introduce the theory VW* in the signature $\mathcal{L}^* = \{\underline{e}, S_1, S_2, S_3, S_4, \sqsubseteq, \sqsupseteq\}$ interpreted in the reduct $\mathfrak{S}^* = <\Sigma^*, e, f_1, f_2, f_3, f_4, \sqsubseteq, \sqsupseteq>$ of $\mathfrak{S}^\dagger$. The canonical $\mathcal{L}^*$-terms are defined in the same way as for $\mathcal{L}^+$. In addition to appropriately formulated identity axioms, the axioms of the first-order theory VW* are (VW1)-(VW5) plus all instances of the schemas

(VW2*) $S_d(\underline{t}) = \underline{S_d(t)}$ for each $t \in \Sigma^*$, $d = 3,4$,

(VW3*) $\underline{e} \sqsupseteq \underline{t}$ for each $t \in \Sigma^*$,

(VW4*) $\forall x\,(x \sqsupseteq \underline{e} \rightarrow x = \underline{e})$

(VW5*) $\forall x\,(x \sqsupseteq S_d(\underline{t}) \leftrightarrow x \sqsupseteq \underline{t} \vee x = S_d(\underline{t}))$ for each $t \in \Sigma^*$, $d = 3,4$,

along with (VW10).

In §7 we introduce the theory $VW_0$ in the signature $\mathcal{L}_0 = \{\underline{c}_0, \underline{c}_1, \underline{c}_2, \ldots, \sqsubseteq, \sqsupseteq\}$ interpreted in the reduct $\mathfrak{S}_0 = <\Sigma^*, s_0, s_1, s_2, \ldots, \sqsubseteq, \sqsupseteq>$ of $\mathfrak{S}^\dagger$. Each constant $\underline{c}_i$ may be identified with a distinct canonical term $\underline{s}$ for $s \in \Sigma^*$. Hence we shall simply write $\underline{s}, \underline{t}, \ldots$ in place of $\underline{c}_i, \underline{c}_j, \ldots$.. The axioms of the first-order theory $VW_0$ are all instances of the schemas

($VW_0$1) $\neg \underline{s} = \underline{t}$ for any distinct $s, t \in \Sigma^*$,

($VW_0$2) $\forall x\,(x \sqsubseteq \underline{t} \leftrightarrow \bigvee_{s \sqsubseteq t} x = \underline{s})$, for each $t \in \Sigma^*$,

($VW_0$3) $\forall x\,(x \sqsupseteq \underline{t} \leftrightarrow \bigvee_{s \sqsupseteq t} x = \underline{s})$, for each $t \in \Sigma^*$,

($VW_0$4) $\forall x\,(\forall y\,(y \sqsubset x \leftrightarrow y \sqsubset \underline{t}) \rightarrow \bigwedge_{d \in D} (\underline{d} \sqsupseteq x \,\&\, \underline{d} \sqsupseteq \underline{t} \rightarrow x = \underline{t}))$ for each $t \in \Sigma^*$,

($VW_0$5) $\forall x\,(\forall y\,(y \sqsupset x \leftrightarrow y \sqsupset \underline{t}) \rightarrow \bigwedge_{d \in D} (\underline{d} \sqsubseteq x \,\&\, \underline{d} \sqsubseteq \underline{t} \rightarrow x = \underline{t}))$ for each $t \in \Sigma^*$,

where $x \sqsubset y \equiv: x \sqsubseteq y \,\&\, x \neq y$ and $x \sqsupset y \equiv: x \sqsupseteq y \,\&\, x \neq y$, plus (VW10).

In §8 we introduce the theory VWDI in the signature $\mathcal{L}_{\sqsubseteq,f,c} = \{\underline{c}_0, \underline{c}_1, \underline{c}_2, \ldots, f, \sqsubseteq\}$ interpreted in the reduct $\mathfrak{S}_{\sqsubseteq,f,c} = <\Sigma^*, s_0, s_1, s_2, \ldots, f, \sqsubseteq>$ of $\mathfrak{S}^\dagger$, having as its axioms, along with the usual axioms for =, all instances of the schemas

(VWDI1) $\neg(\underline{s} = \underline{t})$ for any distinct $s, t \in \Sigma^*$,

(VWDI2) $f(\underline{s}) = \underline{s^{-1}}$ for each $s \in \Sigma^*$,

(VWDI3) $\forall x\,(x \sqsubseteq \underline{t} \leftrightarrow \bigvee_{s \sqsubseteq t} x = \underline{s})$ for each $s \in \Sigma^*$,

(VWDI4) $\forall x\; f(f(x)) = x$

(VWDI5) $\forall x\,(\forall y\,(y \sqsubset x \leftrightarrow y \sqsubset \underline{t}) \rightarrow \bigwedge_{d \in D} (\underline{d} \sqsubseteq f(x) \,\&\, \underline{d} \sqsubseteq f(\underline{t}) \rightarrow x = \underline{t}))$

for each $t \in \Sigma^*$

plus (VW10).

Summarizing the overall plan of the paper, after introducing in §1 some essential tools to be used in subsequent sections, we first show, in §§2-3 that $VW^+$ interprets a very weak theory of dyadic trees WT previously studied in [3]. In §4 we show how to define within $VW^+$ binary concatenation restricted to tallies and pseudo-tallies (strings of alternating digits). In §5 we use some of the tools from §4 to establish that a relational variant of a very weak theory WD of binary concatenation introduced in [6] is interpretable in $VW^+$. In §6 we show that VW* also interprets WD, and in §7 that a relational variant of VW* is in turn interpretable in $VW_0$. In §8 we define a direct interpretation of $VW_0$ in VWDI. Finally, in §9, by combining these together with previously known results we establish the mutual interpretability of all these theories among themselves as well with Robinson's weak arithmetic R. It therefore follows that all the theories of unary dyadic arithmetic considered here are essentially undecidable. In §9 we also show that binary concatenation is definable in the first-order theory $T_1$ of the structure $\langle \Sigma^*, 1, 2, \sqsubseteq, \sqsupseteq \rangle$ as well as the theory of the structure $\mathfrak{S}_{\sqsubseteq,f,c} = \langle \Sigma^*, s_0, s_1, s_2, \ldots, f, \sqsubseteq \rangle$. Thus, somewhat surprisingly, it turns out that when it comes to concatenation theory the fundamental role belongs to the partial order relations $\sqsubseteq$ and $\sqsupseteq$.

## §1. Some VW$^+$ preliminaries

We introduce some useful abbreviations. For any $\mathcal{L}^+$ term t not containing the variable x and any formula A, we let

$$\forall x \sqsubseteq t\, A \;\equiv:\; \forall x\,(x \sqsubseteq t \rightarrow A) \quad \text{and} \quad \exists x \sqsubseteq t\, A \;\equiv:\; \exists x\,(x \sqsubseteq t \,\&\, A).$$

Further, let $x \sqsupseteq y \;\equiv:\; f(x) \sqsubseteq f(y)$, ($\mathfrak{S}^+ \vDash$ "x is an endsegment of y")

$x \subseteq_p y \equiv:\; \exists z \sqsubseteq y\; x \sqsupseteq z$, ($\mathfrak{S}^+ \vDash$ "x is a (sub)segment of y")

Let $x \subseteq_1 y$ abbreviate

$$x \subseteq_p y \,\&\, \forall u,v \sqsubseteq y\,(x \sqsupseteq u \,\&\, x \sqsupseteq v \;\rightarrow\; u = v) \,\&\, \forall u,v \sqsupseteq y\,(x \sqsubseteq u \,\&\, x \sqsubseteq v \;\rightarrow\; u = v)$$

($\mathfrak{S}^+ \vDash$ "x has a single occurrence in y")

$L_1(x) \;=:\; f(S_1(f(x)))$, ($\mathfrak{S}^+ \vDash$ "$L_1(x) = \underline{1}\,\hat{}\,x$")

$L_2(x) \;=:\; f(S_2(f(x)))$, ($\mathfrak{S}^+ \vDash$ "$L_2(x) = \underline{2}\,\hat{}\,x$").

Bounded quantifiers '$\forall x \sqsupseteq$' and '$\exists x \sqsupseteq$' are introduced analogously to those for $\sqsubseteq$. We write '$\underline{1}$' and '$\underline{2}$' as shorthand for '$S_1(\underline{e})$', '$S_2(\underline{e})$'. For the sake of convenience and to highlight the fact that the vocabulary of $\mathcal{L}^+$ allows only for unary concatenation operations, we shall often, when writing $\mathcal{L}^+$ formulae, use the notations '$x\hat{}1$', '$x\hat{}2$', '$1_{\vee}x$', '$2_{\vee}x$', instead of the more cumbersome '$S_1(x)$', '$S_2(x)$', '$L_1(x)$', '$L_2(x)$', resp., 'x' being a variable. We define the class of $\mathcal{L}^+$ <u>concatenation terms</u> $C^+$ inductively as follows: (i) $\underline{e} \in C^+$, (ii) if $t \in C^+$, then $t\hat{}\underline{1}$ and $t\hat{}\underline{2} \in C^+$, (iii) if $t \in C^+$, then $\underline{1}_{\vee}t \in C^+$ and $\underline{2}_{\vee}t \in C^+$. Note that all canonical $\mathcal{L}^+$ terms are $\mathcal{L}^+$ concatenation terms.

We extend the ‘x^t‘ and ‘$t_v x$’ notation to any $\mathcal{L}^+$ concatenation term t by introducing *ad hoc* abbreviations by induction on the length of concatenation terms: ‘x^(t^d)’ is ‘(x^t)^d’ and ‘$(d\text{^}t)_v x$’ stands for ‘$d_v(t_v x)$’, whereas ‘$x\text{^}(d_v t)$’ is ‘(x^d)^t’ and ‘$(t\text{^}d)_v x$’ is ‘$t_v(d_v x)$’, for $d \in D$. When referring to dyadic strings as elements of the domain $\Sigma^*$ of the intended interpretation $\mathfrak{S}^+$ of $\mathcal{L}^+$ we will freely avail ourselves of the binary string concatenation operation $\frown$, sometimes omitting parentheses on account of associativity of $\frown$, writing, e.g., ‘$x \frown y \frown z$’ for ‘$(x \frown y) \frown z$’. Thus, for example, for any $t \in \Sigma^*$, we have

$$\mathrm{val}_{\mathfrak{S}+}\,(\underline{1}_v\underline{t})\text{^}\underline{2} = (\underline{1} \frown t) \frown \underline{2} = \underline{1} \frown (t \frown \underline{2}\,) = \mathrm{val}_{\mathfrak{S}+}\,\underline{1}_v(t\text{^}\,\underline{2}\,).$$

We proceed to establish some elementary facts about $VW^+$. Throughout this section and §§2-5 we let M be an arbitrary model of $VW^+$.

1.1 For each $t \in \Sigma^*$, $VW^+ \vdash \forall x\,(x \sqsubseteq \underline{t} \rightarrow \bigvee_{s \sqsubseteq t} x = \underline{s})$.

<u>Proof</u>: We argue by induction on the length of string t. The case t = e holds by (VW4). Assume the claim holds for t and consider $t \frown \underline{d}$ for $d \in D$. Assume $M \vDash x \sqsubseteq \underline{S_d(t)}$. By (VW2) then $M \vDash x \sqsubseteq S_d(\underline{t})$. By (VW5) we have $M \vDash x \sqsubseteq \underline{t}\ \vee\ x = S_d(\underline{t})$. From the induction hypothesis it follows that $M \vDash \bigvee_{s \sqsubseteq t} x = \underline{s}\ \vee\ x = S_d(\underline{t})$. But $M \vDash S_d(\underline{t}) \sqsubseteq S_d(\underline{t})$ by (VW5). Hence $M \vDash \bigvee_{s \sqsubseteq t \frown d} x = \underline{s}$, as required. □

1.2 Let u be any variable free term of $\mathcal{L}^+$. Then

$$VW^+ \vdash u = \underline{s} \quad \text{for some } s \in \Sigma^* \text{ such that } \mathrm{val}_{\mathfrak{S}+}\, u = s.$$

Proof: We argue by induction on the complexity of variable free terms of $\mathcal{L}^+$. The claim is trivial if u is $\underline{e}$. Assume the claim holds for u and consider $\phi(\underline{u})$, where $\phi$ is any one of the function symbols $S_1$, $S_2$, $f$. By the induction hypothesis we have that $VW^+ \vdash u = \underline{s}$ for some $s \in \Sigma^*$. Hence $VW^+ \vdash \phi(u) = \phi(\underline{s})$ by axioms of identity. But then $val_{\mathfrak{S}+}\, \phi(\underline{s})$ is $s\hat{}\,\underline{1}$, $s\hat{}\,\underline{2}$ or $s^{-1}$, and from (VW2) and (VW8) we have that $VW^+ \vdash \phi(\underline{s}) = \underline{\phi(s)}$. Hence $VW^+ \vdash \phi(u) = \underline{\phi(s)}$ as required since $val_{\mathfrak{S}+}\, \phi(\underline{s}) \in \Sigma^*$. □

We then have:

1.3 Let u, v be any variable-free terms of $\mathcal{L}^+$. If $\mathfrak{S}^+ \vDash u = v$, then $VW^+ \vdash u = v$.

Proof: From $\mathfrak{S}^+ \vDash u = v$ we have that $val_{\mathfrak{S}+}\, u = t = val_{\mathfrak{S}+}\, v$ for some $t \in \Sigma^*$. By 1.2 we have that $VW^+ \vdash u = \underline{t}$ and $VW^+ \vdash v = \underline{t}$, whence $VW^+ \vdash u = v$. □

We define, for each variable-free $\mathcal{L}^+$ term u its length lth'u' inductively: lth'$\underline{e}$' = 1 and lth' $S_1$(u)' = lth' $S_2$(u)' = lth'$f$(u)' = lth'u' + 1.

1.4 For any variable free $\mathcal{L}^+$ term u there is a variable free $\mathcal{L}^+$ term v such that

$$VW^+ \vdash f(u) = v \qquad \text{where} \quad \text{lth'v'} \leq \text{lth'u'}.$$

Proof: Let $s \in \Sigma^*$ be such that $val_{\mathfrak{S}+}\, u = s$. Then $s^{-1} \in \Sigma^*$ whereas $val_{\mathfrak{S}+}\, f(u) = s^{-1}$. Now, $VW^+ \vdash f(u) = \underline{s^{-1}}$ by 1.2, and since $\underline{s}$ has no occurrences of $f$ we have that lth'u' $\geq$ lth'$\underline{s}$' = lth' $\underline{s^{-1}}$', as required. □

1.5 Let u, v be any variable-free $\mathcal{L}^+$ terms. If $\mathfrak{S}^+ \vDash u \sqsubseteq v$, then $VW^+ \vdash u \sqsubseteq v$.

Proof: We argue by induction on lth'v'. Assume $\mathfrak{S}^+ \vDash u \sqsubseteq v$. If v is $\underline{e}$, then u must also be $\underline{e}$, and the claim holds by (VW3). Suppose v is $S_1(w)$. Let $w^* \in \Sigma^*$ be such that $val_{\mathfrak{S}+}\, w = w^*$. Then $\mathfrak{S}^+ \vDash u \sqsubseteq S_1(\underline{w^*})$, whence by (VW5) we have $\mathfrak{S}^+ \vDash u \sqsubseteq \underline{w^*} \lor u = S_1(\underline{w^*})$. If $\mathfrak{S}^+ \vDash u \sqsubseteq \underline{w^*}$ then $VW^+ \vdash u \sqsubseteq \underline{w^*}$ follows by the induction hypothesis because lth'$\underline{w^*}$' < lth'v'. But then $VW^+ \vdash u \sqsubseteq S_1(\underline{w^*})$ follows by (VW5), whereas $VW^+ \vdash w = \underline{w^*}$ by 1.2. So $VW^+ \vdash u \sqsubseteq S_1(w)$ as needed. If $\mathfrak{S}^+ \vDash u = S_1(\underline{w^*})$, then $VW^+ \vdash u = S_1(\underline{w^*})$ by 1.3, and, again by (VW5), we have that $VW^+ \vdash u \sqsubseteq S_1(\underline{w^*})$, whence $VW^+ \vdash u \sqsubseteq S_1(w)$ follows by 1.2. Exactly analogous argument applies if v is $S_2(w)$. If v is $f(w)$, then by 1.4 we have that $VW^+ \vdash f(w) = t$ where lth't' ≤ lth'w'. Since $\mathfrak{S}^+ \vDash u \sqsubseteq t$, we have by the induction hypothesis that $VW^+ \vdash u \sqsubseteq t$. But then $VW^+ \vdash u \sqsubseteq f(w)$. □

1.6 For any variable free $\mathcal{L}^+$ term u,

$$VW^+ \vdash \forall x\, (x \sqsubseteq u \leftrightarrow \textstyle\bigvee_{s \sqsubseteq u^*} x = \underline{s}) \quad \text{where} \quad u^* \in \Sigma^* \text{ and } u^* = val_{\mathfrak{S}+}\, u.$$

Proof: Assume $M \vDash x \sqsubseteq u$. By 1.2 we have that $VW^+ \vdash u = \underline{u^*}$. Hence $VW^+ \vdash x \sqsubseteq \underline{u^*}$, and so we have $VW^+ \vdash \bigvee_{s \sqsubseteq u^*} x = \underline{s}$ by 1.1. Conversely, assume $M \vDash x = \underline{s}$ for some $s \sqsubseteq u^*$. Then $\mathfrak{S}^+ \vDash \underline{s} \sqsubseteq \underline{u^*}$, hence $M \vDash \underline{s} \sqsubseteq \underline{u^*}$ by 1.5. But since $M \vDash u = \underline{u^*}$ by 1.2, it follows that $M \vDash \underline{s} \sqsubseteq u$, hence $M \vDash x \sqsubseteq u$. We thus also have $M \vDash \forall x\, (\bigvee_{s \sqsubseteq v^*} x = \underline{s} \rightarrow x \sqsubseteq u)$. □

1.7 For any variable-free $\mathcal{L}^+$ terms u, v, if $\mathfrak{S}^+ \vDash \neg u \sqsubseteq v$, then $VW^+ \vdash \neg u \sqsubseteq v$.

Proof: Assume $\mathfrak{S}^+ \vDash \neg u \sqsubseteq v$. Then $\mathfrak{S}^+ \vDash \neg \underline{w^*} = u$ for each $w^* \sqsubseteq v^*$ where $v^* \in \Sigma^*$ and $val_{\mathfrak{S}}\, v = v^*$. By (VW1) we have that $VW^+ \vdash \bigwedge_{w^* \sqsubseteq v^*} \neg \underline{u^*} = \underline{w^*}$ where $u^* \in \Sigma^*$ and $val_{\mathfrak{S}+}\, u = u^*$. Then we have $VW^+ \vdash \bigwedge_{w^* \sqsubseteq v^*} \neg u = \underline{w^*}$ by 1.2. But then $VW^+ \vdash \neg u \sqsubseteq v$ follows by 1.6. □

1.8 For any variable-free $\mathcal{L}^+$ terms u, v,

(a) if $\mathfrak{S}^+ \vDash u\ |\leq|\ v$, then $VW^+ \vdash u\ |\leq|\ v$,

(b) if $\mathfrak{S}^+ \vDash \neg u\ |\leq|\ v$, then $VW^+ \vdash \neg u\ |\leq|\ v$.

Proof: Let $u^* = val_{\mathfrak{S}+}\, u$ and let $v^* = val_{\mathfrak{S}+}\, v$. For (a), suppose $\mathfrak{S}^+ \vDash u\ |\leq|\ v$. Then $v^*$ is $2^n$ for some $n \geq 0$, and $lh(u^*) \leq n$. Hence from (VW6) we have that $M \vDash \underline{u^*}\ |\leq|\ \underline{v^*}$, whence $M \vDash u\ |\leq|\ v$ follows by 1.3. For (b), suppose $\mathfrak{S}^+ \vDash \neg u\ |\leq|\ v$. If $v^*$ is $2^n$ for some $n \geq 0$ then $lh(u^*) > n$. So $\mathfrak{S}^+ \vDash \bigwedge_{lh(s) \leq n} \underline{u}^* \neq \underline{s}$, hence by (VW1) we have $VW^+ \vdash \bigwedge_{lh(s) \leq n} \underline{u}^* \neq \underline{s}$, and from (VW6) we obtain $VW^+ \vdash \neg \underline{u}^*\ |\leq|\ \underline{v}$. Then $VW^+ \vdash \neg u\ |\leq|\ v$ follows by 1.3. If $v^*$ is not a 2-tally, then $\mathfrak{S}^+ \vDash \neg\, \underline{2} \sqsubseteq f(\underline{s})$ for some $s \sqsubseteq v^*$ where $s \neq e$. By 1.5, 1.7 and (VW1) we have $VW^+ \vdash \neg\, \underline{2} \sqsubseteq f(\underline{s})\ \&\ \underline{s} \sqsubseteq \underline{v}^*\ \&\ \underline{s} \neq \underline{e}$, hence $VW^+ \vdash \exists x\, (\neg\, \underline{2} \sqsubseteq f(x)\ \&\ x \sqsubseteq \underline{v}^*\ \&\ x \neq \underline{e})$, so $VW^+ \vdash \neg u\ |\leq|\ \underline{v}^*$ by (VW7). Then $VW^+ \vdash \neg u\ |\leq|\ v$ follows from 1.3. □

1.9 For any variable-free $\mathcal{L}^+$ terms u, v,

(a) if $\mathfrak{S}^+ \vDash u \sqsupseteq v$, then $VW^+ \vdash u \sqsupseteq v$,

(b) if $\mathfrak{S}^+ \vDash \neg u \sqsupseteq v$, then $VW^+ \vdash \neg u \sqsupseteq v$.

Proof: (a) Assume $\mathfrak{S}^+ \vDash u \sqsupseteq v$, i.e. $\mathfrak{S}^+ \vDash f(u) \sqsubseteq f(v)$. By 1.5 we have $VW^+ \vdash f(u) \sqsubseteq f(v)$, that is, $VW^+ \vdash u \sqsupseteq v$. For (b), assume $\mathfrak{S}^+ \vDash \neg u \sqsupseteq v$, i.e.

$\mathfrak{S}^+ \vDash \neg f(u) \sqsubseteq f(v)$. Then $VW^+ \vdash \neg f(u) \sqsubseteq f(v)$ by 1.5, that is, $VW^+ \vdash \neg u \sqsupseteq v$, as needed. □

1.10 For any variable-free $\mathcal{L}^+$ term t,

(a) $VW^+ \vdash \forall x\,(x \sqsubseteq t \rightarrow f(f(x)) = x)$.

(b) $VW^+ \vdash \forall x\,(x \sqsupseteq t \rightarrow f(f(x)) = x)$.

Proof: Assume $M \vDash x \sqsubseteq t$ for $x \in M$. Then $M \vDash \bigvee_{s \sqsubseteq t^*} x = \underline{s}$ by 1.1 where $t^* \in \Sigma^*$ and $t^* = \mathrm{val}_{\mathfrak{S}+}\, t$. Now, $M \vDash \bigwedge_{s \sqsubseteq t^*} f(f(\underline{s})) = f(\underline{s}^{-1}) = \underline{s}$ by (VW8). Hence $M \vDash \bigvee_{s \sqsubseteq t^*} (x = \underline{s}\ \&\ f(f(x)) = x)$. But then $M \vDash \forall x\,(x \sqsubseteq t \rightarrow f(f(x)) = x)$. For (b), assume $M \vDash x \sqsupseteq t$ for $x \in M$, that is, $M \vDash f(x) \sqsubseteq f(t)$. By 1.1 we have that $M \vDash \bigvee_{s \sqsubseteq f(t^*)} f(x) = \underline{s}$. Then $M \vDash \bigvee_{s \sqsubseteq f(t^*)} f(f(x)) = f(\underline{s})$. By (VW8) and 1.3 it follows that $M \vDash \bigvee_{s \sqsubseteq f(t^*)} f(x) = f(\underline{s^{-1}})$, whence $M \vDash \bigvee_{s \sqsubseteq f(t^*)} x = \underline{s^{-1}}$ by (VW9). But then, again by (VW8), we obtain $M \vDash f(f(x)) = x$, as needed. Then (b) follows from (a) using (VW9). □

1.11 For any variable-free $\mathcal{L}^+$ term u,

$$VW^+ \vdash \forall x\,(x \sqsupseteq u \leftrightarrow \bigvee_{s \sqsupseteq u^*} x = \underline{s}) \quad \text{where} \quad u^* \in \Sigma^* \text{ and } u^* = \mathrm{val}_{\mathfrak{S}+}\, u.$$

Proof: By 1.6 we have $VW^+ \vdash \forall x\,(x \sqsubseteq f(u) \leftrightarrow \bigvee_{s \sqsubseteq (u^*)\text{-}1} x = \underline{s})$, and further

$$VW^+ \vdash \forall x\,(f(x) \sqsubseteq f(u) \leftrightarrow \bigvee_{s \sqsubseteq (u^*)\text{-}1} f(x) = \underline{s}),$$

which means $VW^+ \vdash \forall x\,(x \sqsupseteq u \leftrightarrow \bigvee_{s \sqsubseteq (u^*)\text{-}1} f(x) = \underline{s})$. From $M \vDash \bigvee_{s \sqsubseteq (u^*)\text{-}1} f(x) = \underline{s}$ we have $M \vDash \bigvee_{s \sqsubseteq (u^*)\text{-}1} f(f(x)) = f(\underline{s})$. Now, assuming $M \vDash x \sqsupseteq u$, from 1.10(b) it follows that $M \vDash \bigvee_{s \sqsubseteq (u^*)\text{-}1} x = f(\underline{s})$. Hence from (VW8) we obtain $M \vDash \bigvee_{s \sqsupseteq u^*} x = \underline{s}$. Then $M \vDash \forall x\,(x \sqsupseteq u \to \bigvee_{s \sqsupseteq u^*} x = \underline{s})$. Conversely, assuming $M \vDash \bigvee_{s \sqsupseteq u^*} x = \underline{s}$ for $x \in M$, we have that $\mathfrak{S}^+ \vDash x \sqsupseteq u$, so $M \vDash x \sqsupseteq u$ by 1.9(a). Therefore also $M \vDash \forall x\,(\bigvee_{s \sqsupseteq u^*} x = \underline{s} \to x \sqsupseteq u)$. □

1.12 Let $\varphi$, $\psi$ be any $\mathcal{L}^+$ formulae with x as sole free variable, and let $t \in \Sigma^*$.

Then (a) $VW^+ \vdash \forall x\,(x \sqsubseteq \underline{t} \to \varphi(x)) \leftrightarrow \bigwedge_{s \sqsubseteq t} \varphi(\underline{s})$,

(b) $VW^+ \vdash \exists x\,(x \sqsubseteq \underline{t}\ \&\ \psi(x)) \leftrightarrow \bigvee_{s \sqsubseteq t} \psi(\underline{s})$.

(c) $VW^+ \vdash \forall x\,(x \sqsupseteq \underline{t} \to \varphi(x)) \leftrightarrow \bigwedge_{s \sqsupseteq t} \varphi(\underline{s})$,

(d) $VW^+ \vdash \exists x\,(x \sqsupseteq \underline{t}\ \&\ \psi(x)) \leftrightarrow \bigvee_{s \sqsupseteq t} \psi(\underline{s})$.

Proof: For (a), assume $M \vDash \forall x\,(x \sqsubseteq \underline{t} \to \varphi(x))$, and let $s \sqsubseteq t$. By 1.5 we have that $M \vDash \underline{s} \sqsubseteq \underline{t}$, whence $M \vDash \varphi(\underline{s})$ follows from hypothesis. Hence $M \vDash \bigwedge_{s \sqsubseteq t} \varphi(\underline{s})$. Conversely, suppose $M \vDash \bigwedge_{s \sqsubseteq t} \varphi(\underline{s})$, and assume $M \vDash x \sqsubseteq \underline{t}$ where $x \in M$. By 1.1 then $M \vDash \bigvee_{s \sqsubseteq t} x = \underline{s}$, whence $M \vDash \varphi(x)$. Hence it follows that $M \vDash \forall x\,(x \sqsubseteq \underline{t} \to \varphi(x))$. This proves (a). Then (b) follows from (a). Exactly analogous arguments establish (c) and (d) using 1.11. □

1.13 For any variable-free $\mathcal{L}^+$ term u,

$VW^+ \vdash \forall x\,(x \subseteq_p u \leftrightarrow \bigvee_{s \sqsubseteq u^*} \bigvee_{r \sqsupseteq s} x = \underline{r})$ where $u^* \in \Sigma^*$ and $u^* = val_{\mathfrak{S}^+}\, u$.

Proof: By 1.12(b) and 1.11 we have

$VW^+ \vdash \forall x\,(x \subseteq_p u \leftrightarrow \exists z \sqsubseteq u\ x \sqsupseteq z \leftrightarrow \bigvee_{s \sqsubseteq u^*} x \sqsupseteq \underline{s} \leftrightarrow \bigvee_{s \sqsubseteq u^*} \bigvee_{r \sqsupseteq s} x = \underline{r})$. □

1.14 Let u, v be any variable-free $\mathcal{L}^+$ terms.

(a) If $\mathfrak{S}^+ \vDash u \subseteq_p v$, then $VW^+ \vdash u \subseteq_p v$.

(b) If $\mathfrak{S}^+ \vDash \neg u \subseteq_p v$, then $VW^+ \vdash \neg u \subseteq_p v$.

Proof: For (a), assume $\mathfrak{S}^+ \vDash u \subseteq_p v$. Then, for $v^* \in \Sigma^*$ such that $val_{\mathfrak{S}+}\, v = v^*$, we have that $\mathfrak{S}^+ \vDash \underline{s} \sqsubseteq \underline{v^*}\ \&\ u \sqsupseteq \underline{s}$ for some $s \in \Sigma^*$. By 1.5 and 1.9(a) we obtain $VW^+ \vdash \underline{s} \sqsubseteq \underline{v^*}\ \&\ u \sqsupseteq \underline{s}$, whence $VW^+ \vdash \exists z(z \sqsubseteq \underline{v^*}\ \&\ u \sqsupseteq z)$, that is $VW^+ \vdash u \subseteq_p \underline{v^*}$. But then $VW^+ \vdash u \subseteq_p v$ follows from 1.2. For (b), assume $\mathfrak{S}^+ \vDash \neg u \subseteq_p v$. Then $\mathfrak{S}^+ \vDash \forall z \sqsubseteq v\ \neg u \sqsupseteq z$, hence $\mathfrak{S}^+ \vDash \bigwedge_{s \sqsubseteq v^*} \neg u \sqsupseteq \underline{s}$, whence, by 1.9(b), we have that $VW^+ \vdash \bigwedge_{s \sqsubseteq v^*} \neg u \sqsupseteq \underline{s}$. Now, $VW^+ \vdash \forall z\, (z \sqsubseteq \underline{v^*} \rightarrow \neg u \sqsupseteq z) \leftrightarrow \bigwedge_{s \sqsubseteq v^*} \neg u \sqsupseteq \underline{s}$ by 1.12(a). We thus derive $VW^+ \vdash \forall z\, (z \sqsubseteq \underline{v^*} \rightarrow \neg u \sqsupseteq z)$, that is, $VW^+ \vdash \neg u \subseteq_p \underline{v^*}$. But then $VW^+ \vdash \neg u \subseteq_p v$ follows by 1.2. □

1.15 Let φ, ψ be any $\mathcal{L}^+$ formulae with x as sole free variable, and u a variable-free $\mathcal{L}^+$ term. where $u^* \in \Sigma^*$ and $val_{\mathfrak{S}+}\, u = u^*$. Then:

(a) $VW^+ \vdash \forall x\, (x \subseteq_p u \rightarrow \varphi(x)) \leftrightarrow \bigwedge_{r \sqsubseteq u^*} \bigwedge_{s \sqsupseteq r} \varphi(\underline{s})$,

(b) $VW^+ \vdash \exists x\, (x \subseteq_p u\ \&\ \psi(x)) \leftrightarrow \bigvee_{r \sqsubseteq u^*} \bigvee_{s \sqsupseteq r} \psi(\underline{s})$.

Proof: For (a), from definitions we have that

$$M \vDash \forall x\, (x \subseteq_p u \rightarrow \varphi(x)) \Leftrightarrow M \vDash \forall x\, (\exists z \sqsubseteq u\ x \sqsupseteq z \rightarrow \varphi(x)) \Leftrightarrow$$

$$\Leftrightarrow M \vDash \forall x\, (\exists z\, (z \sqsubseteq u\ \&\ x \sqsupseteq z) \rightarrow \varphi(x)),\ \text{whence, by logic,}$$

$$\Leftrightarrow M \vDash \forall x\, \forall z\, (z \sqsubseteq u\ \&\ x \sqsupseteq z \rightarrow \varphi(x)) \Leftrightarrow M \vDash \forall z\, \forall x\, (z \sqsubseteq u\ \&\ x \sqsupseteq z \rightarrow \varphi(x)) \Leftrightarrow$$

$$\Leftrightarrow M \vDash \forall z\, (z \sqsubseteq u \rightarrow \forall x\, (x \sqsupseteq z \rightarrow \varphi(x))),\ \text{whence, by 1.12(a),}$$

$$\Leftrightarrow M \vDash \bigwedge_{r \sqsubseteq u^*} \forall x\, (x \sqsupseteq \underline{r} \rightarrow \varphi(x)),\ \text{and by 1.12(a),}$$

$\Leftrightarrow\ M \vDash \bigwedge_{r \sqsubseteq u^*} \bigwedge_{s \sqsupseteq r} \varphi(\underline{s})$.

Therefore $M \vDash \forall x\, (x \subseteq_p u \rightarrow \varphi(x)) \leftrightarrow \bigwedge_{r \sqsubseteq u^*} \bigwedge_{s \sqsupseteq r} \varphi(\underline{s})$.

An analogous argument establishes (b) using 1.12(b) instead. □

1.16 Let φ, ψ be any $\mathcal{L}^+$ formulae with x as sole free variable.

Then, for each $m \geq 0$,

(a) $VW^+ \vdash \forall x\, (x\, |{\leq}|\, \underline{2}^m \rightarrow \varphi(x)) \leftrightarrow \bigwedge_{lh(s) \leq m} \varphi(\underline{s})$,

(b) $VW^+ \vdash \exists x\, (x\, |{\leq}|\, \underline{2}^m\ \&\ \psi(x)) \leftrightarrow \bigvee_{lh(s) \leq m} \psi(\underline{s})$.

Let t be any variable-free $\mathcal{L}^+$ term other than a 2-tally. Then

(c) $VW^+ \vdash \forall x\, (x\, |{\leq}|\, t \rightarrow \varphi(x))$,

(d) $VW^+ \vdash \neg\exists x\, (x\, |{\leq}|\, t\ \&\ \psi(x))$.

Proof: (a) and (b) are proved analogously to 1.12(a)-(b) using (VW6). For (c) and (d), suppose t is a variable free $\mathcal{L}^+$ term that is not a 2-tally. Then, as in the proof of 1.8(b), we have that $VW^+ \vdash \exists y\, (\neg\underline{2} \sqsubseteq f(y)\ \&\ y \sqsubseteq t\ \&\ y \neq \underline{e})$. But then $VW^+ \vdash \forall x\, \neg x\, |{\leq}|\, t$ by (VW7). Then (c) and (d) follow immediately. □

We now introduce some additional useful abbreviatons. Let A be any $\mathcal{L}^+$ formula and t a $\mathcal{L}^+$ term not containing the variable x. Then

$\forall x \subseteq_p t\, A\ \equiv:\ \forall x\, (x \subseteq_p t \rightarrow A)$ and $\exists x \subseteq_p t\, A\ \equiv:\ \exists x\, (x \subseteq_p t\ \&\ A)$,

and $\forall x \,|{\leq}|\, t\, A \equiv: \forall x\,(x \,|{\leq}|\, t \rightarrow A)$ and $\exists x \,|{\leq}|\, t\, A \equiv: \exists x\,(x \,|{\leq}|\, t \,\&\, A)$

We define the class of <u>bounded $\mathcal{L}^+$ formulae</u> inductively as follows

(i) $u = v$, $u \sqsubseteq v$, $u \sqsupseteq v$, $u \subseteq_p v$ and $u \,|{\leq}|\, v$ are bounded $\mathcal{L}^+$ formulae for any $\mathcal{L}^+$ terms $u$, $v$,

(ii) if $\varphi$ and $\psi$ are bounded $\mathcal{L}^+$ formulae, then so are $\neg\varphi$, $(\varphi \,\&\, \psi)$ and $(\varphi \vee \psi)$,

(iii) if $\varphi$ is a bounded $\mathcal{L}^+$ formula and $t$ an $\mathcal{L}^+$ term not containing the variable $x$, then $\exists x \sqsubseteq t\, \varphi$, $\forall x \sqsubseteq t\, \varphi$, $\exists x \sqsupseteq t\, \varphi$, $\forall x \sqsupseteq t\, \varphi$, $\forall x \,|{\leq}|\, t\, \varphi$, $\exists x \,|{\leq}|\, t\, \varphi$, $\forall x \subseteq_p t\, \varphi$ and $\exists x \subseteq_p t\, \varphi$ are bounded $\mathcal{L}^+$ formulae.

1.17 Let $\varphi$ be any bounded $\mathcal{L}^+$ sentence:

(a) if $\mathfrak{S}^+ \vDash \varphi$, then $VW^+ \vdash \varphi$,

(b) if $\mathfrak{S}^+ \nvDash \varphi$, then $VW^+ \vdash \neg\varphi$.

<u>Proof</u>: We argue by induction on the complexity of bounded $\mathcal{L}^+$ sentences. If $\varphi$ is of the form '$u = v$', '$u \sqsubseteq v$', '$u \sqsupseteq v$', '$u \subseteq_p v$' or '$u \,|{\leq}|\, v$', the claim holds by 1.3, (VW1), 1.5, 1.7, 1.9, 1.14 and 1.8, resp.. The cases where $\varphi$ is of the form $\neg\psi$, $(\varphi_1 \,\&\, \varphi_2)$ or $(\varphi_1 \vee \varphi_2)$ follow immediately from the induction hypothesis. Suppose $\varphi$ is of the form '$\exists x \sqsubseteq t\, \psi(x)$' and assume, for (a), that $\mathfrak{S}^+ \vDash \exists x \sqsubseteq t\, \psi(x)$. Then $\mathfrak{S}^+ \vDash \exists x \sqsubseteq \underline{t^*}\, \psi(x)$ where $t^* = val_{\mathfrak{S}+}\, t$, so $\mathfrak{S}^+ \vDash \bigvee_{s \sqsubseteq t^*} \psi(\underline{s})$. From the induction hypothesis $VW^+ \vdash \bigvee_{s \sqsubseteq t^*} \psi(\underline{s})$, whence from 1.12(b) we obtain $VW^+ \vdash \exists x \sqsubseteq \underline{t^*}\, \psi(x)$. But then $VW^+ \vdash \exists x \sqsubseteq t\, \psi(x)$ follows by 1.2. For (b) we argue analogously using 1.12(a). The reasoning proceeds similarly if $\varphi$ is of the form '$\forall x \sqsubseteq t\, \psi(x)$', '$\exists x \sqsupseteq t\, \psi(x)$' or '$\forall x \sqsupseteq t\, \psi(x)$'. If $\varphi$ is of the form '$\forall x \subseteq_p t\, \psi(x)$' and '$\exists x \subseteq_p t\, \psi(x)$' we use 1.15 instead. Suppose $\varphi$ is of the form '$\exists x \,|{\leq}|\, t\, \psi(x)$'. If $t$ is not 2-tally, we apply 1.16(d). Assume $\mathfrak{S}^+ \vDash \exists x\ |\leq|\ t\ \psi(x)$ where t is a 2-tally $2^m$. Then $\mathfrak{S}^+ \vDash \psi(\underline{s})$ for some $s \in \Sigma^*$ such that $lh(s) \leq m$, and by the induction hypothesis we have that $VW^+ \vdash \psi(\underline{s})$. By 1.16(b) we have that $VW^+ \vdash \exists x\ (x\ |\leq|\ \underline{2}^m\ \&\ \psi(x))$. But then $VW^+ \vdash \exists x\ |\leq|\ t\ \psi(x)$ by 1.2, as needed. If $\varphi$ is of the form '$\forall x\ |\leq|\ t\ \psi(x)$' we use 1.16(c) and 1.16(a) instead, resp.. □

We now define the class of <u>$\Sigma$-formulae</u> of $\mathcal{L}^+$ inductively as follows:

(i) Any bounded $\mathcal{L}^+$ formula is a $\Sigma$-formula,

(ii) if $\varphi$ and $\psi$ are $\Sigma$-formulae, so are $(\varphi\ \&\ \psi)$ and $(\varphi \vee \psi)$,

(iii) if $\varphi$ is a $\Sigma$-formula, then so is $\exists x\ \varphi$,

(iv) if $\varphi$ is a $\Sigma$-formula and t an $\mathcal{L}^+$ term not containing the variable x, then $\exists x \sqsubseteq t\ \varphi$, $\forall x \sqsubseteq t\ \varphi$, $\exists x \sqsupseteq t\ \varphi$, $\forall x \sqsupseteq t\ \varphi$, $\exists x \subseteq_p t\ \varphi$, $\forall x \subseteq_p t\ \varphi$, $\exists x\ |\leq|\ t\ \varphi$ and $\forall x\ |\leq|\ t\ \varphi$ are $\Sigma$-formulae.

1.18 Let $\varphi$ be any $\Sigma$-sentence of $\mathcal{L}^+$: if $\mathfrak{S}^+ \vDash \varphi$, then $VW^+ \vdash \varphi$.

<u>Proof</u>: We argue by induction on the complexity of $\Sigma$-sentences using 1.17(a). We only consider the case where $\varphi \equiv \exists x\ \psi(x)$. Assume $\mathfrak{S}^+ \vDash \exists x\ \psi(x)$. Then $\mathfrak{S}^+ \vDash \psi(\underline{t})$ for some $t \in \Sigma^*$. Then $VW^+ \vdash \psi(\underline{t})$ holds by the induction hypothesis and $VW^+ \vdash \exists x\ \psi(x)$. □

1.19 Let t be any variable-free $\mathcal{L}^+$ term. Then

(a) $VW^+ \vdash \forall x \sqsubseteq t\ \forall y \sqsubseteq x\ \ y \sqsubseteq t$.

(b) $VW^+ \vdash \forall x \subseteq_p t\ \forall y \subseteq_p x\ \ y \subseteq_p t$.

This follows immediately from 1.17(a).

## §2. Dyadic Trees as Dyadic Strings

In [3], Kristiansen and Murwanashyaka introduced a first-order theory WT formulated in the vocabulary $\mathcal{L}_T = \{0, (\ ), \sqsubseteq\}$ with a single individual constant 0, a binary operation symbol ( , ) and a 2-place relational symbol $\sqsubseteq$, with -- along with the usual axioms for identity -- infinitely many axioms given by the instances of the two schemas

(WT1) $\neg s = t$ for any distinct variable-free terms s, t of $\mathcal{L}_T$,

(WT2) $\forall x\ (x \sqsubseteq t \leftrightarrow \bigvee_{s \in \mathcal{S}[t]} x = s)$ for each variable-free term t of $\mathcal{L}_T$,

where $\mathcal{S}[t]$ is the set of all subterms of t. The intended interpretation $\mathcal{T} = \langle \mathbf{T}, \underline{0}, \tau, \sqsubseteq_T \rangle$ of the theory WT is the term model where **T** is the set of all variable-free $\mathcal{L}_T$-terms

0, (00), (0(00)), ((00)0), ((00)(00)),…

(We often omit the outermost parentheses.) The elements of **T** may be identified with finite dyadic trees, trees in which every non-terminal node has two immediate descendants. With that in mind we may take the individual constant 0 in $\mathcal{T}$ to stand for the single-node tree $\underline{0}$, the binary operation $\tau : \mathbf{T} \times \mathbf{T} \rightarrow \mathbf{T}$ applied to trees $T_1, T_2$ yields the tree $\tau(T_1, T_2)$ whose root node has as its immediate descendants the root nodes of $T_1$ and $T_2$, and the binary relational symbol $\sqsubseteq$ is meant to express the subtree relation between trees, where subtrees are defined so that for a given tree $T$, any of its nodes x determines a subtree $T_x$ consisting of all and only the descendants of x in $T$ including itself. Kristiansen and Murwanashyaka showed that WT is mutually interpretable with Robinson's very weak arithmetic R.

For the purposes of formally interpreting WT in $VW^+$ we reformulate WT in the expanded relational vocabulary $\mathcal{L}_{T(Rel,c)} = \{b_0, b_1, b_2,..., T, \sqsubseteq\}$ with infinitely many individual constants $b_0, b_1, b_2,...$ and a ternary relational symbol T in place of the binary operation symbol '( , )'. The constants $b_0, b_1, b_2,...$ are intended to denote the distinct elements in **T** – the variable-free terms of $\mathcal{L}_T$ or dyadic trees – and the relational symbol T the graph of the term/tree-forming operation τ. It will be convenient to use the notation '$\underline{r}$', '$\underline{s}$', '$\underline{t}$' for the $\mathcal{L}_{T(Rel,c)}$ constants when 'r', 's', 't' are various variable-free terms of $\mathcal{L}_T$. We then restate the axioms of WT in the following form:

($WT_{Rel}1$) $\neg \underline{s} = \underline{t}$ for any distinct variable-free terms s, t of $\mathcal{L}_T$,

($WT_{Rel}2$) $\forall x\, (x \sqsubseteq \underline{t} \leftrightarrow \bigvee_{s \in \mathcal{S}[t]} x = \underline{s})$ for each variable-free term t of $\mathcal{L}_T$,

($WT_{Rel}3$) $T(\underline{s},\underline{t},\underline{(s,t)})$ for any variable-free terms s, t of $\mathcal{L}_T$,

($WT_{Rel}4$) $\forall x, y\, \exists z\, T(x,y,z)$

($WT_{Rel}5$) $\forall x, y, u, v\; (T(x,y,u)\ \&\ T(x,y,v) \rightarrow u = v)$,

along with appropriately formulated axioms for identity.

We now proceed to construct elements of the desired interpretation of WT in $VW^+$. For this purpose we represent variable-free terms of $\mathcal{L}_T$ by certain dyadic strings: e.g., the terms 0, (00), (00)0, 0(00), (00)(00), 0((00)0) shall be represented by 1, 121, $1212^21$, $12^2121$, $1212^2121$, $12^31212^21$, resp., writing '$1^n$' and '$2^n$' as shorthand for tallies of length n, i.e strings of n consecutive 1's or 2's, resp., and omitting, for transparency, symbols for the concatenation operation. We say that a 2-tally $2^n$, $n \geq 1$, is a <u>separator</u> in a dyadic string σ if the string $1 ^\frown 2^n {}^\frown 1$ occurs as a segment of σ.

Let $\Theta \subseteq \Sigma^*$ be the smallest set satisfying the following conditions: (i) $1 \in \Theta$,

(ii) if $\sigma_1, \sigma_2 \in \Theta$, and $2^k$ is the shortest 2-tally exceeding in length any 2-tally occurring in either $\sigma_1$ or $\sigma_2$, then $\sigma_1\hat{}\, 2^k\hat{}\, \sigma_2 \in \Theta$.

(We call such a 2-tally the <u>minimax separator</u> of $\sigma_1\hat{}\, 2^k\hat{}\, \sigma_2$.) It is clear that each string $\sigma$ in $\Theta$ has a unique decomposition $\sigma_1, \sigma_2, \ldots, \sigma_n$, where the resulting substrings $\sigma_i$, $1 \leq i \leq n$, are all in $\Theta$ and represent the subtrees of a dyadic tree $T_\sigma$. We call the string $\sigma$ <u>the (unique) minimax code</u> of $T_\sigma$.

We inductively define a mapping $\theta$ of variable-free $\mathcal{L}_T$ terms to $\mathcal{L}^+$ concatenation terms that are their minimax codes according to rules (i) and (ii), as follows: $\theta(0) =: \underline{1}$ and $\theta((u,v)) =: (\theta(u)\hat{}\,\underline{2}^k)\hat{}\,\theta(v)$. Then $\theta(u) \in \Theta$ for each variable-free $\mathcal{L}_T$ term u.

We proceed to give an explicit characterization of the strings in $\Theta$. Consider the following conditions:

($\theta$i) x begins with $\underline{1}$ and ends with $\underline{1}$,

($\theta$ii) x contains only single occurrences of $\underline{1}$,

($\theta$iii) x contains a single occurrence of a separator, $s_x$, if any, exceeding in length any other separator in x,

($\theta$iv) if any separators occur in x, then every initial segment of $s_x$ occurs as a separator in x.

It will be convenient to regard the empty string e as the minimax separator of $\underline{1}$, writing $\underline{2}^0$ for $\underline{e}$. Conditions ($\theta$iii) and ($\theta$iv) identify $s_x$ as the shortest 2-tally that functions as a separator in x all of whose initial segments are also separators in x and which is also longer than any other separator in x.

2.1 Let $\sigma \in \Sigma^*$. Then: $\sigma$ satisfies ($\theta$i)-($\theta$iv) $\Leftrightarrow$ $\sigma \in \Theta$.

Proof: Suppose $\sigma \in \Theta$. We argue by induction on the number of applications of rule (ii) in $\sigma$. If $\sigma$ is 1, conditions ($\theta$i) and ($\theta$ii) hold and ($\theta$iii) and ($\theta$iv) hold trivially. Suppose $\sigma = \sigma_1\hat{}\, 2^n\hat{}\, \sigma_2$ where $2^n$ is the minimax separator of $\sigma$ and $\sigma_1, \sigma_2 \in \Theta$. By the induction hypothesis both $\sigma_1, \sigma_2$ satisfy ($\theta$i)-($\theta$iv). Then ($\theta$i) and ($\theta$ii) hold for $\sigma$ by choice of $\sigma$. Let $2^j$ and $2^k$ be the minimax separators of $\sigma_1$ and $\sigma_2$, resp.. By choice of n we have that $n = \max(j,k) + 1$. Then ($\theta$iii) is immediate. Since the initial segments of $2^n$ are those of $2^j$ and $2^k$ plus $2^n$ itself, ($\theta$iv) follows as well. This completes the inductive argument that every string in $\Theta$ satisfies ($\theta$i)-($\theta$iv).

Conversely, assume $\sigma$ satisfies ($\theta$i)-($\theta$iv). We argue by induction on the number of blocks of 2's in $\sigma$. Let $2^n$ be the separator occurring in $\sigma$ exceeding in length any other separator in $\sigma$, which exists by ($\theta$iii). Then $1\hat{}\, 2^n\hat{}\, 1$ is a segment of $\sigma$, and so there are nonempty segments $\sigma_1$ and $\sigma_2$ of $\sigma$ such that $\sigma = \sigma_1\hat{}\, 2^n\hat{}\, \sigma_2$. Since $\sigma_1$ and $\sigma_2$ each have fewer blocks of 2's than $\sigma$, by the induction hypothesis both $\sigma_1$ and $\sigma_2$ are in $\Theta$. By the first part of the proof above $\sigma_1$ and $\sigma_2$ both have properties ($\theta$i)-($\theta$iv). Let $2^j$ and $2^k$ be their minimax separators. Then, since by hypothesis $\sigma$ satisfies ($\theta$iii) we have that $n \geq \max(j,k) + 1$. And since $\sigma$ also satisfies ($\theta$iv), we have that $n \leq \max(j,k) + 1$. Hence $2^n$ is the shortest 2-tally exceeding in length both $2^j$ and $2^k$. But no 2-tally occurring in either $\sigma_1$ or $\sigma_2$ can exceed in length both $2^j$ and $2^k$. Hence $2^n$ is the shortest 2-tally exceeding in length any 2-tally occurring in either $\sigma_1$ or $\sigma_2$, that is, $2^n$ is the minimax separator of $\sigma$, and so $\sigma \in \Theta$, as claimed. □

We now define the relevant concepts in $\mathcal{L}^+$. Let

$$\mathrm{Tally}_1(x) \;\equiv:\; x \neq \underline{e} \;\&\; \neg\underline{2} \subseteq_p x \quad\text{and}\quad \mathrm{Tally}_2(x) \;\equiv:\; x \neq \underline{e} \;\&\; \neg\underline{1} \subseteq_p x,$$

$$\mathrm{Sep}(x,y) \;\equiv:\; (\mathrm{Tally}_2(x) \;\&\; (\underline{1}_{\vee}x){}^{\wedge}\underline{1} \subseteq_p y) \vee (y = \underline{1} \;\&\; x = \underline{e}),$$

$$\mathrm{MaxSep}(x,y) \;\equiv:\; x \subseteq_p y \;\&\; \mathrm{Sep}(x,y) \;\&\; \forall z \subseteq_p y\, (\mathrm{Sep}(z,y) \;\&\; z \neq x \;\rightarrow\; z{}^{\wedge}\underline{2} \sqsubseteq x) \;\&$$

$$\&\; \forall z \subseteq_p y\, (\mathrm{Tally}_2(z) \;\rightarrow\; z \sqsubseteq x)$$

$$\mathrm{MinMaxSep}(x,y) \;\equiv:\; \mathrm{MaxSep}(x,y) \;\&\; \forall z \sqsubseteq x\, \exists u \subseteq_p y\, \mathrm{MaxSep}(z,u))$$

We then have:

2.2 $\quad \mathrm{VW}^+ \vDash \forall u,v,x\, (\mathrm{MinMaxSep}(u,x) \;\&\; \mathrm{MinMaxSep}(v,x) \;\rightarrow\; u = v)$

Proof: Assume $M \vDash \mathrm{MinMaxSep}(u,x) \;\&\; \mathrm{MinMaxSep}(v,x)$ where $u, v, x \in M$. We may suppose that $x \neq \underline{1}$. From the hypothesis we have

$$M \vDash \mathrm{Tally}_2(u) \;\&\; \mathrm{Tally}_2(v) \;\&\; u \subseteq_p x \;\&\; v \subseteq_p x\,.$$

Assume $\mathfrak{S}^+ \vDash u \sqsubseteq v \;\&\; u \neq v$. Then from $M \vDash \mathrm{MaxSep}(u,x) \;\&\; \mathrm{MaxSep}(v,x)$ it follows that $M \vDash u \sqsubseteq v \;\&\; v \sqsubseteq u$., whence $M \vDash u = v$ by (VW10). □

Let $T^*(x,y,w)$ abbreviate the $\mathcal{L}^+$ formula

$$\underline{1} \sqsubseteq w \;\&\; \underline{1} \sqsupseteq w \;\&\; \neg\underline{1}{}^{\wedge}\underline{1} \subseteq_p w \;\&$$

$$\&\; \exists t_1 \subseteq_p x\, \exists t_2 \subseteq_p y\, (\mathrm{MinMaxSep}(t_1,x) \;\&\; t_1 \subseteq_1 x \;\&\; \mathrm{MinMaxSep}(t_2,y) \;\&\; t_2 \subseteq_1 y \;\&$$

$$\&\; \exists t \subseteq_p w(\mathrm{MinMaxSep}(t,w) \;\&\; t \subseteq_1 w \;\&$$

$\& ((t_1 \sqsubseteq t_2 \ \& \ t = t_2{}^\wedge\underline{2}) \ v \ (t_2 \sqsubseteq t_1 \ \& \ t = t_1{}^\wedge\underline{2})) \ \& \ x \sqsubseteq w \ \& \ y \sqsupseteq w \ \&$

$\& \ \exists! w_1 \sqsubseteq w \ \exists! w_2 \sqsupseteq w \ (\underline{1}{}^\wedge t \sqsupseteq w_1 \ \& \ t^\wedge\underline{1} \sqsubseteq w_2 \ \& \ x \sqsubseteq^+ w_1 \ \& \ y \sqsupseteq^+ w_2 \ \&$

$\& \ \forall u \sqsubseteq x \ u \sqsubseteq w \ \& \ \forall u \sqsubseteq w_1 \ \forall u_1 \sqsubseteq u \ u_1 \sqsubseteq w \ \&$

$\& \ \forall z \sqsubseteq w_1 \ (z \sqsubseteq^+ w_1 \ \rightarrow \ z^\wedge\underline{2} \sqsubseteq w_1 \ v \ (\neg(z^\wedge\underline{2} \sqsubseteq w_1) \ \& \ 1^\wedge t \sqsupseteq z)) \ \&$

$\& \ \forall z \sqsupseteq w_2 \ (z \sqsupseteq^+ w_2 \ \rightarrow \ \underline{2}_\vee z \sqsupseteq w_2 \ v \ (\neg(\underline{2}_\vee z \sqsupseteq w_2) \ \& \ t^\wedge 1 \sqsubseteq z)) \ \&$

$\& \ \forall u \sqsubseteq w_2 \ (t \subseteq_1 u \rightarrow \exists! v \sqsubseteq w \ (x \sqsubseteq v \ \& \ t \subseteq_1 v \ \& \ u \sqsupseteq v)))))$

where $x \sqsubseteq^+ u \ \equiv: \ x \sqsubseteq u \ \& \ \forall x_1 \sqsubseteq u \ (\underline{1} \sqsupseteq x_1 \ \rightarrow \ x_1 \sqsubseteq x)$

and $x \sqsupseteq^+ u \ \equiv: \ x \sqsupseteq u \ \& \ \forall x_2 \sqsupseteq u \ (1 \sqsubseteq x_2 \ \& \ x_2 \sqsupseteq u \rightarrow \ x_2 \sqsupseteq x)$, and

$\exists! x \sqsubseteq t \ A(x) \equiv: \exists x \sqsubseteq t \ (A(x) \ \& \ \forall y \sqsubseteq t \ A(y) \rightarrow y = x)$, and analogously for $\sqsupseteq$.
(Here it is assumed that x does not occur in t.)

Let $T^+(x,y,z)$ abbreviate the $\mathcal{L}^+$ formula

$\exists w \ (T^*(x,y,w) \ \& \ \forall v(T^*(x,y,w) \rightarrow w \sqsubseteq v) \ \& \ z = w) \ v$

$v \ (\neg\exists w \ (T^*(x,y,w) \ \& \ \forall v \ (T^*(x,y,w) \rightarrow w \sqsubseteq v)) \ \& \ z = \underline{1})$,

and let $I(x) \ \equiv: \ \underline{1} \sqsubseteq x \ \& \ \underline{1} \sqsupseteq x \ \& \ \neg\underline{1}{}^\wedge\underline{1} \subseteq_p x$.

We note the following obvious fact:

2.3 Let $s, t \in \Sigma^*$. Then

$\mathfrak{S}^+ \ \vDash \forall x \sqsubseteq \underline{s} \ \forall y \sqsupseteq \underline{s} \ (\text{Tally}_2(\underline{t}) \ \& \ \underline{t} \subseteq_1 \underline{s} \ \& \ x^\wedge\underline{t} \sqsubseteq \underline{s} \ \& \ \underline{t}_\vee y \sqsupseteq \underline{s} \ \rightarrow$

$\rightarrow \ \forall z \subseteq_p \underline{s} \ (\text{Tally}_1(z) \ \rightarrow \ z \subseteq_p x \ v \ z \subseteq_p y))$.

2.4 For each variable-free $\mathcal{L}_T$ term t there is an $n \in \mathbb{N}$ such that

$$VW^+ \vdash I(\underline{\theta(t)}) \ \&\ MinMaxSep(\underline{2}^n, \underline{\theta(t)}) \ \&\ (\underline{\theta(t)} \neq \underline{1} \rightarrow \underline{2}^n \subseteq_1 \underline{\theta(t)}) \ .$$

Proof: We argue by induction on the complexity of $\mathcal{L}_T$ terms. If t is 0, we have that $\theta(t)$ is $\underline{1}$. Now $\mathfrak{S}^+ \vDash I(\underline{1})\ \&\ MinMaxSep(\underline{2}^0,\underline{1})$, and we have $VW^+ \vdash I(\underline{1})\ \&\ MinMaxSep(\underline{2}^0,\underline{1})$ by 1.17(a). Suppose t is (u,v) and assume, as the induction hypothesis, that

$$VW^+ \vdash I(\underline{\theta(u)}) \ \&\ MinMaxSep(\underline{2}^j, \underline{\theta(u)}) \ \&\ (\underline{\theta(u)} \neq \underline{1} \rightarrow \underline{2}^j \subseteq_1 \underline{\theta(u)}) \ \&$$

$$\&\ I(\underline{\theta(v)}) \ \&\ MinMaxSep(\underline{2}^k, \underline{\theta(v)}) \ \&\ (\underline{\theta(v)} \neq \underline{1} \rightarrow \underline{2}^k \subseteq_1 \underline{\theta(v)}) \ .$$

Now, $\theta((u,v))$ is $\theta(u)\,\hat{}\,\underline{2}^n\,\hat{}\,\theta(v)$ where $2^n$ is the shortest 2-tally exceeding in length any 2-tally in either $\theta(u)$ or $\theta(v)$. By 2.1 we have that $n = \max(j,k) + 1$.

From $VW^+ \vdash I(\underline{\theta(u)})\ \&\ I(\underline{\theta(v)})$ we have that $\mathfrak{S}^+ \vDash \underline{1} \sqsubseteq \theta(u)\ \&\ \underline{1} \sqsupseteq \theta(v)$, hence $\mathfrak{S}^+ \vDash \underline{1} \sqsubseteq \theta((u,v))\ \&\ \underline{1} \sqsupseteq \theta((u,v))$. Suppose $\mathfrak{S}^+ \vDash \underline{1}\hat{}\underline{1} \subseteq_p \underline{\theta((u,v))}$. Then, since

$$\mathfrak{S}^+ \vDash Tally_2(\underline{2}^n)\ \&\ \underline{\theta(u)}\hat{}\underline{2}^n \sqsubseteq \underline{\theta((u,v))}\ \&\ \underline{2}^n{}_\vee\underline{\theta(v)} \sqsupseteq \underline{\theta((u,v))}\ \&\ \underline{2}^n \subseteq_1 \underline{\theta((u,v))}$$

we have, by 2.3, that $\mathfrak{S}^+ \vDash \underline{1}\hat{}\underline{1} \subseteq_p \underline{\theta(u)}\ \vee\ \underline{1}\hat{}\underline{1} \subseteq_p \underline{\theta(v)}$. But this contradicts $VW^+ \vdash I(\underline{\theta(u)})\ \&\ I(\underline{\theta(v)})$. Hence $\mathfrak{S}^+ \vDash \neg\underline{1}\hat{}\underline{1} \subseteq_p \theta\underline{((u,v))}$. By 1.17(a) it follows that $VW^+ \vdash \underline{1} \sqsubseteq \underline{\theta((u,v))}\ \&\ \underline{1} \sqsupseteq \underline{\theta((u,v))}\ \&\ \neg\underline{1}\hat{}\underline{1} \subseteq_p \underline{\theta((u,v))}$, hence $VW^+ \vdash I(\underline{\theta((u,v))})$.

Now, by choice of n and definition of $\theta((u,v))$ we have that $\mathfrak{S}^+ \vDash \underline{2}^n \subseteq_1 \underline{\theta(t)}$, hence by 1.17(a) it follows that $VW^+ \vdash (\underline{\theta(t)} \neq \underline{1} \rightarrow \underline{2}^n \subseteq_1 \underline{\theta(t)})$ .

It remains to prove that $VW^+ \vdash MinMaxSep(\underline{2}^n,\underline{\theta((u,v))})$. We have that

$\mathfrak{S}^+ \vDash Tally_2(\underline{2}^n)\ \&\ (\underline{1}_\vee\underline{2}^n)\hat{}\underline{1} \subseteq_p \underline{\theta((u,v))}$, so $\mathfrak{S}^+ \vDash Sep(\underline{2}^n,\underline{\theta((u,v))})$. Suppose $\mathfrak{S}^+ \vDash z \subseteq_p \underline{\theta((u,v))}\ \&\ Sep(z,\underline{\theta((u,v))})\ \&\ z \neq \underline{2}^n$. Then $\mathfrak{S}^+ \vDash Tally_2(z)$ and by choice of $\theta((u,v))$ we have that $\underline{2}^n$ exceeds z in length.

We first argue that $\mathfrak{S}^+ \models \mathrm{Sep}(z,\underline{\theta(u)})$ or $\mathfrak{S}^+ \models \mathrm{Sep}(z,\underline{\theta(v)})$. By hypothesis

$\mathfrak{S}^+ \models (\underline{1}_{\vee}z)^\wedge\underline{1} \subseteq_p \underline{\theta((u,v))}$, so $\mathfrak{S}^+ \models \exists w \sqsubseteq \underline{\theta((u,v))}\ (\underline{1}_{\vee}z)^\wedge\underline{1} \sqsupseteq w$. Now we also have that $\mathfrak{S}^+ \models \underline{\theta(u)} \sqsubseteq \underline{\theta((u,v))}$, hence

$$\mathfrak{S}^+ \models \underline{\theta(u)} \sqsubseteq w \ \vee \ w = \underline{\theta(u)} \ \vee \ w \sqsubseteq \underline{\theta(u)}.$$

If $\mathfrak{S}^+ \models w = \underline{\theta(u)} \vee w \sqsubseteq \underline{\theta(u)}$, then $\mathfrak{S}^+ \models (\underline{1}_{\vee}z)^\wedge\underline{1} \subseteq_p \underline{\theta(u)}$, so
$\mathfrak{S}^+ \models \mathrm{Sep}(z,\underline{\theta(u)})$. Suppose, on the other hand, that $\mathfrak{S}^+ \models \underline{\theta(u)} \sqsubseteq w \;\&\; \underline{\theta(u)} \neq w$. Then $\mathfrak{S}^+ \models \exists w_1\, (\underline{\theta(u)}_{\vee}w_1 = w \;\&\; w_1 \neq \underline{e})$. Since

$$\mathfrak{S}^+ \models (\underline{\theta(u)}_{\vee}w_1 = w \sqsubseteq \underline{\theta((u,v))} = (\underline{\theta(u)}^\wedge\underline{2}^n)^\wedge\underline{\theta(v)},$$

we have that $\mathfrak{S}^+ \models w_1 \sqsubseteq \ \underline{2}^{n\wedge}\underline{\theta(v)}$, so $\mathfrak{S}^+ \models \underline{2} \sqsubseteq w_1$. Given that
$\mathfrak{S}^+ \models (\underline{1}_{\vee}z)^\wedge\underline{1} \sqsupseteq w \;\&\; w_1 \sqsupseteq w$, we have

$$\mathfrak{S}^+ \models (\underline{1}_{\vee}z)^\wedge\underline{1} \sqsupseteq w_1 \ \vee \ w_1 = (\underline{1}_{\vee}z)^\wedge\underline{1} \ \vee \ w_1 \sqsupseteq (\underline{1}_{\vee}z)^\wedge\underline{1}.$$

Now, $\mathfrak{S}^+ \models w_1 = (\underline{1}_{\vee}z)^\wedge\underline{1}$ is ruled out because $\mathfrak{S}^+ \models \underline{2} \sqsubseteq w_1$. Suppose
$\mathfrak{S}^+ \models w_1 \sqsupseteq (\underline{1}_{\vee}z)^\wedge\underline{1}$. Then, since $\mathfrak{S}^+ \models w_1 \neq \underline{e} \;\&\; w_1 \sqsubseteq \underline{2}^{n\wedge}\underline{\theta(v)}$, we must have $\mathfrak{S}^+ \models \underline{2}^{n\wedge}\underline{1} \sqsubseteq w_1$. But then $\mathfrak{S}^+ \models \underline{2}^n \subseteq_p (\underline{1}_{\vee}z)^\wedge\underline{1}$, which is impossible.
Therefore, $\mathfrak{S}^+ \models \neg w_1 \sqsupseteq (\underline{1}_{\vee}z)^\wedge\underline{1}$ and we must have $\mathfrak{S}^+ \models (\underline{1}_{\vee}z)^\wedge\underline{1} \sqsupseteq w_1$. But then from $\mathfrak{S}^+ \models w_1 \sqsubseteq \ \underline{2}^{n\wedge}\underline{\theta(v)}$, it follows that $\mathfrak{S}^+ \models (\underline{1}_{\vee}z)^\wedge\underline{1} \subseteq_p \underline{\theta(v)}$. That is,
$\mathfrak{S}^+ \models \mathrm{Sep}(z,\underline{\theta(v)})$.

Hence, $\mathfrak{S}^+ \models \mathrm{Sep}(z,\underline{\theta(u)}) \ \vee \ \mathrm{Sep}(z,\underline{\theta(v)})$, and by choice of n we have that
$\mathfrak{S}^+ \models z^\wedge\underline{2} \sqsubseteq \underline{2}^n$. Therefore,

$$\mathfrak{S}^+ \models \forall z \subseteq_p \underline{\theta((u,v))}\ (\mathrm{Sep}(z,\underline{\theta((u,v)})\ \&\ z \neq \underline{2}^n \ \rightarrow \ z^\wedge\underline{2} \sqsubseteq \underline{2}^n).$$

Also by choice of n, we have that $\mathfrak{S}^+ \models \forall z \subseteq_p \underline{\theta((u,v))}\ (\mathrm{Tally}_2(z) \ \rightarrow \ z \sqsubseteq \underline{2}^n)$.

Thus, $\mathfrak{S}^+ \models \mathrm{MaxSep}(\underline{2}^n,\underline{\theta((u,v)})$.

Finally, suppose $\mathfrak{S}^+ \models z \sqsubseteq \underline{2}^n$. By choice of n, then $\mathfrak{S}^+ \models z \sqsubseteq \underline{2}^{\max(j,k)} \vee z = \underline{2}^n$. If $\mathfrak{S}^+ \models z = \underline{2}^n$, we already have that $\mathfrak{S}^+ \models \mathrm{MaxSep}(\underline{2}^n, \underline{\theta((u,v))})$. If $\mathfrak{S}^+ \models z \sqsubseteq \underline{2}^{\max(j,k)}$, suppose $j = \max(j,k)$ and $\mathfrak{S}^+ \models z \sqsubseteq \underline{2}^j$. Then from the induction hypothesis $\mathrm{VW}^+ \vdash \mathrm{MinMaxSep}(\underline{2}^j, \underline{\theta(u)})$, it follows that $\mathfrak{S}^+ \models \exists u_1 \subseteq_p \underline{\theta(u)}\ \mathrm{MaxSep}(z,u_1)$. But $\mathfrak{S}^+ \models \underline{\theta(u)} \subseteq_p \underline{\theta((u,v))}$, hence this suffices to establish that $\mathfrak{S}^+ \models \exists u_1 \subseteq_p \underline{\theta((u,v))}\ \mathrm{MaxSep}(z,u_1)$. Similarly if $k = \max(j,k)$. Therefore we have also established that

$$\mathfrak{S}^+ \models \forall z \sqsubseteq \underline{2}^n\ \exists u_1 \subseteq_p \underline{\theta((u,v))}\ \mathrm{MaxSep}(z,u_1).$$

From 1.17(a) it follows that $\mathrm{VW}^+ \vdash \mathrm{MinMaxSep}(\underline{2}^n, \underline{\theta((u,v))})$. □

2.5 $\mathrm{VW}^+ \vdash \forall x,y\ (I(x)\ \&\ I(y) \rightarrow \exists z\ (I(z)\ \&\ T^+(x,y,z)))$.

Proof: Assume $M \models I(x)\ \&\ I(y)$ where $x, y \in M$. Suppose $M \models \exists w\ (T^*(x,y,w)\ \&\ \forall v\ (T^*(x,y,v) \rightarrow w \sqsubseteq v))$. Let $z = w$. Then $M \models \underline{1} \sqsubseteq z\ \&\ \underline{1} \sqsupseteq z\ \&\ \neg \underline{1}{}^\wedge\underline{1} \subseteq_p z$ follows from $M \models T^*(x,y,w)$ and so we have $M \models I(z)\ \&\ T^+(x,y,z)$. On the other hand, suppose

$$M \models \neg\exists w\ (T^*(x,y,w)\ \&\ \forall v\ (T^*(x,y,v) \rightarrow w \sqsubseteq v)).$$

Let $z = 1$. We have that $\mathfrak{S}^+ \models \underline{1} \sqsubseteq \underline{1}\ \&\ \underline{1} \sqsupseteq \underline{1}\ \&\ \neg\underline{1}{}^\wedge\underline{1} \subseteq_p \underline{1}$, so by 1.5, 1.9(a) and 1.14(b) it follows that $M \models \underline{1} \sqsubseteq \underline{1}\ \&\ \underline{1} \sqsupseteq \underline{1}\ \&\ \neg\underline{1}{}^\wedge\underline{1} \subseteq_p \underline{1}$, that is, $M \models I(z)\ \&\ T^+(x,y,z)$, as needed. □

2.6 $\mathrm{VW}^+ \vdash \forall x,y,u,v\ (I(x)\ \&\ I(y)\ \&\ I(u)\ \&\ I(v) \rightarrow (T^+(x,y,u)\ \&\ T^+(x,y,v) \rightarrow$

$\rightarrow\ u = v)).$

This is immediate from the definition of $T^+$ and (VW10).

2.7 Let s, t be variable-free $\mathcal{L}^+$-concatenation terms. Then

$$VW^+ \vdash \forall w\,(\forall z \sqsubseteq w\,(z \sqsubseteq^+ w \rightarrow z^\wedge\underline{2} \sqsubseteq w \ \mathrm{v}\ (\neg(z^\wedge\underline{2} \sqsubseteq w)\ \&\ 1^\wedge t \sqsupseteq z))\ \&$$

$$\&\ s \sqsubseteq^+ w\ \&\ 1^\wedge t \sqsupseteq w\ \&\ \underline{1} \sqsupseteq s\ \&\ \mathrm{Tally}_2(t) \rightarrow s^\wedge t \sqsubseteq w).$$

Proof: Assume $M \vDash \forall z \sqsubseteq w\,(z \sqsubseteq^+ w \rightarrow z^\wedge\underline{2} \sqsubseteq w \ \mathrm{v}\ (\neg(z^\wedge\underline{2} \sqsubseteq w)\ \&\ 1^\wedge t \sqsupseteq z))$ along with $M \vDash s \sqsubseteq^+ w\ \&\ 1^\wedge t \sqsupseteq w\ \&\ \underline{1} \sqsupseteq s\ \&\ \mathrm{Tally}_2(t)$ for $w \in M$ and variable-free concatenation terms s, t.

From 1.17 we have that $\mathfrak{S}^+ \vDash \mathrm{Tally}_2(t)$. We first consider the case $t = 2$. Then from $M \vDash s \sqsubseteq^+ w$ we have $M \vDash s \sqsubseteq w$ and further

$$M \vDash s^\wedge\underline{2} \sqsubseteq w \ \mathrm{v}\ (\neg(s^\wedge\underline{2} \sqsubseteq w)\ \&\ \underline{1}^\wedge\underline{2} \sqsupseteq s)$$

from the hypothesis. Suppose, for a reductio, that $M \vDash \underline{1}^\wedge\underline{2} \sqsupseteq s$. Then we have that $M \vDash \underline{1} \sqsupseteq s\ \&\ \underline{1}^\wedge\underline{2} \sqsupseteq s$, whence by 1.9 it follows that $\mathfrak{S}^+ \vDash \underline{1} \sqsupseteq s\ \&\ \underline{1}^\wedge\underline{2} \sqsupseteq s$, which is impossible. Therefore $M \vDash \neg(\underline{1}^\wedge\underline{2} \sqsupseteq s)$, and so $M \vDash s^\wedge\underline{2} \sqsubseteq w$, as needed.

Suppose that t is $2^n$ for some $n > 1$. We now argue, under the above principal hypothesis, that

(*) $\quad M \vDash s^\wedge\underline{2}^k \sqsubseteq^+ w \rightarrow s^\wedge\underline{2}^{k+1} \sqsubseteq^+ w$

for any k, $0 \leq k < n$. Assume that $M \vDash s^\wedge\underline{2}^k \sqsubseteq^+ w$. Then

$$M \vDash s^\wedge\underline{2}^k \sqsubseteq w\ \&\ \forall x_1 \sqsubseteq w\,(\underline{1} \sqsupseteq x_1 \rightarrow x_1 \sqsubseteq s^\wedge\underline{2}^k).$$

Assuming further $M \models \underline{1} \sqsupseteq x_1$ for $M \models x_1 \sqsubseteq w$, we have that $M \models x_1 \sqsubseteq s{}^\wedge\underline{2}^k$. But

$M \models s{}^\wedge\underline{2}^k \sqsubseteq (s{}^\wedge\underline{2}^k){}^\wedge\underline{2}$ by 1.5, and $M \models (s{}^\wedge\underline{2}^k){}^\wedge\underline{2} = s{}^\wedge(\underline{2}^k{}^\wedge\underline{2}) = s{}^\wedge\underline{2}^{k+1}$ by 1.3, and further $M \models \forall y \sqsubseteq s{}^\wedge\underline{2}^k (s{}^\wedge\underline{2}^k \sqsubseteq s{}^\wedge\underline{2}^{k+1} \rightarrow y \sqsubseteq s{}^\wedge\underline{2}^{k+1})$ by 1.18(a). Hence $M \models x_1 \sqsubseteq s{}^\wedge\underline{2}^{k+1}$, and so we have

(#) $\quad M \models \forall x_1 \sqsubseteq w\ (\underline{1} \sqsupseteq x_1 \rightarrow x_1 \sqsubseteq s{}^\wedge\underline{2}^{k+1})$.

From hypothesis $M \models s{}^\wedge\underline{2}^k \sqsubseteq^+ w$ and the principal hypothesis we derive

$M \models (s{}^\wedge\underline{2}^k){}^\wedge\underline{2} \sqsubseteq w \ \vee\ (\neg((s{}^\wedge\underline{2}^k){}^\wedge\underline{2} \sqsubseteq w) \ \&\ \underline{1}{}^\wedge t \sqsupseteq s{}^\wedge\underline{2}^k)$. Suppose, for a reductio, that the latter disjunct holds. Then $M \models \underline{1}{}^\wedge t \sqsupseteq s{}^\wedge\underline{2}^k$, that is, $M \models \underline{1}{}^\wedge\underline{2}^n \sqsupseteq s{}^\wedge\underline{2}^k$, hence, by 1.9, also $\mathfrak{S}^+ \models \underline{1}{}^\wedge\underline{2}^n \sqsupseteq s{}^\wedge\underline{2}^k$. Then $s_0\,\hat{}\,1\,\hat{}\,2^n = s\,\hat{}\,2^k$ for some $s_0 \in \Sigma^*$, and since $n > k$, we have that $s_0\,\hat{}\,1\,\hat{}\,2^j = s$ where $n = k + j$. But this is impossible because $1 \sqsupseteq s$ follows from the principal hypothesis by 1.9. Therefore $M \models 1{}^\wedge t \sqsupseteq s{}^\wedge\underline{2}^k$ is ruled out, and we must have $M \models (s{}^\wedge\underline{2}^k){}^\wedge\underline{2} \sqsubseteq w$, that is, $M \models s{}^\wedge\underline{2}^{k+1} \sqsubseteq w$. But then along with (#) this shows that $M \models s{}^\wedge\underline{2}^{k+1} \sqsubseteq^+ w$, which completes the proof of (*).

Now, let $k + 1 = n$. By applying (*) n times we have from the hypothesis $M \models s \sqsubseteq^+ w$ that $M \models s{}^\wedge\underline{2}^n \sqsubseteq^+ w$, whence $M \models s{}^\wedge t \sqsubseteq w$. $\square$

2.8 Let r, t be variable-free $\mathcal{L}^+$-concatenation terms. Then

$$VW^+ \vdash \forall w\, (\forall z \sqsupseteq w\, (z \sqsupseteq^+ w \rightarrow \underline{2}_\vee z \sqsupseteq w \ \vee\ (\neg(\underline{2}_\vee z \sqsupseteq w) \ \&\ t{}^\wedge\underline{1} \sqsubseteq z)) \ \&$$

$$\&\ r \sqsupseteq^+ w \ \&\ t{}^\wedge 1 \sqsubseteq w \ \&\ \underline{1} \sqsubseteq r \ \&\ \mathrm{Tally}_2(t) \rightarrow t_\vee r \sqsupseteq w).$$

<u>Proof</u>: Assume $M \models \forall z \sqsupseteq w\, (z \sqsupseteq^+ w \rightarrow \underline{2}_\vee z \sqsupseteq w \ \vee\ (\neg(\underline{2}_\vee z \sqsupseteq w) \ \&\ t{}^\wedge\underline{1} \sqsubseteq z))$ along with $M \models r \sqsupseteq^+ w \ \&\ t{}^\wedge 1 \sqsubseteq w \ \&\ \underline{1} \sqsubseteq r \ \&\ \mathrm{Tally}_2(t)$ for $w \in M$ and variable-free concatenation terms r, t.

Consider the case $t = 2$. Again, we have from the hypothesis $M \vDash Tally_2(t)$ and 1.17 that $\mathfrak{S}^+ \vDash Tally_2(t)$. Then from $M \vDash r \sqsupseteq^+ w$ we have $M \vDash r \sqsupseteq w$ and further $M \vDash \underline{2}_{\vee} r \sqsupseteq w \ \text{v}\ (\neg(\underline{2}_{\vee} r \sqsupseteq w) \ \&\ \underline{2}\hat{}\underline{1} \sqsubseteq r))$ from the hypothesis. Suppose, for a reductio, that $M \vDash \underline{2}\hat{}\underline{1} \sqsubseteq r$. Then $M \vDash \underline{1} \sqsubseteq r \ \&\ 2\hat{}\underline{1} \sqsubseteq r$, whence by 1.9 it follows that $\mathfrak{S}^+ \vDash \underline{1} \sqsubseteq r \ \&\ \underline{2}\hat{}\underline{1} \sqsubseteq r$, which is a contradiction. Therefore $M \vDash \neg(\underline{2}\hat{}\underline{1} \sqsubseteq r)$, and so $M \vDash \underline{2}_{\vee} r \sqsubseteq w$, as needed.

Suppose that t is $2^n$ for some $n > 1$. We now argue, under the principal hypothesis, that

(**) $\quad M \vDash \underline{2}^k{}_{\vee} r \sqsupseteq^+ w \rightarrow \underline{2}^{k+1}{}_{\vee} r \sqsupseteq^+ w$

for any k, $0 \leq k < n$. Assume that $M \vDash \underline{2}^k{}_{\vee} r \sqsupseteq^+ w$. Then

$$M \vDash \underline{2}^k{}_{\vee} r \sqsupseteq w \ \&\ \forall x_2 \sqsupseteq w\ (\underline{1} \sqsubseteq x_2 \rightarrow x_2 \sqsupseteq \underline{2}^k{}_{\vee} r).$$

Assuming $M \vDash \underline{1} \sqsubseteq x_2$ for $M \vDash x_2 \sqsupseteq w$ and $x_2 \in M$, we have that
$M \vDash x_2 \sqsupseteq \underline{2}^k{}_{\vee} r$. But $M \vDash \underline{2}^k{}_{\vee} r \sqsupseteq \underline{2}_{\vee}(\underline{2}^k{}_{\vee} r)$ by 1.5, and
$M \vDash \underline{2}_{\vee}(\underline{2}^k{}_{\vee} r) = (\underline{2}_{\vee}\underline{2}^k)_{\vee}\, r = \underline{2}^{k+1}{}_{\vee} r$ by 1.3, and further

$$M \vDash \forall y \sqsupseteq \underline{2}^k{}_{\vee} r\ (\underline{2}^k{}_{\vee} r \sqsupseteq \underline{2}^{k+1}{}_{\vee} r \rightarrow y \sqsupseteq \underline{2}^{k+1}{}_{\vee} r)$$

by 1.17(a). Hence $M \vDash x_2 \sqsupseteq \underline{2}^{k+1}{}_{\vee} r$, and so we have

(##) $\quad M \vDash \forall x_2 \sqsupseteq w\ (\underline{1} \sqsubseteq x_2 \rightarrow x_2 \sqsupseteq \underline{2}^{k+1}{}_{\vee} r).$

From hypothesis $M \vDash \underline{2}^k{}_{\vee} r \sqsupseteq^+ w$ and the principal hypothesis we derive

$M \vDash \underline{2}_{\vee}(\underline{2}^k{}_{\vee} r) \sqsupseteq w \ \text{v}\ (\neg(\underline{2}_{\vee}(\underline{2}^k{}_{\vee} r) \sqsupseteq w) \ \&\ t\hat{}\underline{1} \sqsubseteq \underline{2}^k{}_{\vee} r).$

Suppose, for a reductio, that the latter disjunct holds. Then $M \vDash t\hat{}\underline{1} \sqsubseteq \underline{2}^k{}_{\vee} r$, that is, $M \vDash \underline{2}^n\hat{}\underline{1} \sqsubseteq \underline{2}^k{}_{\vee} r$, hence, by 1.17, also $\mathfrak{S}^+ \vDash \underline{2}^n\hat{}\underline{1} \sqsubseteq \underline{2}^k{}_{\vee} r$. Then $2^n \hat{}\ \underline{1} \hat{}\ r_0 = 2^k \hat{}\ r$ for some $r_0 \in \Sigma^*$, and since $n > k$, we have that $2^j \hat{}\ \underline{1} \hat{}\ r_0 = r$ where $n = k + j$. But this is impossible because $1 \sqsubseteq r$ by the principal hypothesis. Therefore $M \vDash t\hat{}1 \sqsubseteq \underline{2}^k{}_\vee r$ is ruled out, and we must have $M \vDash \underline{2}_\vee(\underline{2}^k{}_\vee r) \sqsupseteq w$, that is, $M \vDash \underline{2}^{k+1}{}_\vee r \sqsupseteq w$. But then along with (##) this shows that $M \vDash \underline{2}^{k+1}{}_\vee r \sqsupseteq^+ w$, which completes the proof of (**). Let $k + 1 = n$. By applying (**) n times we obtain from the hypothesis $M \vDash r \sqsupseteq^+ w$ that $M \vDash \underline{2}^n{}_\vee r \sqsupseteq^+ w$, whence $M \vDash t_\vee r \sqsupseteq w$. □

2.9 Let u and v be any variable-free $\mathcal{L}_T$ terms. Then

$VW^+ \vdash T^*(\underline{\theta(u)}, \underline{\theta(v)},\underline{\theta((u,v))})$ & $\forall w\, (T^*(\underline{\theta(u)}, \underline{\theta(v)},w) \rightarrow \underline{\theta((u,v))} \sqsubseteq w)$.

Proof: We have that $\theta((u,v)) = \theta(u)\hat{}\,\underline{2}^n\hat{}\,\theta(v)$ where n is chosen as in 2.1. By 2.4, $VW^+ \vdash I(\underline{\theta(u)})$ & $MinMaxSep(\underline{2}^j,\underline{\theta(u)})$ & $(\underline{\theta(u)} \neq \underline{1} \rightarrow \underline{2}^j \subseteq_1 \underline{\theta(u)})$ &

& $I(\underline{\theta(v)})$ & $MinMaxSep(\underline{2}^k,\underline{\theta(v)})$ & $(\underline{\theta(v)} \neq \underline{1} \rightarrow \underline{2}^k \subseteq_1 \underline{\theta(v)})$

for appropriate j, k. Again by 2.4,

$\mathfrak{S}^+ \vDash I(\underline{\theta((u,v))})$ & $MinMaxSep(\underline{2}^n,\underline{\theta((u,v))})$ & $\underline{2}^n \subseteq_1 \underline{\theta((u,v))}$,

taking into account that $\mathfrak{S}^+ \vDash \underline{\theta((u,v))} \neq \underline{1}$. By 2.1 we have that $n = \max(j,k)+1$. Hence $\mathfrak{S}^+ \vDash (\underline{2}^j \sqsubseteq \underline{2}^k$ & $\underline{2}^n = \underline{2}^k\hat{}\underline{2}) \vee (\underline{2}^k \sqsubseteq \underline{2}^j$ & $\underline{2}^n = \underline{2}^j\hat{}\underline{2})$. We also have that $\mathfrak{S}^+ \vDash \underline{\theta(u)} \sqsubseteq \underline{\theta((u,v))}$ & $\underline{\theta(v)} \sqsupseteq \underline{\theta((u,v))}$. Finally, we have that $\underline{\theta(u)}\hat{}\underline{2}^n$ and $\underline{2}^n{}_\vee\underline{\theta(v)}$, respectively, uniquely satisfy $\mathfrak{S}^+ \vDash \underline{\theta(u)} \sqsubseteq^+ \underline{\theta(u)}\hat{}\underline{2}^n$ along with $\mathfrak{S}^+ \vDash \underline{\theta(u)}\hat{}\underline{2}^n \sqsubseteq \underline{\theta((u,v))}$ and $\mathfrak{S}^+ \vDash \underline{\theta(v)} \sqsupseteq^+ \underline{2}^n{}_\vee\underline{\theta(v)}$ where $\mathfrak{S}^+ \vDash \underline{2}^n{}_\vee\underline{\theta(v)} \sqsupseteq \underline{\theta((u,v))}$, plus

$\mathfrak{S}^+ \vDash \underline{1}\hat{}\underline{2}^n \sqsupseteq \underline{\theta(u)}\hat{}\underline{2}^n$ & $\underline{2}^n\hat{}\underline{1} \sqsubseteq \underline{2}^n{}_\vee\underline{\theta(v)}$ & $\forall u \sqsubseteq \underline{\theta(u)}\; u \sqsubseteq \underline{\theta((u,v))}$ &

& $\forall u_1 \sqsubseteq \underline{\theta(u)}\hat{}\underline{2}^n\; \forall u_2 \sqsubseteq u_1\; u_2 \sqsubseteq \underline{\theta((u,v))}$ &

$\& \ \forall z \sqsubseteq \underline{\theta(u)}\hat{}\,\underline{2}^n \ (z \sqsubseteq^+ \underline{\theta(u)}\hat{}\,\underline{2}^n \rightarrow$

$$\rightarrow z\hat{}\,\underline{2} \sqsubseteq \underline{\theta(u)}\hat{}\,\underline{2}^n \ \vee \ (\neg z\hat{}\,\underline{2} \sqsubseteq \underline{\theta(u)}\hat{}\,\underline{2}^n \ \& \ \underline{1}\hat{}\,\underline{2}^n \sqsupseteq z)) \ \&$$

$\& \ \forall z \sqsupseteq \underline{2}^n{}_{\vee}\underline{\theta(v)} \ (z \sqsupseteq^+ \underline{2}^n{}_{\vee}\underline{\theta(v)} \rightarrow$

$$\rightarrow \underline{2}_{\vee}z \sqsupseteq \underline{2}^n{}_{\vee}\underline{\theta(v)} \ \vee \ (\neg \underline{2}_{\vee}z \sqsupseteq \underline{2}^n{}_{\vee}\underline{\theta(v)} \ \& \ \underline{2}^n \hat{}\ \underline{1} \sqsupseteq z)) \ \&$$

$\& \ \forall u_2 \sqsubseteq \underline{2}^n{}_{\vee}\underline{\theta(v)} \ (\underline{2}^n \subseteq_1 u_2 \rightarrow \exists! z \sqsubseteq \underline{\theta((u,v))} \ (\underline{\theta(u)} \sqsubseteq z \ \& \ \underline{2}^n \subseteq_1 z \ \& \ u_2 \sqsupseteq z)).$

Then from 1.17(a) it follows that $VW^+ \vdash T^*(\underline{\theta(u)}, \underline{\theta(v)}, \underline{\theta((u,v))})$.

Assume now that $M \vDash T^*(\underline{\theta(u)}, \underline{\theta(v)}, w)$ where $w \in M$. We want to show that

$M \vDash \forall x \sqsubseteq \underline{\theta((u,v))} \ x \sqsubseteq w$. Suppose $M \vDash x \sqsubseteq \underline{\theta((u,v))}$, that is,
$M \vDash x \sqsubseteq (\underline{\theta(u)}\hat{}\,\underline{2}^n)\hat{}\,\underline{\theta(v)}$. Then, by 1.1, we have $M \vDash \bigvee_{p \sqsubseteq \theta((u,v))} x = \underline{p}$. Now, for any such $p \in \Sigma^*$, we have that (a) $p \sqsubseteq \theta(u)$, or (b) $\theta(u) \sqsubseteq p$ and $p \sqsubseteq \theta(u)\hat{}\,\underline{2}^n$, or else (c) $\theta(u)\hat{}\,\underline{2}^n \sqsubseteq p$. If (a), then $\mathfrak{S}^+ \vDash \underline{p} \sqsubseteq \underline{\theta(u)}$, hence, by 1.5, also $M \vDash \underline{p} \sqsubseteq \underline{\theta(u)}$, and so from $M \vDash T^*(\underline{\theta(u)}, \underline{\theta(v)}, w)$ we have $M \vDash \underline{p} \sqsubseteq w$. If (b), then $M \vDash \underline{\theta(u)} \sqsubseteq \underline{p} \ \& \ \underline{p} \sqsubseteq \underline{\theta(u)}\hat{}\,\underline{2}^n$, again by 1.5. From the principal hypothesis $M \vDash T^*(\underline{\theta(u)}, \underline{\theta(v)}, w)$ it follows that

$$M \vDash \exists! w_1 \sqsubseteq w \ (\underline{1}\hat{}\,\underline{2}^n \sqsupseteq w_1 \ \& \ \underline{\theta(u)} \sqsubseteq^+ w_1)$$

along with $M \vDash \underline{1} \sqsupseteq \underline{\theta(u)} \ \& \ Tally_2(\underline{2}^n)$ since $\mathfrak{S}^+ \vDash \underline{1} \sqsupseteq \underline{\theta(u)} \ \& \ Tally_2(\underline{2}^n)$.

(Here we appeal to 1.9(a) and 1.17(a).) Also from the principal hypothesis we have that

$$M \vDash \forall z \sqsubseteq w_1 \ (z \sqsubseteq^+ w_1 \ \rightarrow \ z\hat{}\,\underline{2} \sqsubseteq w_1 \ \vee \ (\neg(z\hat{}\,\underline{2} \sqsubseteq w_1) \ \& \ \underline{1}\hat{}\,\underline{2}^n \sqsupseteq z)).$$

Hence from 2.7, letting s be $\underline{\theta(u)}$ and t be $\underline{2}^n$ and instantiating with $w_1$, we obtain $M \vDash \underline{\theta(u)}\hat{}\,\underline{2}^n \sqsubseteq w_1$. But then $M \vDash \underline{p} \sqsubseteq w$ follows from $M \vDash \underline{p} \sqsubseteq \underline{\theta(u)}\hat{}\,\underline{2}^n$ and the principal hypothesis since $M \vDash \forall z \sqsubseteq w_1 \ \forall z_1 \sqsubseteq z \ z_1 \sqsubseteq w$.

Finally, suppose (c ), that is, $\theta(u)\,\hat{}\;\underline{2}^n \sqsubseteq p$. Then

$$\mathfrak{S}^+ \vDash \exists q \sqsubseteq \underline{2}^n{}_\vee\underline{\theta(v)}\ (\underline{2}^n \subseteq_1 q \ \ \&\ \ q \sqsupseteq \underline{p}),$$

as well as $\mathfrak{S}^+ \vDash \underline{\theta(u)} \sqsubseteq \underline{p} \ \ \&\ \ \underline{2}^n \subseteq_1 \underline{p}$. Hence by 1.17(a),

$$M \vDash \underline{q} \sqsubseteq \underline{2}^n{}_\vee\underline{\theta(v)} \ \ \&\ \ \underline{2}^n \subseteq_1 \underline{q} \ \ \&\ \ \underline{\theta(u)} \sqsubseteq \underline{p} \ \ \&\ \ \underline{2}^n \subseteq_1 \underline{p} \ \ \&\ \ \underline{q} \sqsupseteq \underline{p}$$

and $M \vDash \underline{\theta(v)} \sqsupseteq^+ w_2 \ \&\ \underline{2}^n\,\hat{}\,\underline{1} \sqsubseteq w_2 \ \&\ \underline{1} \sqsubseteq \underline{\theta(v)} \ \&\ \mathrm{Tally}_2(\underline{2}^n)$ for some unique $w_2 \in M$ such that $M \vDash w_2 \sqsupseteq w$. Letting r be $\underline{\theta(v)}$ and t be $\underline{2}^n$, we obtain, from 2.8, by instantiating with $w_2$, that $M \vDash \underline{2}^n{}_\vee\underline{\theta(v)} \sqsupseteq w_2$. Then from hypothesis $M \vDash T^*(\underline{\theta(u)}, \underline{\theta(v)}, w)$ we have that $M \vDash \exists! z \sqsubseteq w\ (\underline{\theta(u)} \sqsubseteq z \ \&\ \underline{2}^n \subseteq_1 z \ \&\ \underline{q} \sqsupseteq z)$. Therefore $M \vDash z = \underline{p}$, so $M \vDash \underline{p} \sqsubseteq w$ as claimed. We have thus established that $M \vDash \forall x \sqsubseteq \underline{\theta((u,v))}\ x \sqsubseteq w$. But $\mathfrak{S}^+ \vDash \underline{\theta((u,v))} \sqsubseteq \underline{\theta((u,v))}$, so by 1.3 we obtain $M \vDash \underline{\theta((u,v))} \sqsubseteq \underline{\theta((u,v))}$. It follows that $M \vDash \underline{\theta((u,v))} \sqsubseteq w$. □

2.10 Let u and v be any variable-free $\mathcal{L}_T$ terms. Then

$$VW^+ \ \vdash T^+(\underline{\theta(u)}, \underline{\theta(v)}, \underline{\theta((u,v))}).$$

<u>Proof</u>: By 2.9, $VW^+ \vdash T^*(\underline{\theta(u)}, \underline{\theta(v)}, \underline{\theta((u,v))})$, and furthermore,

$$VW^+ \ \vdash \ \forall v_1\ (T^*(\underline{\theta(u)}, \underline{\theta(v)}, v_1) \ \rightarrow \ \underline{\theta((u,v))} \sqsubseteq v_1).$$

Hence $VW^+ \vdash \exists w\ (T^*(\underline{\theta(u)}, \underline{\theta(v)}, w)\ \&$

$$\&\ \forall v_1\ (T^*(\underline{\theta(u)}, \underline{\theta(v)}, v_1) \ \rightarrow \ \underline{\theta((u,v))} \sqsubseteq v_1)\ \&\ w = \underline{\theta((u,v))}.$$

But then $VW^+ \ \vdash T^+(\underline{\theta(u)}, \underline{\theta(v)}, \underline{\theta((u,v))})$ follows from the definition of $T^+$. □

We let $\exists! x \subseteq_p t\ A(x) \equiv: \exists x \subseteq_p t\ (A(x)\ \&\ \forall y \subseteq_p t\ (A(y) \rightarrow y = t))$, assuming that the variable x does not occur in t.

We then have:

2.11 Let t be any variable-free $\mathcal{L}_T$ term of the form $(t_1,t_2)$. Then

$$VW^+ \vdash \exists! y_1 \subseteq_p \underline{\theta(t)}\ \exists! y_2 \subseteq_p \underline{\theta(t)}\ (I(y_1) \ \&\ I(y_2) \ \&\ T^*(y_1,y_2,\underline{\theta(t)})).$$

Proof: We first show that

$$\mathfrak{S}^+ \vDash \exists! y_1 \subseteq_p \underline{\theta(t)}\ \exists! y_2 \subseteq_p \underline{\theta(t)}\ (I(y_1) \ \&\ I(y_2) \ \&\ T^*(y_1,y_2,\underline{\theta((u,v))})).$$

Now from 2.4, 2.9 and 1.17 we have that

$$\mathfrak{S}^+ \vDash I(\underline{\theta(t_1)}) \ \ \&\ \ I(\underline{\theta(t_2)}) \ \ \&\ \ T^*(\underline{\theta(t_1)}, \underline{\theta(t_2)},\underline{\theta((t_1,t_2))}).$$

We need to show that

$\mathfrak{S}^+ \vDash \forall y_1,y_2 \subseteq_p \underline{\theta(t)}(\ I(y_1) \ \&\ I(y_2) \ \&\ T^*(y_1,y_2.,\underline{\theta(t)}) \rightarrow y_1 = \underline{\theta(t_1)} \ \&\ y_2 = \underline{\theta(t_2)})$.

So assume $\mathfrak{S}^+ \vDash T^*(y_1,y_2.,\underline{\theta(t)})$ for some $y_1,y_2 \in \Sigma^*$ such that

$\mathfrak{S}^+ \vDash I(y_1) \ \&\ I(y_2) \ \&\ y_1 \subseteq_p \underline{\theta(t)} \ \&\ y_2 \subseteq_p \underline{\theta(t)}$. Let $\mathfrak{S}^+ \vDash MinMaxSep(\underline{2}^n,\underline{\theta(t)})$. Then we have that

$$\mathfrak{S}^+ \vDash y_1 \sqsubseteq \underline{\theta(t)} \ \&\ y_1 \sqsubseteq^+ \ \underline{\theta(t_1)}{}^\wedge\underline{2}^n \ \ \&\ y_2 \sqsupseteq \underline{\theta(t)} \ \&\ y_2 \sqsupseteq^+ \ \underline{2}^n{}_\vee\underline{\theta(t_2)}.$$

Then from $\mathfrak{S}^+ \vDash I(\underline{\theta(t_1)})$ we have $\mathfrak{S}^+ \vDash \underline{\theta(t_1)} \sqsubseteq y_1$. On the other hand, from $\mathfrak{S}^+ \vDash I(y_1) \ \&\ \ y_1 \sqsubseteq \underline{\theta(t_1)}{}^\wedge\underline{2}^n \ \&\ y_1 \sqsubseteq \underline{\theta(t)}$ and $\mathfrak{S} \vDash \underline{\theta(t_1)} \sqsubseteq^+ \ \underline{\theta(t_1)}{}^\wedge\underline{2}^n$ we also obtain $\mathfrak{S}^+ \vDash y_1 \sqsubseteq \underline{\theta(t_1)}$. Hence $\mathfrak{S}^+ \vDash y_1 = \underline{\theta(t_1)}$ by (VW10). Using $\mathfrak{S}^+ \vDash y_2 \sqsupseteq^+ \ \underline{2}^n{}_\vee\underline{\theta(t_2)}$ and $\mathfrak{S}^+ \vDash \underline{\theta(t_2)} \sqsupseteq^+ \ \underline{2}^n{}_\vee\underline{\theta(t_2)}$ we analogously obtain $\mathfrak{S}^+ \vDash y_2 = \underline{\theta(t_2)}$. We then derive 2.11 by applying 1.17(a). □

## §3. Finding Subtrees among Substrings

Note that under minimax coding a string code $\theta(u)$ of a tree u may contain substrings $\theta(v)$ that code trees v that are not subtrees of u: e.g., $\theta((00)(00)) = 12122^2121$ has a substring $12^2121 = \theta(0(00))$, but 0(00) is not a subtree of (00)(00). To interpret the subtree relation of $\mathcal{L}_T$ in $\mathcal{L}^+$ we will need to somehow distinguish from among the substrings of a tree's string code those that code its subtrees. To do this we explicitly represent in $\mathcal{L}^+$ the process by which u is constructed from its component subtrees.

E.g., with the tree (00)(00) we may associate a sequence of ordered pairs such as |(211,1)|(212,1)|(221,1)|(222,1)|(21,121)|(22,121)|(2,$1212^2121$), with each ordered pair standing for a distinct node of the tree, its first term being a dyadic string and second term the minimax code of a subtree of (00)(00). The last ordered pair in the sequence is $(2,\theta((00)(00)))$, and each ordered pair p with second term distinct from 1 – which $= \theta(0)$ – is preceded by ordered pairs $p_1,p_2$ whose first terms are dyadic strings that extend the first term of p with a single digit 1 or 2, resp., and having as their second terms tree codes $\theta(u_1)$, $\theta(u_2)$, resp., of the subtrees $u_1$, $u_2$ of (00)(00). The root nodes of $u_1$, $u_2$ are the immediate descendants in (00)(00) of the node characterized by p, which has as its second term the minimax code $\theta((u_1,u_2))$ of the tree $(u_1,u_2)$. We say that such a sequence of ordered pairs is a tree map of (00)(00). It is clear that distinct dyadic trees have distinct tree maps.

To define tree maps in $\mathcal{L}^+$ we need to code ordered pairs of dyadic strings by single dyadic strings, and then to uniquely represent sequences of such codes by dyadic strings while simultaneously expressing the relation between the pairs $p_1$, $p_2$, p in terms of the appropriate relation between the corresponding tree codes $\theta(u_1)$, $\theta(u_2)$, $\theta((u_1,u_2))$. For the latter we shall use the relation $T^*(x,y,w)$ defined in §2. For the purpose of coding ordered pairs of dyadic strings and sequences thereof we will use suitably chosen 2-tallies, called keys, as separators, of sufficient length depending on the mapped tree.

Let $\quad FrPair(u,x,y_1,y_2) \equiv: \ x \subseteq_1 u \ \& \ y_1 \subseteq_p u \ \& \ y_2 \subseteq_p u \ \&$

$\& \ \exists w_1 \sqsubseteq u \ \exists w_2 \sqsupseteq u \ (\underline{1}_\vee x \sqsupseteq w_1 \ \& \ \underline{1}_\vee(y_1\hat{}\underline{1}) \sqsubseteq^+ w_1 \ \& \ x\hat{}\underline{1} \sqsubseteq w_2 \ \&$

$\& \ (\underline{1}_\vee y_2)\hat{}\underline{1} \sqsupseteq^+ w_2) \ \& \ \underline{1} \sqsupseteq \underline{1}_\vee(y_1\hat{}\underline{1}) \ \& \ \underline{1} \sqsubseteq (\underline{1}_\vee y_2)\hat{}\underline{1} \ \& \ \forall w \sqsubseteq w_1 \ \underline{1} \sqsupseteq w\hat{}\underline{1} \ \&$

$\& \ \forall w \sqsubseteq u \ (\underline{1}_\vee x \sqsupseteq w \rightarrow x \sqsupseteq w) \ \& \ \forall w \sqsubseteq u \ (w \sqsubset u \rightarrow w\hat{}\underline{1} \sqsubseteq u \vee w\hat{}\underline{2} \sqsubseteq u) \ \&$

$\& \ \forall w \sqsubseteq u \ \neg w\hat{}\underline{1} \sqsubseteq w \ \& \ \forall w \sqsubseteq u \ \neg w\hat{}\underline{2} \sqsubseteq w \ \& \ \forall w \sqsubseteq u \ \forall v \sqsubseteq w \ v \sqsubseteq u \ \&$

$\& \ \forall u_1,u_2 \subseteq_p u \ (\underline{1}_\vee(u_1\hat{}\underline{1}) = \underline{1}_\vee(u_2\hat{}\underline{1}) \rightarrow \ u_1 = u_2) \ \&$

$$\& \ \forall u_1,u_2 \subseteq_p u \ ((\underline{1}_\vee u_1)\hat{}\underline{1} = (\underline{1}_\vee u_2)\hat{}\underline{1} \rightarrow \ u_1 = u_2)),$$

reading "u is ordered pair of $y_1,y_2$ framed by x".

3.1 $\quad VW^+ \vDash \forall x,t,y_1,y_2,z_1,z_2 \ (FrPair(x,t,y_1,y_2) \ \& \ FrPair(x,t,z_1,z_2) \ \rightarrow$

$$\rightarrow \ y_1 = z_1 \ \& \ y_2 = z_2).$$

<u>Proof</u>: Suppose $M \vDash FrPair(x,t,y_1,y_2) \ \& \ FrPair(x,t,z_1,z_2)$ where $x,t,y_1,y_2,z_1,z_2 \in M$. Then

$M \vDash w_1 \sqsubseteq x \ \&\ w_3 \sqsubseteq x \ \&\ \underline{1}_\vee t \sqsupseteq w_1 \ \&\ \underline{1}_\vee t \sqsupseteq w_3 \ \&\ t \subseteq_1 x \ \&$

$\&\ \underline{1}_\vee(y_1\hat{}\underline{1}) \sqsubseteq^+ w_1 \ \&\ \underline{1}_\vee(z_1\hat{}\underline{1}) \sqsubseteq^+ w_3$

for some $w_1,w_3 \in M$. Then $M \vDash \forall v \sqsubseteq w_1 (\underline{1} \sqsupseteq v \to v \sqsubseteq \underline{1}_\vee(y_1\hat{}\underline{1}))$ along with

$M \vDash t \sqsupseteq w_1 \ \&\ t \sqsupseteq w_3 \ \&\ \underline{1} \sqsupseteq \underline{1}_\vee(z_1\hat{}\underline{1}) \ \&\ \underline{1}_\vee(z_1\hat{}\underline{1}) \sqsubseteq w_3$ . Now from

$$M \vDash w_1 \sqsubseteq x \ \&\ w_3 \sqsubseteq x \ \&\ t \sqsupseteq w_1 \ \&\ t \sqsupseteq w_3 \ \&\ t \subseteq_1 x$$

it follows that $M \vDash w_1 = w_3$. Hence we also have $M \vDash \underline{1}_\vee( z_1\hat{}\underline{1}) \sqsubseteq^+ w_1$. But then $M \vDash \underline{1}_\vee( z_1\hat{}\underline{1}) \sqsubseteq \underline{1}_\vee(y_1\hat{}\underline{1})$. A symmetric argument shows that also

$M \vDash \underline{1}_\vee(y_1\hat{}\underline{1}) \sqsubseteq \underline{1}_\vee(z_1\hat{}\underline{1})$. Hence $M \vDash \underline{1}_\vee(y_1\hat{}\underline{1}) = \underline{1}_\vee(z_1\hat{}\underline{1})$ by (VW10).

From the principal hypothesis we have that $M \vDash y_1 \subseteq_p x \ \&\ z_1 \subseteq_p x$, hence it follows that $M \vDash y_1 = z_1$, again from the principal hypothesis. An analogous argument shows that $M \vDash y_2 = z_2$. □

Let

$Key_m(x,y) \ \equiv: \ \exists t \subseteq_p y(MinMaxSep(t,y) \ \&\ x = t\hat{}\underline{2}^2)$, "x is minor key for y",

$Key_M(x,y) \ \equiv: \ \exists t \subseteq_p y(MinMaxSep(t,y) \ \&\ x = t\hat{}\underline{2}^3$, "x is major key for y".

3.2 Let $r_1, r_2, t_1, t_2$ be variable-free $\mathcal{L}^+$-terms. Then

$\mathfrak{S}^+ \vDash \forall z_1,z_2 (FrPair(z_1,t_1,r_1,r_2) \ \&\ FrPair(z_2,t_2,r_1,r_2) \ \&\ Key_m(t_1,r_2) \ \&$

$\&\ Key_m(t_2,r_2) \to t_1 = t_2 \ \&\ z_1 = z_2)$.

<u>Proof</u>: Assume $\mathfrak{S}^+ \vDash FrPair(z_1,t_1,r_1,r_2) \ \&\ FrPair(z_2,t_2,r_1,r_2)$ where $\mathfrak{S}^+ \vDash Key_m(t_1,r_2) \ \&\ Key_m(t_2,r_2)$ and $z_1,t_1,z_2,t_2,r_1,r_2 \in \Sigma^*$. From 2.2 and the definition of $\mathrm{Key}_m$ we have that $\mathfrak{S}^+ \vDash t_1 = t_2$. We have from the hypothesis that $\mathfrak{S}^+ \vDash t_1 \subseteq_1 z_1$ & $\underline{1}_{\vee}t_1 \sqsupseteq w_1$ & $\underline{1}_{\vee}(r_1\hat{}\,\underline{1}) \sqsubseteq^+ w_1$ & $w_1 \sqsubseteq z_1$ &

$$\& \; t_1 \subseteq_1 z_2 \;\&\; \underline{1}_{\vee}t_1 \sqsupseteq w_3 \;\&\; \underline{1}_{\vee}(r_1\hat{}\,\underline{1}) \sqsubseteq^+ w_3 \;\&\; w_3 \sqsubseteq z_2$$

for some $w_1,w_3 \in \Sigma^*$. We argue that for no $r \in \Sigma^*$,

$$\mathfrak{S}^+ \vDash \underline{r} \neq \underline{e} \;\&\; ((\underline{1}_{\vee}(r_1\hat{}\,\underline{1}))\hat{}\,\underline{r})\hat{}\,t_1 = w_1.$$

Suppose for, a reductio, that $\mathfrak{S}^+ \vDash ((\underline{1}_{\vee}(r_1\hat{}\,\underline{1}))\hat{}\,\underline{r})\hat{}\,t_1 = w_1$ for some $r \in \Sigma^*$, $r \neq e$. Assume $\mathfrak{S}^+ \vDash \underline{2} \sqsupseteq \underline{r}$. But then $\mathfrak{S}^+ \vDash \underline{2}\hat{}\,t_1 \sqsupseteq w_1$ along with $\mathfrak{S}^+ \vDash \underline{1}\hat{}\,t_1 \sqsubseteq t_1$, whence $\mathfrak{S}^+ \vDash \mathrm{Tally}_2(\underline{2}\hat{}\,t_1) \subseteq_p r_2$ and $\mathfrak{S}^+ \vDash \underline{1} = \underline{2}$, a contradiction. Therefore $\mathfrak{S}^+ \vDash \underline{1} \sqsupseteq \underline{r}$. But then $\mathfrak{S}^+ \vDash (\underline{1}_{\vee}(r_1\hat{}\,\underline{1}))\hat{}\,\underline{r} \sqsubseteq w_1$ & $\underline{1} \sqsupseteq \underline{r}$, whereas from $\mathfrak{S}^+ \vDash \underline{1}_{\vee}(r_1\hat{}\,\underline{1}) \sqsubseteq^+ w_1$ it follows that $\mathfrak{S}^+ \vDash (\underline{1}_{\vee}(r_1\hat{}\,\underline{1}))\hat{}\,\underline{r} \sqsubseteq \underline{1}_{\vee}(r_1\hat{}\,\underline{1})$, a contradiction. Hence there is no such $r \in \Sigma^*$, and it follows that $\mathfrak{S}^+ \vDash (\underline{1}_{\vee}(r_1\hat{}\,\underline{1}))\hat{}\,t_1 = w_1$.

The same argument with $z_1$ and $w_1$ replaced by $z_2$ and $w_3$ shows that $\mathfrak{S}^+ \vDash (\underline{1}_{\vee}(r_1\hat{}\,\underline{1}))\hat{}\,t_1 = w_3$. Hence $\mathfrak{S}^+ \vDash w_1 = w_3$.

On the other hand, from the hypothesis we also have that

$\mathfrak{S}^+ \vDash t_1 \subseteq_1 z_1$ & $\underline{t_1}_{\vee}\underline{1} \sqsubseteq w_2$ & $\underline{1}_{\vee}(r_2\hat{}\,\underline{1}) \sqsupseteq^+ w_2$ & $w_2 \sqsupseteq z_1$ &

$$\& \; t_1 \subseteq_1 z_2 \;\&\; t_1\hat{}\,\underline{1} \sqsubseteq w_4 \;\&\; \underline{1}_{\vee}(r_2\hat{}\,\underline{1}) \sqsubseteq^+ w_4 \;\&\; w_4 \sqsupseteq z_2$$

for some $w_2,w_4 \in \Sigma^*$. Then exactly analogous arguments show that $\mathfrak{S}^+ \vDash w_2 = t_1\hat{}\,((\underline{1}_{\vee}r_2)\hat{}\,\underline{1}) = w_4$.

But then we have that $\mathfrak{S}^+ \vDash z_1 = ((\underline{1}_{\vee}(r_1\hat{}\,\underline{1}))\hat{}\,t_1)\hat{}\,((\underline{1}_{\vee}r_2)\hat{}\,\underline{1}) = z_2$. □

Next, we let $\mathrm{Mapframe}(x,x_1,x_2)$ abbreviate the $\mathcal{L}^+$ formula

$Tally_2(x_1)$ & $Tally_2(x_2)$ & $x_2\hat{}\underline{1} \sqsubseteq x$ & $x_2 = x_1\hat{}\underline{2}$ & $\neg x_2\hat{}\underline{2} \subseteq_p x$ &

& $\forall y \subseteq_p x\ (x_2\hat{}\underline{1} \sqsubseteq y\ \&\ \underline{1}_{\mathrm{v}}x_2 \sqsupseteq y\ \&\ \neg(\underline{1}_{\mathrm{v}}x_2)\hat{}\underline{1} \subseteq_p y \rightarrow (\underline{1}_{\mathrm{v}}x_1)\hat{}\underline{1} \subseteq_1 y)$,

reading "x is a map frame with keys $x_1,x_2$".

Finally, let Treemap(x,y) abbreviate the $\mathcal{L}^+$ formula

$\exists x_1,x_2 \subseteq_p x\ (Key_m(x_1,y)\ \&\ Key_M(x_2,y)\ \&\ Mapframe(x,x_1,x_2)$ & (tm1)-(tm6)

where: (tm1) $\exists z \sqsupseteq x\ (FrPair(z,x_1,\underline{2},y)$,

(tm2) $\forall z \subseteq_p x\ \forall x_0 \subseteq_p z\ \forall y_1 |\leq| x_1\ \forall y_2 \subseteq_p z\ (FrPair(z,x_0,y_1,y_2) \rightarrow x_0 = x_1\ \&$

$\&\ \forall w \subseteq_p x\ \forall y_3 \subseteq_p x\ (FrPair(w,x_1,y_1,y_3) \rightarrow y_3 = y_2\ \&\ w = z))$,

(tm3) $\forall z \subseteq_p x\ \forall y_1 |\leq| x_1\ \forall y_2 \subseteq_p z\ (FrPair(z,x_1,y_1,y_2) \rightarrow (y_2 = \underline{1}\ \&\ \underline{2} \sqsubseteq y_1)\ \mathrm{v}$

$\mathrm{v}\ ((y_2 \neq \underline{1}\ \&\ I(y_2)\ \&\ \exists z_1,z_2 \subseteq_p x\ \exists u \subseteq_p z_1\ \exists v \subseteq_p z_2\ (FrPair(z_1,x_1,y_1\hat{}\underline{1},u)\ \&$

$\&\ FrPair(z_2,x_1,y_1\hat{}\underline{2},v)\ \&\ I(u)\ \&\ I(v)\ \&\ T^*(u,v,y_2))))$,

(tm4) $\forall z \subseteq_p x\ \forall y_1 |\leq| x_1\ \forall y_2 \subseteq_p z\ (FrPair(z,x_1,y_1,y_2)\ \&\ y_2 \neq y \rightarrow$

$\rightarrow \exists w,z' \subseteq_p x\ \exists w_1 |\leq| x_1\ \exists w_2 \subseteq_p w\ \exists w_3,w_4 \subseteq_p w\ (FrPair(w,x_1,w_1,w_2)\ \&$

$\&\ FrPair(z',x_1,w_3,w_4)\ \&\ (((y_1 = w_1\hat{}\underline{1}\ \&\ w_3 = w_1\hat{}\underline{2}\ \&\ T^*(y_2,w_4,w))\ \mathrm{v}$

$\mathrm{v}\ (w_3 = w_1\hat{}\underline{1}\ \&\ y_1 = w_1\hat{}\underline{2}\ \&\ T^*(w_4,y_2,w)))))$,

(tm5) $\forall w \subseteq_p x\ \forall y_3 |\leq| x_1\ \forall y_4 \subseteq_p w\ (FrPair(w,x_1,y_3,y_4) \rightarrow (y_3 = \underline{2} \leftrightarrow y_4 = y))$,

(tm6) $\forall z \subseteq_p x\ \forall y_1,y_2 \subseteq_p z\ (FrPair(z,x_1,y_1,y_2) \rightarrow \underline{2} \sqsubseteq y_1\ \&\ y_1 |\leq| x_1)$

reading "x is a tree map for y". Then we let $u \ll_x v$ abbreviate

$u \subseteq_p v\ \&\ \exists x_1 \subseteq_p x\ (Treemap(x,v)\ \&\ Key_m(x_1,v)\ \&$

$\&\ \exists z_1,z_2 \subseteq_p x\ \exists u' |\leq| x_1\ (FrPair(z_1,x_1,u',u)\ \&\ FrPair(z_2,x_1,\underline{2},v)\ \&\ \underline{2} \sqsubseteq u'))$.

E.g., if t is $\underline{0}$, then $\theta(t) = \underline{1}$, which has $\underline{2}^0 = e$ as its minimax separator, and $\underline{2}^2$ and $\underline{2}^3$ as the minor and major keys, resp.. Then $2^3\,\hat{}\,1\,\hat{}\,2\,\hat{}\,1\,\hat{}\,2^2\,\hat{}\,1\,\hat{}\,1\,\hat{}\,1\,\hat{}\,2^3$ is a tree map for $\theta(t)$ containing (2,1) as its sole framed pair. (To identify the members of the framed pair, ignore the first 1's immediately following and immediately preceding occurrences of the major and minor keys.) For $t \equiv (00)(00)$, we have that $\theta(t)$ is $1\,\hat{}\,2\,\hat{}\,1\,\hat{}\,2^2\,\hat{}\,1\,\hat{}\,2\,\hat{}\,1$, which has $\underline{2}^4$ and $\underline{2}^5$ as its minor and major keys, resp.. To obtain a tree map $q^t$ of $\theta(t)$ in which each framed ordered pair occurs exactly once we first obtain framed ordered pairs $q^{t,s}$ corresponding to each of the subtrees s of t and then juxtapose them together in appropriate order in a single string using the major key. We have that $q^{t,2} = 1\,\hat{}\,2\,\hat{}\,1\,\hat{}\,2^4\,\hat{}\,1\,\hat{}\,\theta(t)\,\hat{}\,1$ is the ordered pair $(2,\theta(t))$ framed by $2^4$. Now, (00) and (00) are the immediate subtrees of t, and $\theta((0,0)) = 1\,\hat{}\,2\,\hat{}\,1$. Then

$$q^{t,21} = 1\,\hat{}\,2\,\hat{}\,1\,\hat{}\,1\,\hat{}\,2^4\,\hat{}\,1\,\hat{}\,\theta((0,0))\,\hat{}\,1 \text{ and } q^{t,22} = 1\,\hat{}\,2\,\hat{}\,2\,\hat{}\,1\,\hat{}\,2^4\,\hat{}\,1\,\hat{}\,\theta((0,0))\,\hat{}\,1$$

are the ordered pairs of $21,\theta((0,0))$ and $22,\theta((0,0))$, resp., both framed by $2^4$. Finally, 0 and 0 are the immediate subtrees of (00), and $\theta(0) = \underline{1}$. Then

$$q^{t,211} = 1\,\hat{}\,2\,\hat{}\,1\,\hat{}\,1\,\hat{}\,1\,\hat{}\,2^4\,\hat{}\,1\,\hat{}\,\theta(0)\,\hat{}\,1, \quad q^{t,212} = 1\,\hat{}\,2\,\hat{}\,1\,\hat{}\,2\,\hat{}\,1\,\hat{}\,2^4\,\hat{}\,1\,\hat{}\,\theta(0)\,\hat{}\,1,$$

$$q^{t,221} = 1\,\hat{}\,2\,\hat{}\,2\,\hat{}\,1\,\hat{}\,1\,\hat{}\,2^4\,\hat{}\,1\,\hat{}\,\theta(0)\,\hat{}\,1, \quad q^{t,222} = 1\,\hat{}\,2\,\hat{}\,2\,\hat{}\,2\,\hat{}\,1\,\hat{}\,2^4\,\hat{}\,1\,\hat{}\,\theta(0)\,\hat{}\,1$$

are the ordered pairs of $211,\theta(0)$, $212,\theta(0)$, $221,\theta(0)$ and $222,\theta(0)$, all framed by the minor key of $\theta(t)$, namely $2^4$. Then

$$2^5\,\hat{}\,q^{t,221}\,\hat{}\,2^5\,\hat{}\,q^{t,222}\,\hat{}\,2^5\,\hat{}\,q^{t,211}\,\hat{}\,2^5\,\hat{}\,q^{t,212}\,\hat{}\,2^5\,\hat{}\,q^{t,21}\,\hat{}\,2^5\,\hat{}\,q^{t,22}\,\hat{}\,2^5\,\hat{}\,q^{t,2}$$

is a tree map for $\theta((00)(00))$. On account of this straightforward algorithm, we may assert:

3.3 For each variable-free $\mathcal{L}_T$ term t there is a $q^{[t]} \in \Sigma^*$ such that

$$\mathfrak{S}^+ \vDash \text{Treemap}(\underline{\theta(t)},\underline{q}^{[t]}).$$

Furthermore, from 1.17 we have that $VW^+ \vdash \text{Treemap}(\underline{\theta(t)},\underline{q}^{[t]})$.

Let t be a variable-free $\mathcal{L}_T$ term, and $\mathcal{S}[t]$ the set of its subterms. We let $\mathcal{S}^{\sqsubseteq}[t]$ be the smallest subset of $\Sigma^* \times \mathcal{S}[t]$ satisfying the following conditions:

(i) $<\underline{2},t> \in \mathcal{S}^{\sqsubseteq}[t]$,

(ii) if $<p,s> \in \mathcal{S}^{\sqsubseteq}[t]$ and $s = (s_1,s_2)$, then $<p\hat{\ }\underline{1},s_1> \in \mathcal{S}^{\sqsubseteq}[t]$ and $<p\hat{\ }\underline{2},s_2> \in \mathcal{S}^{\sqsubseteq}[t]$.

A straightforward induction shows that for each $s \in \mathcal{S}[t]$, $<p,s> \in \mathcal{S}^{\sqsubseteq}[t]$ for some dyadic string $p \in \Sigma^*$ with $\underline{2}$ as its initial digit. Moreover, for each $p \in \Sigma^*$ there is at most one $s \in \mathcal{S}[t]$ such that $<p,s> \in \mathcal{S}^{\sqsubseteq}[t]$.

3.4 (a) Let t be a variable-free $\mathcal{L}_T$ term and s any subterm of t. For any dyadic string p, if $<p,s> \in \mathcal{S}^{\sqsubseteq}[t]$, then

$$VW^+ \vdash \exists! z \subseteq_p \underline{q}^{[t]}\ \exists! x_1 \subseteq_p z\ \text{FrPair}(z,x_1,\underline{p},\underline{\theta(s)}).$$

(b) For each $p,u \in \Sigma^*$, if $VW^+ \vdash \exists z \subseteq_p \underline{q}^{[t]}\ \exists x_1 \subseteq_p z\ \text{FrPair}(z,x_1,\underline{p},\underline{u})$, then $u = \theta(s)$ for some $<p,s> \in \mathcal{S}^{\sqsubseteq}[t]$.

Proof: We first show that

$$\mathfrak{S}^+ \vDash \exists! z \subseteq_p \underline{q}^{[t]}\ \exists! x_1 \subseteq_p z\ (\text{FrPair}(z,x_1,\underline{p},\underline{\theta(s)})\ \&\ z \subseteq_1 \underline{q}^{[t]}).$$

Fix t, and let $2^m$ be the minor key for t. We argue by induction on the generating relation of $\mathcal{S}^{\sqsubseteq}[t]$. Suppose p = 2 and s is t. Then $<2,t> \in \mathcal{S}^{\sqsubseteq}[t]$, and for $z = 1\hat{}2\hat{}1\hat{}2^m\hat{}1\hat{}\theta(t)\hat{}1$ we have $\mathfrak{S}^+ \vDash \mathrm{FrPair}(\underline{z},\underline{2}^m,\underline{2},\underline{\theta(t)})$. By 3.2 such z and $x_1$ are unique. From $\mathfrak{S}^+ \vDash \mathrm{Treemap}(\underline{\theta(t)},\underline{q}^{[t]})$ we obtain $\mathfrak{S}^+ \vDash \underline{z} \subseteq_p \underline{q}^{[t]}$ & $\underline{z} \subseteq_1 \underline{q}^{[t]}$, as needed. Suppose now that

$$\mathfrak{S}^+ \vDash \exists! z \subseteq_p \underline{q}^{[t]}\ \exists! x_1 \subseteq_p z(\mathrm{FrPair}(z,x_1,\underline{r},\underline{\theta(s)}) \ \&\ z \subseteq_1 \underline{q}^{[t]})$$

for some $s \in \mathcal{S}[t]$ where $<r,s> \in \mathcal{S}^{\sqsubseteq}[t]$ and $s = (s_1,s_2)$. Then for

$z_1 = 1\hat{}r\hat{}1\hat{}1\hat{}2^m\hat{}1\hat{}\theta(s_1)\hat{}1$ and $z_2 = 1\hat{}r\hat{}2\hat{}1\hat{}2^m\hat{}1\hat{}\theta(s_2)\hat{}1$ we have that $\mathfrak{S}^+ \vDash \mathrm{FrPair}(\underline{z_1},\underline{2}^m,\underline{r}\hat{}\underline{1},\underline{\theta(s_1)})$ and $\mathfrak{S}^+ \vDash \mathrm{FrPair}(\underline{z_2},\underline{2}^m,\underline{r}\hat{}\underline{2},\underline{\theta(s_2)})$. We have that $\theta(s) \neq 1$ because $s = (s_1,s_2)$. Hence from the inductive hypothesis and 3.3 we obtain that $\mathfrak{S}^+ \vDash z_1 \subseteq_p \underline{q}^{[t]}$ & $z_2 \subseteq_1 \underline{q}^{[t]}$ as claimed. Again, uniqueness follows by 3.2. Then 3.4(a) follows by 1.17(a). Part (b) follows from 1.17 and the choice of $q^{[t]}$. □

3.5 Let t be a variable-free $\mathcal{L}_T$ term. Then

$$\mathrm{VW}^+ \vdash \forall x,x_1 (\mathrm{Treemap}(x,\underline{\theta(t)}) \ \&\ \mathrm{Key}_m(x_1,\underline{\theta(t)}) \rightarrow$$

$$\rightarrow \forall y_1 \,|{\leq}|\, x_1 \forall y_2 \subseteq_p \underline{\theta(t)} (\exists z \subseteq_p x\ \mathrm{FrPair}(z,x_1,y_1,y_2) \leftrightarrow$$

$$\leftrightarrow \exists z \subseteq_p \underline{q}^{[t]}\ \mathrm{FrPair}(z,x_1,y_1,y_2))).$$

Proof: Fix t, and assume $M \vDash \mathrm{Treemap}(x,\underline{\theta(t)})$ & $\mathrm{Key}_m(x_1,\underline{\theta(t)})$ for some $x, x_1 \in M$. Let $x_1 = \underline{2}^m$. We first argue from right to left. Suppose $M \vDash \mathrm{FrPair}(z,x_1,y_1,y_2)$ for some $M \vDash z \subseteq_p \underline{q}^{[t]}$ where $y_1,y_2 \in M$. By 1.13 we have that $M \vDash \bigvee_{s \sqsubseteq q[t]} \bigvee_{r \sqsupseteq s} z = \underline{r}$. Hence $z \in \Sigma^*$ and, by 1.19, so are $x_1,y_1,y_2$. Then $\mathfrak{S}^+ \vDash \mathrm{FrPair}(\underline{z},\underline{x_1},\underline{y_1},\underline{y_2})$ by 1.17, and from 3.3 and 3.4 we have that $\mathfrak{S}^+ \vDash \mathrm{Key}_m(x_1,\underline{\theta(t)})$ & $y_1 = \underline{p}$ & $y_2 = \underline{\theta(s)}$ for some $<p,s> \in \mathcal{S}^{\sqsubseteq}[t]$. We now argue by induction on the generating relation of $\mathcal{S}^{\sqsubseteq}[t]$. If s is t, then $\mathfrak{S}^+ \vDash \underline{p} = \underline{2}$, and from hypothesis $M \vDash \text{Treemap}(x,\underline{\theta(t)})$ we have

$M \vDash \exists w \subseteq_p x\ \text{FrPair}(w,\underline{2}^m,\underline{p},\underline{\theta(t)})$, that is, $M \vDash \exists z \subseteq_p x\ \text{FrPair}(z,x_1,y_1,y_2)$ as claimed. Suppose $s = (s_1,s_2)$ and assume

$$M \vDash \exists z_1 \subseteq_p \underline{q}^{[t]}\ \text{FrPair}(z_1,x_1,\underline{p}\hat{\ }\underline{1},\underline{\theta(s_1)})\ \&\ \exists z_2 \subseteq_p \underline{q}^{[t]}\ \text{FrPair}(z_2,x_1,\underline{p}\hat{\ }\underline{2},\underline{\theta(s_2)}).$$

By 1.17 we have that

$$\mathfrak{S}^+ \vDash \exists z_1 \subseteq_p \underline{q}^{[t]}\ \text{FrPair}(z_1,\underline{2}^m,\underline{p}\hat{\ }\underline{1},\underline{\theta(s_1)})\ \&\ \exists z_2 \subseteq_p \underline{q}^{[t]}\ \text{FrPair}(z_2,\underline{2}^m,\underline{p}\hat{\ }\underline{2},\underline{\theta(s_2)}).$$

By choice of $\underline{q}^{[t]}$ it follows that $\mathfrak{S}^+ \vDash \exists z \subseteq_p \underline{q}^{[t]}\ \text{FrPair}(z,\underline{2}^m,\underline{p},\underline{\theta((s_1,s_2))})$, whence, again by 1.17(a), $M \vDash \exists z \subseteq_p \underline{q}^{[t]}\ \text{FrPair}(z,\underline{2}^m,\underline{p},\underline{\theta((s_1,s_2))})$. From the induction hypothesis it then follows that $M \vDash \exists z \subseteq_p x\ \text{FrPair}(z,\underline{2}^m,\underline{p},\underline{\theta((s_1,s_2))})$. Because $s = (s_1,s_2)$ we have that $\theta(s) \neq 1$. From the principal hypothesis $M \vDash \text{Treemap}(x,\underline{\theta(t)})$ we have that

$M \vDash \exists z_1,z_2 \subseteq_p x\ \exists u_1,u_2 \subseteq_p z_1\ \exists v_1,v_2 \subseteq_p z_2\ (\text{FrPair}(z_1,\underline{2}^m,u_1,u_2)\ \&$
$\&\ \text{FrPair}(z_2,\underline{2}^m,v_1,v_2)\ \&\ u_1 = \underline{p}\hat{\ }\underline{1}\ \&\ v_1 = \underline{p}\hat{\ }\underline{2}\ \&\ I(u_2)\ \&\ I(v_2)\ \&\ T^*(u_2,v_2,\underline{\theta(t)}))$

By 2.9 we have that $M \vDash T^*(\underline{\theta(s_1)}, \underline{\theta(s_2)},\underline{\theta((s_1,s_2))})$, and from 2.11 it follows that $M \vDash u_2 = \underline{\theta(s_1)}\ \&\ v_2 = \underline{\theta(s_2)}$. But then

$$M \vDash \text{FrPair}(z_1,x_1,\underline{p}\hat{\ }\underline{1},\underline{\theta(s_1)})\ \&\ \text{FrPair}(z_2,x_1,\underline{p}\hat{\ }\underline{2},\underline{\theta(s_2)})$$

with $M \vDash z_1 \subseteq_p x\ \&\ M \vDash z_2 \subseteq_p x$. This completes the induction step of the proof of the right to left part of 3.5.

For the left to right part, assume $M \vDash \exists z \subseteq_p x\ \text{FrPair}(z,x_1,y_1,y_2)$ for some $y_1,y_2 \in M$ such that $M \vDash y_1 \mathbin{|\leq|} x_1\ \&\ y_2 \subseteq_p \underline{\theta(t)}$. From $M \vDash y_1 \mathbin{|\leq|} \underline{2}^m$ we have by (VW6) that $M \vDash \bigvee_{lh(r)\leq m} y_1 = \underline{r}$. We now argue by induction on the length of strings $p \in \Sigma^*$ such that $2 \sqsubseteq p$ that if $M \vDash \exists z \subseteq_p x\ \text{FrPair}(z,\underline{2}^m,\underline{p},y_2)$, then

$M \vDash \exists w \subseteq_p \underline{q}^{[t]}\ FrPair(w,\underline{2}^m,\underline{p},y_2)$. If $p = 2$, then $M \vDash y_2 = \underline{\theta(t)}$, and we are done because $M \vDash \exists w \sqsupseteq \underline{q}^{[t]} FrPair(w,\underline{2}^m,\underline{2},\underline{\theta(t)})$ follows from 3.3. Suppose now $M \vDash \exists z \subseteq_p x\ FrPair(z,\underline{2}^m,\underline{p}\hat{\ }\underline{1},y_2)$. Then $M \vDash y_2 \neq y$ because $M \vDash \underline{p}\hat{\ }\underline{1} \neq \underline{2}$. From principal hypothesis $M \vDash Treemap(x,\underline{\theta(t)})$ we obtain

$$M \vDash FrPair(w,\underline{2}^m,w_1,w_2)\ \&\ FrPair(z',\underline{2}^m,w_3,w_4)$$

for some unique $w,w_1,w_2,z',w_3,w_4 \in M$ such that
$M \vDash w \subseteq_p x\ \&\ w_1\ |{\leq}|\ \underline{2}^m$ and either
$M \vDash \underline{p}\hat{\ }\underline{1} = w_1\hat{\ }\underline{1}\ \&\ w_3 = w_1\hat{\ }\underline{2}\ \&\ T^*(y_2,w_4,w_2)$, or
$M \vDash w_3 = w_1\hat{\ }\underline{1}\ \&\ \underline{p}\hat{\ }\underline{1} = w_1\hat{\ }\underline{2}\ \&\ T^*(w_4,y_2,w_2)$.

From $M \vDash w_1\ |{\leq}|\ \underline{2}^m$ we have by (VW6) that $M \vDash \bigvee_{lth(r)\leq m} w_1 = \underline{r}$. Hence $M \vDash \underline{p}\hat{\ }\underline{1} = w_1\hat{\ }\underline{2}$ is ruled out by 1.17and so the first disjunct must be true. But then $M \vDash w_1 = \underline{p}$ by 1.17. From the induction hypothesis it follows that $M \vDash \exists z \subseteq_p \underline{q}^{[t]}\ FrPair(z,\underline{2}^m,\underline{p},w_2)$, so $M \vDash w_2 \subseteq_p z$ and further $M \vDash w_2 \subseteq_p \underline{q}^{[t]}$ by 1.19(b) where $w_2 \in \Sigma^*$. From 3.3 we then have

$$M \vDash \exists z_1 \subseteq_p \underline{q}^{[t]}\ \exists z_2 \subseteq_p \underline{q}^{[t]}\ \exists u_1,u_2 \subseteq_p z_1\ \exists v_1,v_2 \subseteq_p z_2\ (FrPair(z_1,x_1,u_1,u_2)\ \&$$

$$\&\ FrPair(z_2,x_1,v_1,v_2)\ \&\ u_1 = w_1\hat{\ }\underline{1}\ \&\ v_1 = w_1\hat{\ }\underline{2}\ \&\ I(u_2)\ \&\ I(v_2)\ \&$$

$$\&\ T^*(u_2,v_2,w_2)\ \&\ u_2 \subseteq_p w_2\ \&\ v_2 \subseteq_p w_2).$$

By 3.4(b) we have that $M \vDash w_2 = \underline{\theta(s)}$ for some subterm s of t. Now, from $M \vDash I(u_2)\ \&\ I(v_2)\ \&\ T^*(y_2,w_4,w_2)$ it follows by 2.11 that $M \vDash u_2 = y_2\ \&\ v_2 = w_4$. Hence $M \vDash \exists w \subseteq_p \underline{q}^{[t]}\ FrPair(w,x_1,\underline{p}\hat{\ }\underline{1},y_2)$, as needed. Similarly, if $M \vDash \exists z \subseteq_p x\ FrPair(z,\underline{2}^m,\underline{p}\hat{\ }\underline{2},y_2)$. This completes the left-to-right part of the proof. □

We now proceed to set up the desired interpretation $I$ of $WT_{Rel}$ in W.

We let the $\mathcal{L}^+$ formula I(x) define the domain of the interpretation. If $b_i$, $i \in N$, is an individual constant of $\mathcal{L}_{T(Rel,c)}$ and $b_i$ is $\underline{s}$ for a variable-free $\mathcal{L}_T$-term s, let

$$[b_i]^I =: \theta(s).$$

We let $T^I(x,y,z) \equiv: T^+(x,y,z)$ and $x \sqsubseteq^I y \equiv: \exists z\, x \ll_z y$.

3.6 Let t be a variable-free $\mathcal{L}_T$-term. Then

$$VW^+ \vdash [\forall x\, (x \sqsubseteq \underline{t} \leftrightarrow \bigvee_{s \in \mathcal{S}[t]} x = \underline{s})]^I.$$

Proof: We need to show that

$$VW^+ \vdash \forall x\, (I(x) \to (\exists y\, x \ll_y \underline{\theta(t)} \leftrightarrow \bigvee_{s \sqsubseteq t} x = \theta(\underline{s}))).$$

Suppose $M \vDash I(x)$ where $x \in M$, and assume $M \vDash \exists y\, x \ll_y \underline{\theta(t)}$. Then

$M \vDash x \subseteq_p \underline{\theta(t)}\ \&\ \mathrm{Treemap}(y,\underline{\theta(t)})\ \&\ \exists x_1 \subseteq_p y\, (\mathrm{Key}_m(x_1,\underline{\theta(t)})\ \&$

$\&\ \exists z_1,z_2 \subseteq_p y\, \exists v\, |{\le}|\, x_1\, (\mathrm{FrPair}(z_1,x_1,v,x)\ \&\ \mathrm{FrPair}(z_2,x_1,\underline{2},\underline{\theta(t)})\ \&\ \underline{2} \sqsubseteq v))$

for some $y \in M$. Let $M \vDash \mathrm{Key}_m(\underline{2}^n,\underline{\theta(t)})$. Then by 3.5,
$M \vDash \exists z_1 \subseteq_p \underline{q}^{[t]}\, \mathrm{FrPair}(z_1,\underline{2}^n,v,x)$, for some $v,x \in M$, and by 1.13 we have

$$M \vDash \bigvee_{lh(p)\le m} \exists z \subseteq_p \underline{q}^{[t]}(x = \underline{u}\ \&\ \mathrm{FrPair}(z,\underline{2}^n,\underline{p},\underline{u}))$$

for some $u \in \Sigma^*$. By 3.4(b) it follows that $u = \theta(s)$ for some $<p,s> \in \mathcal{S}^{\sqsubseteq}[t]$. But then $M \vDash \bigvee_{s \sqsubseteq t} x = \underline{\theta(s)}$.

Conversely, assume $M \vDash \bigvee_{s \sqsubseteq t} x = \underline{\theta(s)}$ where $x \in M$. Then $<p,s> \in \mathcal{S}^{\sqsubseteq}[t]$ for some dyadic string p such that $2 \sqsubseteq p$ and $p\, |{\le}|\, 2^n$. By 3.4(a) and 1.17 it follows that $M \vDash \exists z \subseteq_p \underline{q}^{[t]}\ \exists x_1 \subseteq_p z\ (FrPair(z,x_1,\underline{p},\underline{\theta(s)})\ \&\ Key_m(x_1,\underline{\theta(t)}))$ and $M \vDash \exists w \subseteq_p \underline{q}^{[t]}\ \exists x_1 \subseteq_p w\ (FrPair(w,x_1,\underline{2},\underline{\theta(t)})\ \&\ Key_m(x_1,\underline{\theta(t)})$. But

$\mathfrak{S}^+ \vDash \underline{\theta(s)} \subseteq_p \underline{\theta(t)}$ so by 1.14(a) we have $M \vDash \underline{\theta(s)} \subseteq_p \underline{\theta(t)}$. From all this it follows that $M \vDash \underline{\theta(s)} \ll_{q[t]} \underline{\theta(t)}$, hence $M \vDash x \ll_{q[t]} \underline{\theta(t)}$, that is, $M \vDash x \sqsubseteq^I \underline{t}^I$, as needed. □

Theorem 3.7 $WT_{Rel}$ is interpretable in $VW^+$.

Proof: That $VW^+ \vdash I(\underline{\theta(t)})$ for each $\mathcal{L}_T$ term t was proved in 2.4, hence $VW^+ \vdash I([b_i]^I)$ for each constant $b_i$ of $WT_{Rel}$. That $VW^+ \vdash [WT_{Rel}1]^I$ follows from (VW1) and the fact that the map θ is 1-1. That $VW^+ \vdash [WT_{Rel}2]^I$ was proved in 3.6. For any variable-free $\mathcal{L}_T$ term s, t we have that $VW^+ \vdash T^+((\underline{\theta(s)}, \underline{\theta(t)},\underline{\theta((s,t))}))$, hence also $VW^+ \vdash [WT_{Rel}3]^I$. That $VW^+ \vdash [WT_{Rel}4]^I$ and $VW^+ \vdash [WT_{Rel}5]^I$ follows from 2.5 and 2.6. □

## §4. Defining Binary Concatenation for Tallies and Pseudotallies

Let $Dpl(x) \equiv:\ x = \underline{1}\hat{}\underline{1}\ \vee\ x = \underline{2}\hat{}\underline{2}$, ("x is a duplicate"),

$PsTally(x) \equiv:\ x \neq \underline{e}\ \&\ x \neq \underline{1}\ \&\ x \neq \underline{2}\ \&\ \forall y \subseteq_p x\ \neg Dpl(y)$,

(x is a pseudotally")

$PsTally_{21}(x) \equiv:\ PsTally(x)\ \&\ \underline{2}\hat{}\underline{1} \sqsubseteq x\ \&\ \underline{2}\hat{}\underline{1} \sqsupseteq x$, (x is a 2,1-pseudotally").

We have that:

4.1 For each $s \in \Sigma^*$, $\mathfrak{S}^+ \vDash PsTally_{21}(\underline{s}) \Leftrightarrow s = (\underline{2}\,\hat{}\,\underline{1})^n$ for some $n \geq 1$.

Proof: First note that every 2,1-pseudotally in $\Sigma^*$ must have even length. For, every 2,1-pseudotally must have length at least 2. Suppose, for a reductio, that a 2,1-pseudotally of odd length exists, and let $u \in \Sigma^*$ be a shortest such. Then because $\underline{2}\,\hat{}\,\underline{1} \sqsupseteq u$, there must be a $u_0 \in \Sigma^*$, $u_0 \neq e$, such that $u = u_0\,\hat{}\,2\,\hat{}\,1$. Then $u_0$ is of odd length, and it cannot be a 2,1-pseudotally because $u_0$ is shorter than u. But then $u_0 = 1$, or $u_0 = 2$, or $u_0$ contains a subsegment that is a duplicate. The last possibility is ruled out because otherwise u would contain a duplicate as a subsegment, contradicting the hypothesis. And $u_0 = 1$ is ruled out because u, and hence, $u_0$ must have 2 as its initial digit. But we cannot have $u_0 = 2$ because then u would begin with $\underline{2}\,\hat{}\,\underline{2}$, contradicting the hypothesis that u is a 2,1-pseudotally. Hence no such $u \in \Sigma^*$ exists.

Now let $s \in \Sigma^*$, and let $lh(s) = 2n$, for $n \geq 1$. Suppose $\mathfrak{S}^+ \models PsTally_{21}(\underline{s})$. If $n = 1$, then s is $(\underline{2}\,\hat{}\,\underline{1}) = (\underline{2}\,\hat{}\,\underline{1})^1$. Assume as the induction hypothesis that the claim holds for n, and suppose s has length $2(n + 1) = 2n + 2$, with $n \geq 1$. Then $s = s_0\,\hat{}\,2\,\hat{}\,1$ where $\mathfrak{S}^+ \models PsTally_{21}(\underline{s_0})$. From the induction hypothesis then $s_0 = (\underline{2}\,\hat{}\,\underline{1})^k$ for some $k \geq 1$. But then $s = s_0\,\hat{}\,2\,\hat{}\,1 = (\underline{2}\,\hat{}\,\underline{1})^{k+1}$.

Conversely, suppose $s = (\underline{2}\,\hat{}\,\underline{1})^k$ for some $k \geq 1$. We argue by induction on k. Suppose, for a reductio, that $\mathfrak{S}^+ \models \underline{d} \subseteq_p \underline{s}\ \&\ Dpl(\underline{d})$. Then we must have $k > 1$, so $s = (\underline{2}\,\hat{}\,\underline{1})^j\,\hat{}\,(\underline{2}\,\hat{}\,\underline{1})$. Then $\mathfrak{S}^+ \models \underline{d} \subseteq_p (\underline{2}\,\hat{}\,\underline{1})^j$, or $\mathfrak{S} \models \underline{d} \subseteq_p \underline{2}\,\hat{}\,\underline{1}$, or $(\underline{2}\,\hat{}\,\underline{1})^j$ ends with the first digit of duplicate d and $\underline{2}\,\hat{}\,\underline{1}$ begins with the second digit of d. The first two scenarios are ruled out by the induction hypothesis since then both $(\underline{2}\,\hat{}\,\underline{1})^j$ and $\underline{2}\,\hat{}\,\underline{1}$ are 2,1-pseudotallies. The fact that the last digit of $(\underline{2}\,\hat{}\,\underline{1})^{j+1}$ is $\underline{1}$ and the first digit of $\underline{2}\,\hat{}\,\underline{1}$ is $\underline{2}$ rules out the third scenario. Hence $s = (\underline{2}\,\hat{}\,\underline{1})^{j+1}$ must also be a 2,1-pseudotally. □

For brevity, we sometimes write '21' for '$\underline{2}$ ˆ $\underline{1}$'.

Let $t_1,t_2$ be 2-tallies. We say that $s_1,\ldots,s_n$ is a <u>tally concatenation sequence for $t_1,t_2$</u> if each $s_i$, $1 \le i \le n$, is an ordered pair of strings $u_i$, $v_i \in \Sigma^*$, such that $u_1 = 2$ and $v_1 = t_1 \hat{}\, 2$, and $u_i$ is a proper initial segment of $t_2$ with $u_{i+1} = u_i \hat{}\, 2$ and $v_{i+1} = v_i \hat{}\, 2$, for each i, $1 \le i < n$. Then the last term, $s_n$, of the sequence is the pair $u_n$, $v_n$ where $u_n = t_2$ and $v_n = t_1 \hat{}\, t_2$. Hence for any 2-tallies x, y and a string $z \in \Sigma^*$, $z = x \hat{}\, y$ if and only if z is the second term of the last member of a (unique) tally concatenation sequence for x, y. We use this idea to express in $VW^+$ binary concatenation for tallies.

Let

$$TallyPair(x,y,z) \equiv: Tally_2(y) \ \&\ Tally_2(z) \ \&\ y\hat{}21 \sqsubseteq x \ \&\ 21_{\vee}z \sqsupseteq x \ \&\ 21 \subseteq_1 x,$$

reading "x is a tally-pair of y and z". Then we have:

4.2(a) For any $t_1,t_2 \in \Sigma^*$,

$$\mathfrak{S}^+ \vDash Tally_2(\underline{t}_1) \ \&\ Tally_2(\underline{t}_2) \rightarrow \exists x \subseteq_p (\underline{t}_1\hat{}21)\hat{}\underline{t}_2\ TallyPair(x,\underline{t}_1,\underline{t}_2).$$

(b) For any $s \in \Sigma^*$,

$$\mathfrak{S}^+ \vDash \forall y_1,y_2,z_1,z_2 \subseteq_p \underline{s}\ (TallyPair(\underline{s},y_1,z_1) \ \&\ TallyPair(\underline{s},y_2,z_2) \rightarrow y_1 = y_2 \ \&\ z_1 = z_2).$$

We let TPS(x) abbreviate the $\mathcal{L}_W$-formula

$$(21)^2 \sqsubseteq x \ \&\ (21)^2 \sqsupseteq x \ \&\ \neg(21)^3 \subseteq_p x \ \&$$

$$\&\ \forall y \subseteq_p x\ ((21)^2 \sqsubseteq y \ \&\ (21)^2 \sqsupseteq y \ \&\ \neg(\underline{2}_{\vee}(21)^2)\hat{}\underline{2} \subseteq_p y \rightarrow$$

$$\rightarrow \exists u \subseteq_p y\ \exists v,w \subseteq_p u\ (TallyPair(u,v,w) \ \&\ y = ((21)^2{}_{\vee}u)\hat{}(21)^2))),$$

reading "x is a tally-pair sequence".

We then let $TCSeq(x,t_1,t_2)$ abbreviate the conjunction

$$TPS(x) \ \& \ Tally_2(t_1) \ \& \ Tally_2(t_2) \ \& \ (s1) - (s7),$$

with

(s1) $\exists u \subseteq_p x\, (TallyPair(u,\underline{2},t_1{}^\frown\underline{2}) \ \& \ ((21)^2{}_\vee u)^\frown(21)^2 \sqsubseteq x)$,

(s2) $\exists u \subseteq_p x\, \exists z \subseteq_p u\, (TallyPair(u,t_2,z) \ \& \ ((21)^2{}_\vee u)^\frown(21)^2 \sqsupseteq x \ \&$

$\& \ ((21)_\vee z)^\frown(21)^2 \sqsupseteq x)$,

(s3) $\forall w_1 \sqsubseteq x\, \forall u \subseteq_p x\, \forall t \sqsubseteq t_2\, \forall y \subseteq_p u\, (u^\frown(21)^2 \sqsupseteq w_1 \ \& \ t \neq t_2 \ \& \ t \neq \underline{e} \ \&$

$\& \ TallyPair(u,t,y) \ \rightarrow \ \exists w_2 \sqsubseteq x\, \exists v \subseteq_p x\, (w_1 \sqsubseteq w_2 \ \& \ v^\frown(21)^2 \sqsupseteq w_2 \ \&$

$\& \ TallyPair(v,t^\frown\underline{2},y^\frown\underline{2}) \ \& \ \neg v \subseteq_p w_1 \ \&$

$\& \ \forall v' \subseteq_p x\, \forall t',y' \subseteq_p v'\, (TallyPair(v',t',y') \ \& \ \neg v' \subseteq_p w_1 \ \& \ v' \subseteq_p w_2 \ \rightarrow \ v' = v)))$

(s4) $\forall t \sqsubseteq t_2\, \forall u,v \subseteq_p x\, \forall y,z \subseteq_p x (TallyPair(u,t,y) \ \& \ TallyPair(v,t,z) \ \rightarrow$

$\rightarrow \ u = v \ \& \ y = z)$,

(s5) $\forall u \subseteq_p x\, \forall t,y \subseteq_p u\, (TallyPair(u,t,y) \rightarrow \ t \neq \underline{e} \ \& \ t \sqsubseteq t_2 \ \& \ u \subseteq_1 x)$,

(s6) $\forall t \sqsubseteq t_2\, (t \neq \underline{e} \ \rightarrow \exists u \subseteq_p x\, \exists y \subseteq_p u\, TallyPair(u,t,y))$,

(s7) $\forall t,t' \sqsubseteq t_2\, \forall w,w' \sqsubseteq x\, \forall u,v,y,y' \subseteq_p x (t \neq \underline{e} \ \& \ t \sqsubseteq t' \ \& \ t \neq t' \ \& \ TallyPair(u,t,y) \ \&$

$\& \ TallyPair(v,t',y') \ \& \ u^\frown(21)^2 \sqsupseteq w \ \& \ v^\frown(21)^2 \sqsupseteq w' \ \rightarrow \ w \sqsubseteq w')$.

We now have:

4.3 For any $s,t_1,t_2 \in \Sigma^*$,

$\mathfrak{S}^+ \vDash TCSeq(\underline{s},\underline{t}_1,\underline{t}_2) \Leftrightarrow$ $t_1$ and $t_2$ are 2-tallies, $lh(t_2) = m$ for some $m \geq 1$ and

$s = (21)^2 \hat{}\ s_1 \hat{}\ (21)^2 \hat{} \ldots \hat{}\ (21)^2 \hat{}\ s_m \hat{}\ (21)^2$ where $s_1,\ldots,s_m$ is a (unique) tally concatenation sequence for $t_1,t_2$.

Proof: From right to left, let $t_1,t_2$ be 2-tallies, and let $lh(t_2) = m \geq 1$. We argue by induction on m. Assume m = 1 and let $s = (21)^2 \hat{}\ s_1 \hat{}\ (21)^2$ where $s_1$ is a single-term tally concatenation sequence for $t_1,t_2$. Then $s_1 = (u_1 \hat{}\ (21)) \hat{}\ v_1$ where $u_1 = 2$ and $v_1 = t_1 \hat{}\ 2$. Letting x = s, we verify that $\mathfrak{S}^+ \vDash TCSeq(x,t_1,\underline{2})$. Hence $\mathfrak{S}^+ \vDash TCSeq(\underline{s},\underline{t}_1,\underline{t}_2)$. Assume that the claim holds for m = k, and suppose

$$s = (21)^2 \hat{}\ s_1 \hat{}\ (21)^2 \hat{} \ldots \hat{}\ (21)^2 \hat{}\ s_k \hat{}\ (21)^2 \hat{}\ s_{k+1} \hat{}\ (21)^2$$

where $s_1,\ldots,s_k,s_{k+1}$ is a (unique) tally concatenation sequence for $t_1,t_2$ and $lh(t_2) = k+1$. Letting $t_2 = t_3 \hat{}\ 2$, we have that $s_{k+1} = (u_{k+1} \hat{}\ (21)) \hat{}\ v_{k+1}$ where $u_{k+1} = t_3 \hat{}\ 2$ and $v_{k+1} = t_1 \hat{}\ (t_3 \hat{}\ 2)$. Then $s_1,\ldots,s_k$ is a (unique) tally concatenation sequence for $t_1,t_3$, with $s_k = (u_k \hat{}\ (21)) \hat{}\ v_k$ where $u_k = t_3$ and $v_k = t_1 \hat{}\ t_3$. From the induction hypothesis it follows that $\mathfrak{S}^+ \vDash TCSeq(\underline{s^*},\underline{t}_1,\underline{t}_3)$ where

$s^* = (21)^2 \hat{}\ s_1 \hat{}\ (21)^2 \hat{} \ldots \hat{}\ (21)^2 \hat{}\ s_k \hat{}\ (21)^2$. We then proceed to verify that indeed $\mathfrak{S} \vDash TCSeq(((\underline{s^*}\hat{}\underline{s}_{k+1})\hat{}(21)^2,\underline{t}_1,\underline{t}_2)$.

For the left-to-right part, assume $\mathfrak{S} \vDash TCSeq(\underline{s},\underline{t}_1,\underline{t}_2)$. Then $\mathfrak{S}^+ \vDash Tally_2(\underline{t}_1)\ \&\ Tally_2(\underline{t}_2)$ We argue by induction on the number $k \geq 2$ of occurrences of the 2,1-pseudotally $(21)^2$ in $\underline{s}$. Suppose k = 2. From $\mathfrak{S}^+ \vDash TPS(\underline{s})$ we have that $\mathfrak{S}^+ \vDash (21)^2 \sqsubseteq \underline{s}\ \&\ (21)^2 \sqsupseteq \underline{s}$, hence

$\mathfrak{S}^+ \vDash \exists u\ (u \neq \underline{e}\ \&\ \underline{s} = ((21)^2{}_\vee u)\hat{}(21)^2)$. Then $\mathfrak{S}^+ \vDash \neg(21)^2 \sqsubseteq_p u$ and, further, $\mathfrak{S}^+ \vDash \neg(\underline{2}_\vee(21)^2)\hat{}\underline{2} \sqsubseteq_p \underline{s}$. Hence, again from $\mathfrak{S}^+ \vDash TPS(\underline{s})$, it follows that $\mathfrak{S}^+ \vDash TallyPair(u_0,v,w)\ \&\ \underline{s} = ((21)^2{}_\vee u_0)\hat{}(21)^2$ for some $u_0,v,w \in \Sigma^*$. Then $\mathfrak{S}^+ \vDash u = u_0$. From hypothesis $\mathfrak{S}^+ \vDash TCSeq(\underline{s},\underline{t}_1,\underline{t}_2)$ by (s1) we have that $\mathfrak{S}^+ \vDash ((21)^2{}_\vee u')^\wedge(21)^2 \sqsubseteq \underline{s}$ for some $u' \in \Sigma^*$ such that $\mathfrak{S}^+ \vDash \text{TallyPair}(u',\underline{2},\underline{t}_1{}^\wedge\underline{2})$. Then $\mathfrak{S}^+ \vDash u' = u_0 = u$ & $v = \underline{2}$ & $w = \underline{t}_1{}^\wedge\underline{2}$ by 4.2(b). We claim that $\mathfrak{S}^+ \vDash \underline{t}_2 = \underline{2}$. Suppose, for a reductio, that $\mathfrak{S}^+ \vDash \underline{t}_2 \neq \underline{2}$. Then $\mathfrak{S}^+ \vDash \text{Tally}_2(\underline{t}_2)$ & $\underline{t}_2 \neq \underline{e}$. By (s6) then $\mathfrak{S}^+ \vDash \text{TallyPair}(u_1,t_2,y)$ & $u_1 \subseteq_p \underline{s}$ for some $u_1,y \in \Sigma^*$. But then $\mathfrak{S}^+ \vDash u_1 = u_0$ & $\underline{t}_2 = v = \underline{2}$, contradicting hypothesis. Hence indeed $\mathfrak{S}^+ \vDash \underline{t}_2 = \underline{2}$, and so $lh(t_2) = 1$ and $u = u_0$ is a single-term tally concatenation sequence for $t_1,t_2$, as claimed.

Suppose now that the claim holds for $k = j$ and that s has $j+1$ occurrences of $(21)^2$. From hypothesis $\mathfrak{S}^+ \vDash \text{TCSeq}(\underline{s},\underline{t}_1,\underline{t}_2)$ we have by (s2) that $\mathfrak{S}^+ \vDash ((21)^2{}_\vee u)^\wedge(21)^2 \sqsupseteq \underline{s}$ & $\text{TallyPair}(u,t_2,z)$ for some $u, z \in \Sigma^*$. Then $\mathfrak{S}^+ \vDash \neg(21)^2 \subseteq_p u$, and since $j+1 > k \geq 2$, we have that $\mathfrak{S}^+ \vDash \underline{s} \neq ((21)^2{}_\vee u)^\wedge(21)^2$ and $s = s^{*\frown} ((21)^2{}_\vee u)^\wedge(21)^2$ for some $s^* \in \Sigma^*$, $s^* \neq e$. By (s1) then $\mathfrak{S}^+ \vDash ((21)^{2\frown} u_0)^\frown (21)^2 \sqsubseteq \underline{s}$, and by (s5) we must have $\mathfrak{S}^+ \vDash u_0 \neq u$. But then from (s4) it follows that $\mathfrak{S}^+ \vDash \underline{t}_2 \neq \underline{2}$, that is, we must have $t_2 = t_3{}^\frown \underline{2}$ for some $t_3 \neq e$. We claim that

(*) $\qquad \mathfrak{S}^+ \vDash \text{TCSeq}(\underline{s^{**}},\underline{t}_1,\underline{t}_3)$

where $s^{**} = s^{*\frown} (21)^2$. Now, $\mathfrak{S}^+ \vDash \text{TPS}(\underline{s}^{**})$ follows from $\mathfrak{S}^+ \vDash \text{TPS}(\underline{s})$. We have that $s^{**}$ has $j \geq 2$ occurrences of $(21)^2$. Then

$$\mathfrak{S}^+ \vDash \exists y \sqsupseteq s^{**} ((21)^2 \sqsubseteq y \ \&\ (21)^2 \sqsupseteq y \ \&\ \neg(\underline{2}_\vee(21)^2)^\wedge\underline{2} \subseteq_p y).$$

From $\mathfrak{S}^+ \vDash \text{TPS}(\underline{s}^{**})$ it follows that $\mathfrak{S}^+ \vDash y = ((21)^2{}_\vee u)^\wedge(21)^2$ for some $u,v,w \in \Sigma^*$ such that $\mathfrak{S}^+ \vDash \text{TallyPair}(u,v,w)$. From (s3), (s5) and (s6) in hypothesis $\mathfrak{S}^+ \vDash \text{TCSeq}(\underline{s},\underline{t}_1,\underline{t}_2)$ it follows that $v = t_3$, hence condition (s2) in $\text{TCSeq}(\underline{s^{**}},\underline{t}_1,\underline{t}_3)$ holds. Condition (s3) for $\text{TCSeq}(\underline{s^{**}},\underline{t}_1,\underline{t}_3)$ is immediate from hypothesis $\mathfrak{S}^+ \vDash \text{TCSeq}(\underline{s},\underline{t}_1,\underline{t}_2)$, and conditions (s4)-(s7) for $\text{TCSeq}(\underline{s^{**}},\underline{t}_1,\underline{t}_3)$ likewise follow. This completes the argument for claim (*).

Now, since $s^{**}$ has j occurrences of $(21)^2$, from the induction hypothesis we have that $lh(t_3) = m$ and $s^{**} = (21)^2\hat{}\, s_1\hat{}\,(21)^2\hat{}\ldots\hat{}\,(21)^2\hat{}\, s_m\hat{}\,(21)^2$ where $s_1,\ldots,s_m$ is a (unique) tally concatenation sequence for $t_1,t_3$. Then $s_m$ is the ordered pair of $u_m = t_3$ and $v_m= t_1\hat{}\,t_3$. But $s = s^{**}\hat{}\,u\hat{}\,(21)^2$ where $\mathfrak{S}^+ \vDash TallyPair(u,t_2,z)$ for some $z \in \Sigma^*$. Now, from $\mathfrak{S}^+ \vDash TallyPair(\underline{s}_m,\underline{t}_3,\underline{t}_1\hat{}\underline{t}_3)$ we have by condition (s3) in hypothesis $\mathfrak{S}^+ \vDash TCSeq(\underline{s},\underline{t}_1,\underline{t}_2)$ that

$$\mathfrak{S}^+\vDash \exists w \sqsubseteq \underline{s}(TallyPair(v,\underline{t}_3\hat{}\underline{2},(\underline{t}_1\hat{}\underline{t}_3)\hat{}\underline{2}) \;\&\; \underline{s^{**}} \sqsubseteq w \;\&\; v\hat{}(21)^2 \sqsupseteq w \;\&\; \neg v \subseteq_p \underline{s^{**}})$$

for some unique $v \in \Sigma^*$. But then $\mathfrak{S}^+ \vDash v = u \;\&\; z = (\underline{t}_1\hat{}\underline{t}_3)\hat{}\underline{2} = \underline{t}_1\hat{}\underline{t}_2$. So s is a tally concatenation sequence for $t_1,t_2$, as claimed. □

We now let $\quad T_2C(x,y,z) \;\equiv:\; \exists u\, (TCSeq(u,x,y) \;\&\; ((21)^2{}_\vee z)\hat{}(21)^2 \sqsupseteq u).$

4.4 Let $t_1,t_2 \in \Sigma^*$ be any 2-tallies. Then (a) $VW^+ \vdash T_2C(\underline{t}_1,\underline{t}_2,\underline{t_1\hat{}\,t_2})$.

(b) $VW^+ \vdash \forall z_1,z_2\, (T_2C(\underline{t}_1,\underline{t}_2,z_1) \;\&\; T_2C(\underline{t}_1,\underline{t}_2,z_2) \rightarrow z_1 = z_2)$.

Proof: Let $s \in \Sigma^*$ be a (unique) tally concatenation sequence for $t_1,t_2$ with the ordered pair $t_2$, $t_1\hat{}t_2$ as its last term. By 4.3 we have that $\mathfrak{S}^+ \vDash TCSeq(\underline{s},\underline{t}_1,\underline{t}_2)$, and by (s2) we have that

$$\mathfrak{S}^+ \vDash \exists u \subseteq_p \underline{s}\; \exists z \subseteq_p u(TallyPair(u,\underline{t}_2,z) \;\&\; ((21)^2{}_\vee u)\hat{}(21)^2 \sqsupseteq \underline{s} \;\&$$

$$\&\; ((21)_\vee z)\hat{}(21)^2 \sqsupseteq \underline{s}).$$

From (s3)-(s7) it follows that the successive initial segments $t_1\hat{}\,t$ of $t_1\hat{}\,t_2$ map 1-1 in an ordered preserving way onto the successive initial segments of s ending with tally-pairs of t and $t_1\hat{}\,t$. Then $\mathfrak{S}^+ \vDash z = \underline{t}_1\hat{}\underline{t}_2$ and we have that $\mathfrak{S}^+ \vDash \text{TCSeq}(\underline{s},\underline{t}_1,\underline{t}_2)\ \&\ ((21)_\vee(\underline{t}_1\verb|^|\underline{t}_2))\verb|^|(21)^2 \sqsupseteq \underline{s},$

By 1.3 and 1.17 it follows that

$\text{VW}^+ \vdash \text{TCSeq}(\underline{s},\underline{t}_1,\underline{t}_2)\ \&\ ((21)_\vee(\underline{t_1\hat{}\,t_2}))\verb|^|(21)^2 \sqsupseteq \underline{s}.$

But then $\text{VW}^+ \vdash T_2C(\underline{t}_1,\underline{t}_2,\underline{t_1\hat{}\,t_2})$. (b) follows from the above proof. □

We let L(x,y) abbreviate the $\mathcal{L}^+$ formula

$(x = \underline{e}\ \&\ y = \underline{e})\ v\ (x \neq \underline{e}\ \&\ \text{Tally}_2(y)\ \&\ x\ |{\leq}|\ y\ \&$

$\&\ \forall u \sqsubseteq y\ (u \neq y\ \&\ u \neq e\ \rightarrow\ \neg x\ |{\leq}|\ u)).$

Then we have:

4.5 For each s, t ∈ Σ*, $\mathfrak{S}^+ \vDash L(\underline{s},\underline{t}) \Leftrightarrow$ for some $m \geq 0$, $t = 2^m$ and $lh(s) = m$.

<u>Proof:</u> Assume $\mathfrak{S}^+ \vDash L(\underline{s},\underline{t})$. If s = e and t = e, then $t = 2^0$ and $lh(s) = 0$, as claimed. Suppose $t = 2^m$ for $m > 0$. Then $lh(s) \leq m$ by hypothesis. Also from hypothesis, $lh(s) \nleq k$ for $2^k \sqsubseteq 2^m = t$ for each $k < m$. So $lh(s) = m$. The converse is immediate. □

4.6 For each m, n ∈ N,

$\text{VW}^+ \vdash \forall x\ |{\leq}|\ 2^m\ \forall y\ |{\leq}|\ 2^n\ \forall z_1,z_2\ |{\leq}|\ 2^{m+n+2}\ (L(x,2^m)\ \&\ L(y,2^n) \rightarrow$

$\rightarrow (x\verb|^|21 \sqsubseteq z_1\ \&\ x\verb|^|21 \sqsubseteq z_2\ \&\ 21_\vee y \sqsupseteq z_1\ \&\ 21_\vee y \sqsupseteq z_2\ \rightarrow z_1 = z_2)).$

<u>Proof:</u> We reason in $\mathfrak{S}^+$. Let $q,r,u_1,u_2 \in \Sigma^*$ be such that $\mathfrak{S}^+ \vDash L(\underline{q},2^m)\ \&\ L(\underline{r},2^n)$ along with $\mathfrak{S}^+ \vDash \underline{q}\verb|^|21 \sqsubseteq \underline{u}_1\ \&\ \underline{q}\verb|^|21 \sqsubseteq \underline{u}_2\ \&\ 21_\vee\underline{r} \sqsupseteq u_1\ \&\ 21_\vee\underline{r} \sqsupseteq u_2$ where

$\mathfrak{S}^+ \vDash \underline{u}_1 |\leq| 2^{m+n+2}$ & $\underline{u}_2 |\leq| 2^{m+n+2}$. From 4.5 and the hypothesis we have that lh(q) = m and lh(r ) = n, and further that, if $1 \leq i \leq m+2$ or

$lh(u_1) - (n+2) \leq i \leq lh(u_1)$ or $lh(u_2) - (n+2) \leq i \leq lh(u_2)$ then $u_{1,i} = u_{2,i}$, where $u_{1,i}$ and $u_{2,I}$ are the i-th digits of $u_1, u_2$, resp.. Assume, for a reductio, that $u_1 \neq u_2$. Then there is a smallest $j \leq m+n+2$ such that $u_{1,j} \neq u_{2,j}$. But then $m+2 < j$ and $j < lh(u_1) - (n+2) \leq m+n+2 - (n+2) = m$, which is impossible. Hence $u_1 = u_2$, and the desired claim follows by 1.17(a). □

Let $p_1,p_2$ be 2,1-pseudotallies. We say that $q_1,\ldots,q_n$ is a <u>pseudotally concatenation sequence for $p_1,p_2$</u> if each $q_i$, $1 \leq i \leq n$, is an ordered pair of strings $u_i, v_i \in \Sigma^*$, such that $u_1 = 21$ and $v_1 = p_1 ^\frown 21$, and $u_i$ is a proper initial segment of $p_2$ with $u_{i+1} = u_i ^\frown 21$ and $v_{i+1} = v_i ^\frown 21$, for each i, $1 \leq i < n$. Then the last term, $q_n$, of the sequence is the pair $u_n, v_n$ where $u_n = p_2$ and $v_n = p_1 ^\frown p_2$. Hence for any 2,1-pseudotallies x, y and a string $z \in \Sigma^*$, $z = x ^\frown y$ if and only if z is the second term of the last member of a (unique) pseudotally concatenation sequence for x, y.

We also introduce the "dual" pairing of 2,1-pseudotallies: let PsTallyPair(x,y,z) abbreviate

$$PsTally_{21}(y) \ \& \ PsTally_{21}(z) \ \& \ y\hat{}2^2 \sqsubseteq x \ \& \ 2_{\vee}z \sqsupseteq x \ \& \ 2^2 \subseteq_1 x,$$

reading "x is a pseudotally-pair of y and z". Then we have:

4.7(a) For any $p_1,p_2 \in \Sigma^*$,

$\mathfrak{S}^+ \vDash PsTally_{21}(\underline{p}_1)$ & $PsTally_{21}(\underline{p}_2) \rightarrow$

$\rightarrow \exists x \subseteq_p (\underline{p}_1\hat{}2^2)\hat{}\underline{p}_2 \ \mathrm{PsTallyPair}_{21}(x,\underline{p}_1,\underline{p}_2).$

(b) For any $s \in \Sigma^*$,

$\mathfrak{S}^+ \vDash \forall y_1,y_2,z_1,z_2 \subseteq_p \underline{s} \ (\mathrm{PsTallyPair}_{21}(\underline{s},y_1,z_1) \ \& \ \mathrm{PsTallyPair}_{21}(\underline{s},y_2,z_2) \rightarrow$

$\rightarrow y_1 = y_2 \ \& \ z_1 = z_2).$

Proof: For (b), we have from $\mathfrak{S}^+ \vDash \mathrm{PsTally}_{21}(z_1)$ by 4.1 that $\mathfrak{S}^+ \vDash 21 \sqsubseteq z_1$, so $\mathfrak{S}^+ \vDash \exists z_3 \subseteq_p z_1 \ 21\hat{}z_3 = z_1$. Similarly, $\mathfrak{S}^+ \vDash \exists z_4 \subseteq_p z_2 \ 21\hat{}z_4 = z_2$. From hypothesis we have that $\mathfrak{S}^+ \vDash \underline{2}\hat{}(21\hat{}z_3) \sqsupseteq \underline{s} \ \& \ \underline{2}\hat{}(21\hat{}z_4) \sqsupseteq \underline{s}$, hence $\mathfrak{S}^+ \vDash (\underline{2}\hat{}\underline{2})\hat{}(\underline{1}\hat{}z_3) \sqsupseteq \underline{s} \ \& \ (\underline{2}\hat{}\underline{2})\hat{}(\underline{1}\hat{}z_4) \sqsupseteq \underline{s}$. But then from $2^2 \subseteq_1 s$ it follows that $\mathfrak{S}^+ \vDash \underline{1}\hat{}z_3 = \underline{1}\hat{}z_4$, hence $\mathfrak{S}^+ \vDash z_3 = z_4$. Then

$$\mathfrak{S}^+ \vDash z_1 = 21\hat{}z_3 = \underline{2}\hat{}(\underline{1}\hat{}z_3) = \underline{2}\hat{}(\underline{1}\hat{}z_4) = 21\hat{}z_4 = z_2.$$

On the other hand, from $\mathfrak{S}^+ \vDash y_1\hat{}2^2 \sqsubseteq s \ \& \ y_2\hat{}2^2 \sqsubseteq s \ \& \ 2^2 \subseteq_1 s$ we have $\mathfrak{S}^+ \vDash y_1 = y_2$. □

We let PsTPS(x) abbreviate the $\mathcal{L}^+$ formula

$2^3 \sqsubseteq x \ \& \ 2^2 \sqsupseteq x \ \& \ \neg 2^4 \subseteq_p x \ \&$

$\& \ \forall y \subseteq_p x \ (2^3 \sqsubseteq y \ \& \ 2^2 \sqsupseteq y \ \& \ \neg(\underline{1}_\vee 2^3)\hat{}\underline{1} \subseteq_p y \ \rightarrow$

$\rightarrow \exists u \subseteq_p y \ \exists v,w \subseteq_p u \ (\mathrm{PsTallyPair}_{21}(u,v,w) \ \& \ y = (2^2{}_\vee u)\hat{}2^2)),$

reading "x is a pseudotally-pair sequence".

Let $\mathrm{PsTCSeq}(x,p_1,p_2)$ abbreviate the conjunction

$$\mathrm{PsTPS}(x) \ \& \ \mathrm{PsTally}_{21}(\underline{p}_1) \ \& \ \mathrm{PsTally}_{21}(\underline{p}_2) \ \& \ (\mathrm{ps1}) - (\mathrm{ps7}),$$

with

(ps1) $\exists u \subseteq_p x\ (\mathrm{PsTallyPair}_{21}(u,21,p_1{}^\frown 21)\ \&\ (2^2{}_\vee u)^\frown 2^2 \sqsubseteq x)$,

(ps2) $\exists u \subseteq_p x\ \exists z \subseteq_p u\ (\mathrm{PsTallyPair}_{21}(u,p_2,z)\ \&\ (2^2{}_\vee u)^\frown 2^2 \sqsupseteq x\ \&\ (2_\vee z)^\frown 2^2 \sqsupseteq x)$,

(ps3) $\forall w_1 \sqsubseteq x\ \forall u \subseteq_p x\ \forall p \sqsubseteq p_2\ \forall y \subseteq_p u\ (u^\frown 2^2 \sqsupseteq w_1\ \&\ p \neq p_2\ \&\ p \neq \underline{e}\ \&$

$\&\ \mathrm{PsTallyPair}_{21}(u,p,y) \rightarrow \exists w_2 \sqsubseteq x\ \exists v \subseteq_p x\ (w_1 \sqsubseteq w_2\ \&\ v^\frown 2^2 \sqsupseteq w_2\ \&$

$\&\ \mathrm{PsTallyPair}_{21}(v,p^\frown 21,y^\frown 21)\ \&\ \neg v \subseteq_p w_1\ \&$

$\&\ \forall v' \subseteq_p x\ \forall p',y' \subseteq_p v'(\mathrm{PsTallyPair}_{21}(v',t',y')\ \&\ \neg v' \subseteq_p w_1\ \&\ v' \subseteq_p w_2 \rightarrow v' = v)))$

(ps4) $\forall p \sqsubseteq p_2\ \forall u,v \subseteq_p x\ \forall y,z \subseteq_p x(\mathrm{PsTallyPair}_{21}(u,p,y)\ \&\ \mathrm{PsTallyPair}_{21}(v,p,z) \rightarrow$

$\rightarrow u = v\ \&\ y = z)$,

(ps5) $\forall u \subseteq_p x\ \forall p,y \subseteq_p u\ (\mathrm{PsTallyPair}_{21}(u,p,y) \rightarrow p \neq \underline{e}\ \&\ p \sqsubseteq p_2\ \&\ u \subseteq_1 x)$,

(ps6) $\forall p \sqsubseteq p_2\ (p \neq \underline{e}\ \&\ \mathrm{PsTally}_{21}(p) \rightarrow \exists u \subseteq_p x\ \exists y \subseteq_p u\ \mathrm{PsTallyPair}_{21}(u,p,y))$,

(ps7) $\forall p,p' \sqsubseteq t_2\ \forall w,w' \sqsubseteq x\ \forall u,v,y,y' \subseteq_p x\ (\mathrm{PsTally}_{21}(p)\ \&\ \mathrm{PsTally}_{21}(p')\ \&$

$\&\ p \neq \underline{e}\ \&\ p \sqsubseteq p'\ \&\ p \neq p'\ \&\ \mathrm{PsTallyPair}_{21}(u,p,y)\ \&\ \mathrm{PsTallyPair}_{21}(v,p',y')\ \&$

$\&\ u^\frown 2^2 \sqsupseteq w\ \&\ u'^\frown 2^2 \sqsupseteq w' \rightarrow w \sqsubseteq w')$.

We can then prove, analogously to 4.3:

4.8. For any $q, p_1, p_2 \in \Sigma^*$,

$\mathfrak{S}^+ \vDash \mathrm{PsTCSeq}(\underline{q},\underline{p}_1,\underline{p}_2) \Leftrightarrow$ $p_1$ and $p_2$ are 2,1-pseudotallies and $p_2 = (21)^m$ for some $m \geq 1$ and $q = 2^2{}^\frown q_1{}^\frown 2^2{}^\frown \ldots{}^\frown 2^2{}^\frown q_m{}^\frown 2^2$ where $q_1,\ldots,q_m$ is a (unique) pseudotally concatenation sequence for $p_1,p_2$.

Let $\quad \mathrm{PsTC}(x,y,z) \equiv: \exists u\,(\mathrm{PsTCSeq}(u,x,y)\ \&\ (2^2{}_{\vee}z)^\wedge 2^2 \sqsupseteq u)$.

From the definitions, we then have that, for any m, n, k ≥ 1,

$$\mathfrak{S}^+ \vDash \mathrm{PsTC}((21)^k,(21)^m,(21)^n) \;\Leftrightarrow\; (21)^{k\frown}(21)^m = (21)^n.$$

Note that:

4.9 $\mathfrak{S}^+ \vDash \forall x,y,z\,(\neg\mathrm{PsTC}(x,y,z) \leftrightarrow \exists w\,(\mathrm{PsTC}(x,y,z)\ \&\ \neg w = z))$.

Furthermore:

4.10 (a) Let $p_1, p_2 \in \Sigma^*$ be any 2,1-pseudotallies. Then

$$\mathrm{VW}^+ \vdash \mathrm{PsTC}(\underline{p_1},\underline{p_2},\underline{p_1{}^\frown p_2}).$$

(b) $\mathrm{VW}^+ \vdash \forall z_1,z_2\,(\mathrm{PsTC}(\underline{p_1},\underline{p_2},z_1)\ \&\ \mathrm{PsTC}(\underline{p_1},\underline{p_2},z_2) \rightarrow z_1 = z_2)$.

This is proved analogously to 4.4 using 4.7 and definition of PsTCSeq.

## §5. Interpreting Binary Concatenation Directly

In [6] Murwanashyaka introduced a very weak theory of concatenation WD formulated in the language $\mathcal{L}_C = \{0,1,\circ,\preccurlyeq\}$, where 0 and 1 are individual constants, $\circ$ a binary operation symbol and $\preccurlyeq$ a binary relational symbol, interpreted in the structure $\mathfrak{D} = (B^*,0,1,{}^\frown,\preccurlyeq^{\mathfrak{D}})$, where the domain B* consists of nonempty finite binary strings of 0's and 1's (called there 'bit strings'), $^\frown$ is the binary operation of concatenation or juxtaposition of two strings, and $\preccurlyeq^{\mathfrak{D}}$ is the initial segment relation on strings. ($\mathfrak{D}$ can be thought of as a free semigroup with two generators, extended with the initial segment relation.) Each bit string $s \in B^*$ can be recursively associated with a unique canonical $\mathcal{L}_C$ term $\underline{s}$, called its 'biteral', where $\underline{0}$ is '0', $\underline{1}$ is '1', and $\underline{s^\frown 0} =: \underline{s} \circ 0$ and $\underline{s^\frown 1} =: \underline{s} \circ 1$. Then the axioms of WD are given as the instances of the schemas

(WD1) $\underline{s} \circ \underline{t} = \underline{s^\frown t}$ for any strings s, t ∈ B*,

(WD2) $\neg(\underline{s} = \underline{t})$ for any distinct strings s, t ∈ B*,

(WD3) $\forall x\, (x \preccurlyeq \underline{t} \leftrightarrow \bigvee_{r \in I[t]} x = \underline{r})$ for each string t ∈ B*,

where $I[t] = \{ s \in \Sigma^* \mid s \preccurlyeq^{\mathfrak{D}} t\}$. Murwanashyaka proved that WD is mutually interpretable with R. Since, as shown in [3], the very weak theory of dyadic trees WT described in §2 is also mutually interpretable with R, it follows from our Theorem 3.7 that the very weak theory of binary concatenation WD is interpretable in $VW^+$.

We reformulate WD as a relational theory with infinitely many individual constants $\underline{c}_0, \underline{c}_1, \underline{c}_2, \ldots$ in the vocabulary $\mathcal{L}_{WD(R)} = \{ \underline{c}_0, \underline{c}_1, \underline{c}_2, \ldots, \preccurlyeq, C\}$, where $\preccurlyeq$ is a binary relation symbol for the initial segment relation and C a ternary relation symbol expressing the graph of the binary concatenation operation $x^\frown y = z$. Here $\underline{c}_0$, $\underline{c}_1$ stand for digits 0, 1, resp., and each $\underline{c}_i$ may be identified with a distinct canonical $\mathcal{L}_C$ term $\underline{s}$ for a finite binary (in 0's and 1's) string s ∈ B*. Hence we shall simply write $\underline{s}, \underline{t}, \ldots$ in place of $\underline{c}_i, \underline{c}_j, \ldots$. Then the axioms of the relational variant $WD_R$ of WD are given by the schemas

($WD_R1$) $C(\underline{s}, \underline{t}, \underline{s^\frown t})$ for each s, t ∈ B*,

($WD_R2$) $\neg(\underline{s} = \underline{t})$ for any distinct s, t ∈ B*,

($WD_R3$) $\forall x\, (x \preccurlyeq \underline{t} \leftrightarrow \bigvee_{s \sqsubseteq t} x = \underline{s})$ for each t ∈ B*,

($WD_R4$) $\forall x,y\ \exists z\ C(x,y,z)$

($WD_R5$) $\forall x,y,u,v\ (C(x,y,u)\ \&\ C(x,y,v) \rightarrow u = v)$.

We now proceed to investigate the possibility of interpreting binary concatenation in a theory of unary concatenation in a more direct manner by setting up an interpretation of WD in $VW^+$. For this purpose we consider a homomorphism ψ of the binary concatenation structure $<\Sigma^*, \underline{1}, \underline{2}, ^\frown>$ into itself, defined as follows:

$$\psi(\underline{1}) = \underline{1}^\frown \underline{1},$$

$$\psi(\underline{2}) = \underline{2}^\frown \underline{2},$$

$$\psi(s^\frown t) = \psi(s)^\frown 21^\frown \psi(t) \qquad \text{for any } s, t \in \Sigma^*.$$

We will use ψ to set up a translation of $\mathcal{L}_C$ into $\mathcal{L}^+$ suitable for constructing the desired interpretation. Note that, for example, ψ((12)1) = (112122)2111 whereas ψ(1(21)) = 1121(222111). The pseudotallies 21 will serve as separators of translations ψ(u) of subterms u of $\mathcal{L}_C$-terms v.

As the domain of our interpretation we intend to use the strings $\{ \psi(s) \mid s \in \Sigma^*\}$ that constitute the range of the 1-1 map $\psi: \Sigma^* \rightarrow \Sigma^*$, for which we need an explicit definition. Roughly speaking, these are sequences of consecutive duplicate 11 or 22 separated by the pseudotallies 21. We call the latter 'separators'.

We say that $s \in \Sigma^*$ is a <u>duplicate code</u> iff

(c1) s is a duplicate 11 or 22, or

(c2a) s begins with a (initial) duplicate immediately followed by a separator,

(c2b) s ends with a (terminal) duplicate immediately preceded by a separator,

(c2c) every separator in s is immediately preceded by a duplicate,

(c2d) every non-initial duplicate in s is immediately preceded by a separator.

(c2e) every non-terminal duplicate in s is immediately followed by a separator.

We then have:

5.1(a) For any $s_1, s_2 \in \Sigma^*$, if $s_1$ and $s_2$ are duplicate codes, so is $s_1 \hat{} 21 \hat{} s_2$.

(b) For any $s \in \Sigma^*$, $s \neq e$, if $s \hat{} 21 \hat{} D$ is a duplicate code and D is a duplicate string 11 or 22, then s is a duplicate code.

(c ) Suppose $t \in \Sigma^*$ is a duplicate code and $s \sqsubseteq t$. If $r \sqsupseteq s$ and r is a duplicate code that is not a duplicate string, then s is a duplicate code.

<u>Proof</u>: (a) Let $s = s_1 \hat{} 21 \hat{} s_2$. Suppose $s_1$ is 11 or 22. The (c2a) holds for s immediately, and (c2b) follows from the hypothesis that $s_2$ is a duplicate code. The same applies to conditions (c2c) and (c2d) and the case where $s_2$ is a duplicate. If neither $s_1$ nor $s_2$ are duplicates, (c2a) and (c2b) are immediate from the hypothesis. For (c2c), the separators in s are those in $s_1$ and $s_2$ plus the one that separates the two in s. That (c2c) holds for s follows from the hypothesis that $s_1$ and $s_2$ are both duplicate codes. Likewise for (c2d) because every non-initial duplicate in s is either a non-initial duplicate of $s_1$ or a non-initial duplicate of $s_2$, in which case the claim follows from the induction hypothesis for $s_1$ and $s_2$, or else the initial duplicate of $s_2$, which is immediately preceded in s by the separator that separates $s_1$ and $s_2$ in s. For (c2e), any non-terminal duplicate of $s_1\hat{}21\hat{}s_2$ is either a non-terminal duplicate of $s_2$, a non-terminal duplicate of $s_1$, or else a terminal duplicate of $s_1$. That (c2e) holds for $s_1\hat{}21\hat{}s_2$follows then immediately from the hypothesis.

(b) Suppose $s = s_0\hat{}21\hat{}D$ is a duplicate code and D is a duplicate. We may assume that $s_0$ is not a duplicate for otherwise we are done. Then (c2c), (c2d) and (c2e) for $s_0$ follow immediately from the hypothesis for s. Now by (c2c) for s the separator separating $s_0$ and D in s is immediately preceded in s by a duplicate $D_1$. Because $s_0$ is not a duplicate, $D_1$ is not an initial duplicate of $s_1$, hence by (c2d) for s we have that $D_1$, the terminal duplicate of $s_0$, is immediately preceded in s, and hence also in $s_0$, by a separator. Hence (c2b) also holds for $s_0$. For (c2a) we have from the hypothesis about s that s begins with a duplicate $D_0$. Since $s_0$ contains at least one separator (namely $D_1$), the initial duplicate $D_0$ of $s_0$, being also the initial duplicate of s, is immediately followed by a separator in s, hence also in $s_0$, then (c2a) also holds for $s_0$.

(c ) Suppose $s \sqsubseteq t$ and $r \sqsupseteq s$ where r and t are duplicate codes and r is not a duplicate string. Then $s = r_1\hat{}r$ for some $r_1 \in \Sigma^*$. We may assume that $r_1 \neq e$. Now, s is not a single digit and t is not a duplicate string. Since t is a duplicate code, by (c2a) it must begin with a duplicate string dd immediately followed by a separator. Since $lh(r) \geq 6$, then $dd\hat{}21 \sqsubseteq s$, so (c2a) holds for s. On the other hand, (c2b) holds for s because by hypothesis it holds for r. That conditions (c2c) and (c2d) hold for s follows from the fact that they hold for t by hypothesis. Therefore, s must also be a duplicate code. □

5.2(a) For each $s \in \Sigma^*$, $\psi(s)$ is a duplicate code.

(b) If $t \in \Sigma^*$ is a duplicate code, then $t = \psi(s)$ for some (unique) $s \in \Sigma^*$.

Proof: (a) follows from 5.1(a) by a straightforward induction on the length of strings in $\Sigma^*$. For (b), we argue by induction on the number n of separators in t. If $n = 0$, then t is a duplicate 11 or 22, and s is the single digit 1 or 2. Suppose t has n+1 separators. Then by (c2b) we have that $t = t_0 \hat{} 21 \hat{} D$ for some duplicate D. By 5.1(a) it follows that $t_0$ is a duplicate code, and since $t_0$ has one fewer separator than t, from the induction hypothesis it follows that $t_0 = \psi(s_0)$ for some unique string $s_0 \in \Sigma^*$. But then $t = \psi(s)$ for $s = s_0 \hat{} 1$ or $s = s_0 \hat{} 2$, depending on D. □

We then have:

5.3 For any $t \in \Sigma^*$,

t is a duplicate code $\Leftrightarrow$ t is a duplicate, or $t = D_0 \hat{} 21 \hat{} D_1 \hat{} \ldots \hat{} 21 \hat{} D_m$ for some $m \geq 1$ and duplicates $D_i$, $0 \leq i \leq m$.

Proof: From left to right. we argue by induction on the number m of separators in t. If $m = 0$, then t is a duplicate, and the claim is immediate. Suppose t has k+1 separators and the claim holds for k. Then by (c2b) we have that t has a terminal duplicate D immediately preceded by a separator, that is, $t = t_0 \hat{} 21 \hat{} D$ for some $t_0 \in \Sigma^*$, $t \neq \underline{e}$. But then $t_0$ is a duplicate code by 5.1(b) with k separators, hence $t_0 = D_0 \hat{} (\frown_{1 \leq i \leq k} (21 \hat{} D_i))$. So for $D_{k+1} = D$ we have $t = D_0 \hat{} (\frown_{1 \leq i \leq k+1} (21 \hat{} D_i))$, as claimed.

Conversely, suppose $t_0 = D_0 \hat{} (\frown_{1 \leq i \leq m} (21 \hat{} D_i))$. Let $d_j = d$ if $D_j = d_j \hat{} d_j$ for $d = 1,2$. We argue by induction on $m \geq 1$. Assume $t = D_0 \hat{} 21 \hat{} D_1$. Then for $s = d_0 \hat{} d_1$ we have $\psi(s) = \psi(d_0 \hat{} d_1) = \psi(d_0) \hat{} 21 \hat{} \psi(d_1) = D_0 \hat{} 21 \hat{} D_1 = t$ and so clearly t is a duplicate code. Let $t = D_0{}^\frown(\frown_{1\leq i\leq k+1}(21^\frown D_i))$, assuming as the induction hypothesis that the claim holds for k. By 5.2(b) we have that $D_0{}^\frown(\frown_{1\leq i\leq k}(21^\frown D_i)) = \psi(s_0)$ for some $s_0 \in \Sigma^*$. But then

$$D_0{}^\frown(\frown_{1\leq i\leq k+1}(21^\frown D_i)) = D_0{}^\frown(\frown_{1\leq i\leq k}(21^\frown D_i))^\frown(21^\frown D_{k+1}) =$$

$$= \psi(s_0)^\frown 21^\frown D_{k+1} = \psi(s_0)^\frown 21^\frown \psi(d_{k+2}) = \psi(s_0{}^\frown d_{k+2}).$$

Then the claim holds for t by 5.2(a). ◻

5.4 For any $s, t \in \Sigma^*$, $s \sqsubseteq t \Leftrightarrow \psi(s) \sqsubseteq \psi(t)$.

Proof: Suppose $s \sqsubseteq t$. We argue by induction on lh(t). If $t = e$, then $s = e$, and $\psi(s) = e \sqsubseteq e = \psi(t)$. Assume, as the induction hypothesis, that the claim holds for t, and suppose $s \sqsubseteq t^\frown d$, for $d = 1,2$. Then $s \sqsubseteq t$ or $s = t^\frown d$, whence, taking into account the induction hypothesis, $\psi(s) \sqsubseteq \psi(t)$ or $\psi(s) = \psi(t^\frown d)$. But $\psi(t^\frown d) = \psi(t)^\frown 21^\frown \psi(d)$, by definition of $\psi$, so $\psi(t) \sqsubseteq \psi(t^\frown d)$, as needed. For the converse, assume that $\psi(s) \sqsubseteq \psi(t)$ and suppose, for a reductio, that $s \not\sqsubseteq t$. Then for some $s_0 \sqsubseteq s$ we have that $s_0 \sqsubseteq t$ whereas $s_0{}^\frown d_0 \sqsubseteq s$ and $s_0{}^\frown d_1 \sqsubseteq t$ where $d_0 \neq d_1$, with $d_0, d_1 \in D = \{1,2\}$. From the first part of the proof it then follows that $\psi(s_0{}^\frown d_0) \sqsubseteq \psi(s)$ and $\psi(s_0{}^\frown d_1) \sqsubseteq \psi(t)$ whereas from the hypothesis $\psi(s) \sqsubseteq \psi(t)$ we derive that $\psi(s_0{}^\frown d_0) \sqsubseteq \psi(t)$. But by definition of $\psi$, we have that

$$\psi(s_0{}^\frown d_0) = \psi(s_0)^\frown 21^\frown \psi(d_0) \text{ and } \psi(s_0{}^\frown d_1) = \psi(s_0)^\frown 21^\frown \psi(d_1).$$

This is a contradiction because $\psi(d_0) \neq \psi(d_1)$ since $d_0 \neq d_1$. Therefore, there is no such $s_0$ and we must have $s \sqsubseteq t$ if $\psi(s) \sqsubseteq \psi(t)$. ◻

For our purposes it is essential that the explicit characterization of duplicate codes in (c1)-(c2e) is expressible in $\mathcal{L}^+$ by a bounded formula. Let $D = \{ 11, 22\}$ and let $J_0(x)$ abbreviate the conjunction of (a)-(e), where

(a) $\bigvee_{r \in D} \underline{r}\hat{\ }21 \sqsubseteq x$,

(b) $\bigvee_{r \in D} 21_{\vee}\underline{r} \sqsupseteq x$,

(c) $\forall y \sqsubseteq x\, (21 \sqsupseteq y \;\rightarrow\; \bigvee_{r \in D} \underline{r}\hat{\ }21 \sqsupseteq y)$,

(d) $\bigwedge_{r \in D} \forall y \sqsubseteq x\, (y \neq \underline{e}\ \&\ y\hat{\ }\underline{r} \sqsubseteq x \;\rightarrow\; 21 \sqsupseteq y)$,

(e) $\bigwedge_{r \in D} \forall y \sqsubseteq x\, (\underline{r} \sqsupseteq y\ \&\ y \neq x \rightarrow\; y\hat{\ }21 \sqsubseteq x)$

and let $J(x) \;\equiv:\; J_0(x) \;\mathrm{v}\; x = 11 \;\mathrm{v}\; x = 22$.

5.5 For each $t \in \Sigma^*$, $\mathrm{VW}^+ \vdash \forall x\, (J(x) \rightarrow (x \sqsubseteq \underline{\psi(t)} \leftrightarrow x\hat{\ }21 \sqsubseteq \underline{\psi(t)} \;\mathrm{v}\; x = \underline{\psi(t)}))$.

<u>Proof</u>: Suppose $M \vDash J(x)$ where $x \in M$. Assume $M \vDash x \sqsubseteq \underline{\psi(t)}$. Then $M \vDash \bigvee_{s \sqsubseteq \psi(t)} x = \underline{s}$ by 1.1 and $x \in \Sigma^*$. We proceed to reason in $\mathfrak{S}^+$. Suppose that $x \neq \psi(t)$. From $\mathfrak{S}^+ \vDash J(\underline{x})$ we have that

$$\mathfrak{S}^+ \vDash \underline{x} = 11 \;\mathrm{v}\; \underline{x} = 22 \;\mathrm{v}\; J_0(\underline{x}).$$

If $x = 11$, then $\psi(t) \neq 11$, and $\psi(t) = 22$ is ruled out from hypothesis $M \vDash \underline{x} \sqsubseteq \underline{\psi(t)}$ by 1.7. But $\mathfrak{S}^+ \vDash J(\underline{\psi(t)})$ by 5.2(a), hence $\mathfrak{S}^+ \vDash J_0(\underline{\psi(t)})$. By (a) we have that $\mathfrak{S}^+ \vDash \bigwedge_{r \in D} \underline{r}\hat{\ }21 \sqsubseteq \underline{\psi(t)}$. Now, $r = 22$ is ruled out by $\mathfrak{S}^+ \vDash \underline{x} \sqsubseteq \underline{\psi(t)}$ and $x = 11$. Then we must have $r = 11$, so $x = r$ and $\mathfrak{S}^+ \vDash \underline{x}\hat{\ }21 \sqsubseteq \underline{\psi(t)}$, as needed. Analogously if $x = 22$. Finally, suppose that $\mathfrak{S}^+ \vDash J_0(\underline{x})$. Then $\mathfrak{S}^+ \vDash \bigvee_{r \in D} 21_{\vee}\underline{r} \sqsupseteq \underline{x}$ by (b), whence $\mathfrak{S}^+ \vDash \underline{r} \sqsupseteq \underline{x}$. But then $\mathfrak{S}^+ \vDash \underline{x}\hat{\ }21 \sqsubseteq \underline{\psi(t)}$ follows from hypothesis $\mathfrak{S}^+ \vDash \underline{x} \sqsubseteq \underline{\psi(t)}\ \&\ \underline{x} \neq \underline{\psi(t)}$. Thus we have established that

$$\mathfrak{S}^+ \vDash \forall x \sqsubseteq \underline{\psi(t)}\ (J(x) \to (x \sqsubseteq \underline{\psi(t)} \to\ x\hat{\ }21 \sqsubseteq \underline{\psi(t)} \vee x = \underline{\psi(t)})$$

and $M \vDash \forall x\ (J(x) \to (x \sqsubseteq \underline{\psi(t)} \to\ x\hat{\ }21 \sqsubseteq \underline{\psi(t)} \vee x = \underline{\psi(t)}))$ follows by 1.17.

Conversely, suppose $M \vDash x\hat{\ }21 \sqsubseteq \underline{\psi(t)} \vee x = \underline{\psi(t)}$ for $x \in M$.
If $M \vDash x\hat{\ }21 \sqsubseteq \underline{\psi(t)}$, then $M \vDash \bigvee_{s \sqsubseteq \psi(t)} x\hat{\ }21 = \underline{s}$ by 1.1, so $x \in \Sigma^*$ and $\mathfrak{S}^+ \vDash \underline{x}\hat{\ }21 \sqsubseteq \underline{\psi(t)}$ by 1.17. But then $\mathfrak{S}^+ \vDash \underline{x} \sqsubseteq \underline{\psi(t)}$ and $M \vDash \underline{x} \sqsubseteq \underline{\psi(t)}$ by 1.5.
Likewise $M \vDash x \sqsubseteq \underline{\psi(t)}$ if $M \vDash x = \underline{\psi(t)}$. □

5.6 $VW^+$ interprets $WD_R$.

Proof: Let the formula $J(x)$ define the domain $D(J)$ of the interpretation. For each individual constant $\underline{c}_i$ of $WD_R$, we let $[\underline{c}_i]^J = \underline{\psi(s_i)}$ where $s_i \in \Sigma^*$ is the dyadic rendering of the string in $B^*$ denoted by $\underline{c}_i$ in $WD_R$.

Let $C_0(x,y,z)$ abbreviate

$$x\hat{\ }21 \sqsubseteq z \ \&\ 21_\vee y \sqsupseteq z \ \&\ \exists v\ \exists u \sqsubseteq v\ \exists u_1,u_2 \sqsubseteq u\ (L(x,u_1)\ \&\ L(y,u_2)\ \&\ z\ |{\le}|\ v\ \&$$

$$\&\ v = (u\hat{\ }\underline{2})\hat{\ }\underline{2}\ \&\ T_2C(u_1.u_2,u)),$$

and let $C^J(x,y,z)$ abbreviate

$$\exists! w\ (C_0(x,y,w)\ \&\ J(w)\ \&\ z = w)\ \vee\ \neg\exists! w\ (C_0(x,y,w)\ \&\ J(w)\ \&\ z = \underline{1}\hat{\ }\underline{1}).$$

Finally, let $\quad x \preccurlyeq^J y\ \equiv:\ (x \sqsubseteq y\ \&\ x\hat{\ }21 \sqsubseteq y) \vee x = y.$

We proceed to show that this translation yields a formal interpretation of $WD_R$ in $VW^+$.

That $VW^+ \vdash J([\underline{c}_i]^J)$ for each $i \in N$ follows from 5.2(a) and 1.17(a). Since the map $\psi$ is 1-1 we have that $VW^+ \vdash [WD_R2]^J$. That the J-translations of ($WD_R4$) and ($WD_R5$) are derivable in $VW^+$ follows from the definition of $C^J$. For ($WD_R3$), let $t \in \Sigma^*$. Assume $M \models x \sqsubseteq \underline{\psi(t)}$ & $x\hat{}21 \sqsubseteq \underline{\psi(t)}$ where $M \models J(x)$ and $x \in M$. Then $M \models \bigvee_{s \sqsubseteq \psi(t)} x = \underline{s}$ by 1.1. Consider such an $s \in \Sigma^*$. From $M \models J(\underline{s})$ and 1.17 we have that s is a duplicate code, hence $s = \psi(r)$ for some $r \in \Sigma^*$, by 5.2(b). But then $M \models \bigvee_{\psi(r) \sqsubseteq \psi(t)} x = \underline{\psi(r)}$, that is, we have

$$M \models \forall x\, (J(x) \rightarrow (x \sqsubseteq \underline{\psi(t)}\ \&\ x\hat{}21 \sqsubseteq \underline{\psi(t)} \vee x = \underline{\psi(t)} \rightarrow \bigvee_{\psi(r) \sqsubseteq \psi(t)} x = \underline{\psi(r)}).$$

Conversely, assume $M \models x \sqsubseteq \underline{\psi(r)}$ where $x \in M$ and $\psi(r) \sqsubseteq \psi(t)$. Then $M \models \underline{\psi(r)} \sqsubseteq \underline{\psi(t)}$, by 1.5, so $M \models x \sqsubseteq \underline{\psi(t)}$. Hence $x \in \Sigma^*$. Suppose $x \neq \psi(t)$. Then, since x is a duplicate code, we have that $\mathfrak{S}^+ \models \underline{x}\hat{}21 \sqsubseteq \underline{\psi(t)}$, whence, by 1.5, we obtain $M \models x\hat{}21 \sqsubseteq \underline{\psi(t)}$. Thus we also have

$$M \models \forall x\, (J(x) \rightarrow (\bigvee_{\psi(r) \sqsubseteq \psi(t)} x = \underline{\psi(r)} \rightarrow (x \sqsubseteq \underline{\psi(t)}\ \&\ x\hat{}21 \sqsubseteq \underline{\psi(t)}) \vee x = \underline{\psi(t)}).$$

Therefore

$$VW^+ \vdash \forall x\, (J(x) \rightarrow ((x \sqsubseteq \underline{\psi(t)}\ \&\ x\hat{}21 \sqsubseteq \underline{\psi(t)}) \vee x = \underline{\psi(t)}) \leftrightarrow \bigvee_{\psi(s) \sqsubseteq \psi(t)} x = \underline{\psi(s)})$$

and from 5.5 we have $VW^+ \vdash \forall x\, (J(x) \rightarrow (x \preccurlyeq^J [\underline{t}]^J \leftrightarrow \bigvee_{\psi(s) \sqsubseteq \psi(t)} x = [\underline{s}]^J))$.
Given 5.4 we have in fact derived $VW^+ \vdash \forall x\, (J(x) \rightarrow (x \preccurlyeq^J [\underline{t}]^J \leftrightarrow \bigvee_{s \sqsubseteq t} x = [\underline{s}]^J))$, that is, $VW^+ \vdash [WD_R3]^J$, as needed.

Finally, regarding ($WD_R1$), note that

$$\mathfrak{S}^+ \models [\underline{s\hat{}t}]^J = \underline{\psi(s\hat{}t)} = \underline{\psi(s)\hat{}21\hat{}\psi(t)} = \underline{s}^J\underline{\hat{}21\hat{}t}^J.$$

We also have that $\mathfrak{S}^+ \models \underline{s^J\hat{}21} \sqsubseteq \underline{s^J\hat{}21\hat{}t^J}\ \&\ \underline{21\hat{}t^J} \sqsupseteq \underline{s^J\hat{}21\hat{}t^J}$. Now, let m, n be such that $lh(s^J) = m$ and $lh(t^J) = n$. Then $lh(s^J\hat{}21\hat{}t^J) = m+n+2$. Hence by 4.5 we also have that $\mathfrak{S}^+ \models L(\underline{s}^J, 2^m)\ \&\ L(\underline{t}^J, 2^n)$ and further that

$$\mathfrak{S}^+ \models T_2C(\underline{2}^m, \underline{2}^n, \underline{2}^{m+n})\ \ \&\ \ \underline{s^J\hat{}21\hat{}t^J}\ |{\leq}|\ \underline{2}^{m+n+2}\ \ \&\ \ \underline{2}^{m+n+2} = (\underline{2}^{m+n})\hat{}\underline{2}\hat{}\underline{2}.$$

But then from 1.17(a), 4.4 and (VW2) it follows that $VW^+ \vdash C_0(\underline{s}^J, \underline{t}^J, \underline{s^J\hat{}21\hat{}t^J})$, that is, $VW^+ \vdash C_0(\underline{s}^J, \underline{t}^J, [\underline{s\hat{}t}]^J)$.

Assume now that $M \vDash C_0(\underline{s}^J,\underline{t}^J,z)$ for $z \in M$. Then we have that

$M \vDash \underline{s}^J{}^\wedge 21 \sqsubseteq z \;\&\; 21_\vee \underline{t}^J \sqsupseteq z \;\&\; z\,|{\le}|\,2^{m+n+2}$. But then by 4.6 it follows that $M \vDash z = [\underline{s^\frown t}]^J$. Thus we derived $VW^+ \vdash \forall z\,(C_0(\underline{s}^J,\underline{t}^J,z) \to z = [\underline{s^\frown t}]^J)$. Therefore $VW^+ \vdash C^J(\underline{s}^J,\underline{t}^J,[\underline{s^\frown t}]^J)$, that is, also $VW^+ \vdash [WD_R1]^J$. □

## §6. Interpreting WD in VW*

We now turn to the theory VW* formulated in the vocabulary $\mathcal{L}^*$. For variable-free $\mathcal{L}_{W^*}$-terms t and a variable x we introduce abbreviations $x^\wedge t$ and $t_\vee x$ inductively as in §1. We write '$\underline{1}$' and '$\underline{2}$' for '$\underline{e}^\wedge\underline{1}$' and '$\underline{e}^\wedge\underline{2}$', resp.. For any $\mathcal{L}^*$ term not containing the variable x and formula A we introduce abbreviations '$\exists x \sqsubseteq t\, A$', '$\forall x \sqsubseteq t\, A$', '$\exists x \sqsupseteq t\, A$' and '$\forall x \sqsupseteq t\, A$' in the usual way. We define $\subseteq_p$ and $\subseteq_1$ in terms of $\sqsubseteq$ and $\sqsupseteq$ as in §1. Throughout this section we let M be any model of VW*. We let $D = \{1, 2\}$.

6.1(a) For each $t \in \Sigma^*$, $VW^* \vdash \forall x\,(x \sqsubseteq \underline{t} \to \bigvee_{s\sqsubseteq t} x = \underline{s})$.

(b) Let u be any variable-free $\mathcal{L}^*$ term. Then

$$VW^* \vdash u = \underline{s} \quad \text{for some } s \in \Sigma^* \text{ such that } val_{\mathfrak{S}^*}\, u = s.$$

(c) Let u, v be any variable-free $\mathcal{L}^*$ terms. If $\mathfrak{S}^* \vDash u = v$, then $VW^* \vdash u = v$.

Proof: (a) This is proved exactly as 1.1, and (b) is proved in the same way from 6.1(a) as 1.3. Then (c) follows from (b). □

6.2(a) For each $t \in \Sigma^*$ and $d \in D$ there is a $t^* \in \Sigma^*$ and $d^* \in D$ such that

$$VW^* \vdash \underline{t \frown d} = \underline{d^* \frown t^*}.$$

(b) For each $t \in \Sigma^*$ and $d \in D$ there is a $t^* \in \Sigma^*$ and $d^* \in D$ such that

$$VW^* \vdash \underline{d \frown t} = \underline{t^* \frown d^*}.$$

Proof: (a) We argue by induction on lh(t). If t is e, let $t^* = e$ and $d^* = d$. Suppose the claim holds for $s \in \Sigma^*$ and let t be $s \frown d_1$ where $d_1 \in D$. From the induction hypothesis then $(s \frown d_1) \frown d = (d_1^* \frown s^*) \frown d$ for some $d_1^* \in D$. But $(d_1^* \frown s^*) \frown d = d_1^* \frown (s^* \frown d)$, so $t \frown d = d_1^* \frown t^*$ for $t^* = s^* \frown d$. Then the claim follows by 6.1(c). For (b), we likewise argue by induction on lh(t). If t is e, let $t^* = e$ and $d^* = d$. Suppose the claim holds for $s \in \Sigma^*$ and let t be $s \frown d_1$ where $d_1 \in D$. Then $d \frown (s \frown d_1) = (d \frown s) \frown d_1 = (s^* \frown d_2) \frown d_1$ for some $s^* \in \Sigma^*$ and $d_2 \in D$. So $d \frown t = t^* \frown d^*$, letting $t^* = s^* \frown d_2$ and $d^* = d_1$. Then use 6.1(c). □

6.3 For each $t \in \Sigma^*$, $VW^* \vdash \forall x\, (x \sqsupseteq \underline{t} \rightarrow \bigvee_{s \sqsupseteq t} x = \underline{s})$.

Proof: If t is e, we use (VW*4). Suppose $M \vDash x \sqsupseteq \underline{t \frown d}$. By 6.2(a) we have that $M \vDash x \sqsupseteq \underline{d^* \frown t^*}$ for some $t^* \in \Sigma^*$ and $d^* \in D$ such that $t \frown d = d^* \frown t^*$. By (VW*2) then $M \vDash x \sqsupseteq \underline{d^*}_\vee \underline{t^*}$, whence $M \vDash x \sqsupseteq \underline{t^*}\ \text{v}\ x = \underline{d^*}_\vee \underline{t^*}$ by (VW*5). From the induction hypothesis it then follows that $M \vDash \bigvee_{s \sqsupseteq t^*} x = \underline{s}\ \text{v}\ x = \underline{d^*}_\vee \underline{t^*}$, which means that $M \vDash \bigvee_{s \sqsupseteq d^* \frown t^*} x = \underline{s}$. But then $M \vDash \bigvee_{s \sqsupseteq t \frown d} x = \underline{s}$. □

6.4 Let u, v be any variable-free $\mathcal{L}^*$ terms. (a) If $\mathfrak{S}^* \vDash u \sqsubseteq v$, then $VW^* \vdash u \sqsubseteq v$.

(b) If $\mathfrak{S}^* \vDash u \sqsupseteq v$, then $VW^* \vdash u \sqsupseteq v$.

Proof: We argue by induction on lth'v'. Assume $\mathfrak{S}^* \vDash u \sqsubseteq v$. If v is $\underline{e}$, then u must also be $\underline{e}$, and the claim holds by (VW*3). We argue as in 1.5 if v is $w \frown \underline{d}$ for $d \in D$. Suppose v is $\underline{d}_\vee w$, and assume $\mathfrak{S}^* \vDash u \sqsubseteq \underline{d}_\vee w$. Let $w^* \in \Sigma^*$ be such that $val_{\mathfrak{S}^*}\, w = w^*$. Then $\mathfrak{S}^* \vDash u \sqsubseteq \underline{d \frown w^*}$ from 6.1(c ), whence $\mathfrak{S}^* \vDash u \sqsubseteq \underline{w_1 \frown d_1}$ for some $w_1 \in \Sigma^*$, $d_1 \in D$ such that $d \frown w^* = w_1 \frown d_1$, by 6.2(b). Then $\mathfrak{S}^* \vDash u \sqsubseteq \underline{w_1}$^$\underline{d_1}$, whence $\mathfrak{S}^* \vDash u \sqsubseteq \underline{w_1}$ v $u = \underline{w_1}$^$\underline{d_1}$. If $\mathfrak{S}^* \vDash u \sqsubseteq \underline{w_1}$ then $VW^* \vdash u \sqsubseteq \underline{w_1}$ follows by the induction hypothesis since lth'$\underline{w_1}$' < lth'$\underline{d}_\vee$w'. But then $VW^* \vdash u \sqsubseteq \underline{w_1}$^$\underline{d_1}$ follows by (VW5), whence $VW^* \vdash u \sqsubseteq \underline{d}_\vee w$ by 6.1(c). If $\mathfrak{S}^* \vDash u = \underline{w_1}$^$\underline{d_1}$, then $VW^* \vdash u = \underline{w_1}$^$\underline{d_1}$ by 6.1(c), and by (VW5), we obtain that $VW^* \vdash u \sqsubseteq \underline{w_1}$^$\underline{d_1}$. But then $VW^* \vdash u \sqsubseteq \underline{d}_\vee w$ follows by 6.1(c). The proof of (b) is entirely analogous to that of (a). □

6.5 Let u be any variable-free $\mathcal{L}^*$-term, and let $u^* \in \Sigma^*$ and $u^* = val_{\mathfrak{S}^*} u$. Then:

(a) $VW^* \vdash \forall x\,(x \sqsubseteq u \leftrightarrow \bigvee_{s \sqsubseteq u^*} x = \underline{s})$,

(b) $VW^* \vdash \forall x\,(x \sqsupseteq u \leftrightarrow \bigvee_{s \sqsupseteq u^*} x = \underline{s})$.

The proof proceeds exactly as that of 1.6, appealing to 6.1, 6.3 and 6.4. □

6.6 For any variable-free $\mathcal{L}^*$-terms u, v:

(a) If $\mathfrak{S}^* \vDash \neg u \sqsubseteq v$, then $VW^* \vdash \neg u \sqsubseteq v$.

(b) If $\mathfrak{S}^* \vDash \neg u \sqsupseteq v$, then $VW^* \vdash \neg u \sqsupseteq v$.

Proof: Both proofs follows the same pattern of the proof of 1.7 making use of 6.5. □

The reformulation of 1.12 for $\mathcal{L}^*$-formulae can now be stated and proved analogously. Likewise, we can state the corresponding reformulations of 1.13, 1.14 and 1.15, resp., for $\mathcal{L}^*$-formulae and $\mathcal{L}^*$-terms with analogous proofs. After introducing the abbreviations '$\exists x \subseteq_p t\, A$' and '$\forall x \subseteq_p t\, A$' for $\mathcal{L}^*$-formulae and $\mathcal{L}^*$-terms, we define the class of bounded $\mathcal{L}^*$-formulae inductively in the same way as for $\mathcal{L}^+$ omitting the clauses relating to $|\leq|$, and likewise the corresponding class of $\Sigma$-formulae for $\mathcal{L}^*$. We then have:

6.7 Let φ be any bounded $\mathcal{L}^*$-sentence. Then: (a) if $\mathfrak{S}^* \vDash \varphi$, then $VW^* \vdash \varphi$,

(b) if $\mathfrak{S}^* \not\models \varphi$, then VW* $\vdash \neg\varphi$.

6.8 For any $\Sigma$-sentence $\varphi$, if $\mathfrak{S}^* \models \varphi$, then VW* $\vdash \varphi$.

These are proved analogously to 1.17 and 1.18.

We now turn to construction of our formal interpretation of WD in VW*. We assume that needed definitions previously given in $\mathcal{L}^+$ cast in terms of the non-logical vocabulary of $\mathcal{L}^+$ are appropriately reformulated as definitions in $\mathcal{L}^*$.

Let $z \underline{\in}_u x \equiv$: $PsTally_{21}(z)$ & $u \sqsubseteq x$ & $z \sqsupseteq u$ & $\neg 21\hat{\ }z \sqsupseteq u$ & $\neg u\hat{\ }21 \sqsubseteq x$, reading "z is a <u>locally maximal 2,1-pseudotally</u> (l.m.p.) in x determined by the initial segment u".

6.9 Let s be a variable-free $\mathcal{L}^*$-term. Then:

(a) VW* $\vdash \forall x,y \subseteq_p s\ (PsTally_{21}(x)\ \&\ PsTally_{21}(y) \rightarrow \bigwedge_{d \in D} \neg(x \sqsupseteq s\ \&\ y\hat{\ }\underline{d} \sqsupseteq s))$.

(b) VW* $\vdash \forall u \sqsubseteq s\ \forall x,y \subseteq_p u\ (x \underline{\in}_u s\ \&\ y \underline{\in}_u s \rightarrow x = y)$.

<u>Proof:</u> We reason in $\mathfrak{S}^*$. For (a), suppose for a reductio that $\mathfrak{S}^* \models x \sqsupseteq s\ \&\ y\hat{\ }\underline{d} \sqsupseteq s$, where $s^* \in \Sigma^*$ such that $s^* = val_{\mathfrak{S}^*} s \in \Sigma^*$, and $d \in D$ with x and y being 2,1-pseudotallies. Then $w_1 \hat{\ } u_1 \hat{\ } 21 = s^* = w_2 \hat{\ } y \hat{\ } d$, where $x = u_1 \hat{\ } 21$ for some (possibly empty) $u_1, w_1, w_2 \in \Sigma^*$. Then $d = 1$, and $w_1{}^\frown u_1{}^\frown 2 = w_2{}^\frown y$, contradicting the hypothesis that y is a 2,1-pseudotally. Hence $\mathfrak{S}^* \vDash \neg(x \sqsupseteq s \;\&\; y^\frown\underline{d} \sqsupseteq s)$. For (b), assume $\mathfrak{S}^* \vDash x \underline{\in}_u s \;\&\; y \underline{\in}_u s$ where $u \sqsubseteq s$. Then $x \sqsupseteq u$ and $y \sqsupseteq u$ are 2,1-pseudotallies such that $\neg 21^\frown x \sqsupseteq u$ and $\neg 21^\frown y \sqsupseteq u$. Suppose, for a reductio, that $x \sqsubseteq y$ and $x \neq y$. Then $x = (21)^j$ and $y = (21)^{j+k}$ for some $j \geq 1$, $k \geq 1$, by 4.1. But then

$$(21)^k{}^\frown x = (21)^k{}^\frown (21)^j = (21)^{j+k} = y \sqsupseteq u,$$

whence $21^\frown x \sqsupseteq u$, contradicting hypothesis. Likewise if $y \sqsupseteq x$ and $y \neq x$. Hence we must have $x = y$, as claimed. To derive (a) and (b) we now apply 6.7(a). □

Let Block(w) abbreviate the formula

$$21 \sqsubseteq w \;\&\; 21 \sqsupseteq w \;\&\; \neg 1^3 \subseteq_p w \;\&\; \neg 2^3 \subseteq_p w \;\&\; \bigvee_{d \in D} (\underline{d}^\frown\underline{d} \subseteq_1 w \;\&$$

$$\& \bigwedge_{d^* \in D} (\underline{d}^* \neq \underline{d} \rightarrow \neg\, \underline{d}^{*\frown}\underline{d}^* \subseteq_p w)),$$

reading "w is a 2,1-pseudotally block". Note that w contains a single occurrence of a duplicate string. We then have:

6.10 For any $p \in \Sigma^*$,

$\mathfrak{S}^* \vDash \mathrm{Block}(\underline{p}) \Leftrightarrow$ there are (unique) 2,1-pseudotallies $p_1, p_2$ such that

$$p = p_1{}^\frown d^\frown p_2 \text{ for some digit d.}$$

<u>Proof:</u> Assume $\mathfrak{S}^* \vDash \mathrm{Block}(\underline{p})$. Then $dd \subseteq_p p$ and $dd \sqsupseteq u$ for some $u \sqsubseteq p$. Since $21 \sqsubseteq p$, we have $u \neq dd$, so $u_1{}^\frown dd = u$ for some $u_1 \neq e$, and $u_1$ has no duplicate parts because dd is the only duplicate occurrence in p. From $u \sqsubseteq p$ we have that $u^\frown u_2 = p$ for some $u_2 \in \Sigma^*$. Since $21 \sqsupseteq p$, we have $u_2 \neq e$, and since dd is the only duplicate occurrence in p, there are no duplicate segments in $u_2$.

Case 1: $d = 1$.

We must have $2 \sqsupseteq u_1$ because $\neg 1^3 \subseteq_p p$. We claim that both $u_1\hat{}\,d$ and $u_2$ are 2,1-pseudotallies which will suffice because $(u_1\hat{}\,d)\hat{}\,d\hat{}\,u_2 = p$. If $u_1 = 2$, then $u_1\hat{}\,d = 21$, as claimed. Suppose $u_1 \neq 2$. Then from $21 \sqsubseteq p$ and $u_1 \sqsubseteq p$ we have $21 \sqsubseteq u_1$. Let k be the largest integer such that $(21)^k \sqsubseteq u_1$. We cannot have $u_1 = (21)^k$ because $u_1\hat{}\,11 = u$ and $\neg 1^3 \subseteq_p p$. Hence $(21)^k \sqsubseteq u_1$ and $u_1 \neq (21)^k$, and so $(21)^k\hat{}\,1 \sqsubseteq u_1$ or $(21)^k\hat{}\,2 \sqsubseteq u_1$. But $\neg(21)^k\hat{}\,1 \sqsubseteq u_1$ because otherwise $dd \subseteq_p u_1$, which is impossible because $u_1$ contains no duplicates. Hence we must have $(21)^k\hat{}\,2 \sqsubseteq u_1$. We claim that $(21)^k\hat{}\,2 = u_1$. Otherwise $(21)^k\hat{}\,21 \sqsubseteq u_1$, which is ruled out by the choice of k, or $(21)^k\hat{}\,22 \sqsubseteq u_1$, which is impossible because $u_1$ has no duplicates. But then $u_1\hat{}\,d = (21)^k\hat{}\,21$ is a 2,1-pseudotally, as claimed. Then $p = u\hat{}\,u_2 = (u_1\hat{}\,dd)\hat{}\,u_2$. Again, since $\neg 1^3 \subseteq_p p$, we have $2 \sqsubseteq u_2$, and since $21 \sqsupseteq p$, we cannot have $u_2 = 2$. Now, since 11 is the only duplicate in p, we must have $21 \sqsubseteq u_2$. But $u_2$ has no duplicate parts, so it must be a 2,1-pseudotally.

Case 2: $d = 2$.

Since $\neg 2^3 \subseteq_p p$, we must have $1 \sqsupseteq u_1$. We claim that both $u_1$ and $d\hat{}\,u_2$ are 2,1-pseudotallies which will suffice because $u_1\hat{}\,d\hat{}\,(d\hat{}\,u_2) = p$. We have from $21 \sqsubseteq p$ that $2 \sqsubseteq u_1$, and from $1 \sqsupseteq u_1$, that $u_1 \neq 2$. So $21 \sqsubseteq u_1$. Since $u_1$ contains no duplicates, it must be a 2,1-pseudotally. On the other hand, since $\neg 2^3 \subseteq_p p$, we also have that $\neg 2 \sqsubseteq u_2$, that is, $1 \sqsubseteq u_2$, hence $d\hat{}\,1 = 21 \sqsubseteq d\hat{}\,u_2$. We also have from $21 \sqsupseteq p$ that $1 \sqsupseteq u_2$. If $u_2 = 1$, then $d\hat{}\,u_2$ is 21 and we are done.

Otherwise, let j be the largest integer such that $(21)^j \sqsupseteq u_2$. We cannot have $(21)^j = u_2$, because $p = u_1 \hat{} d \hat{} d \hat{} u_2$ and $\neg 2^3 \subseteq_p p$. Then $(21)^j \sqsupseteq u_2$ and $(21)^j \neq u_2$, and hence $1 \hat{} (21)^j \sqsupseteq u_2$ or $2 \hat{} (21)^j \sqsupseteq u_2$. The latter is ruled out because $u_2$ contains no duplicates. So we must have $1 \hat{} (21)^j \sqsupseteq u_2$. We claim that $u_2 = 1 \hat{} (21)^j$. Otherwise $1\hat{}1\hat{}(21)^j \sqsupseteq u_2$, which is impossible because $u_2$ has no duplicates, or $2 \hat{} 1 \hat{} (21)^j \sqsupseteq u_2$, which is ruled out by the choice of j. Hence $u_2 = 1 \hat{} (21)^j$, and $d \hat{} u_2 = 21 \hat{} (21)^j$ is a 2,1-pseudotally, as claimed.

For uniqueness, let $p_1 \hat{} d \hat{} p_2 = p = p_3 \hat{} d \hat{} p_4$, and suppose, for a reductio, that $p_1 \neq p_3$, say $p_1 \sqsubseteq p_3$. Then $p_1 = (21)^j \sqsubseteq (21)^{j+k} = p_3$ where $k \geq 1$.

From $(21)^j \hat{} d \hat{} p_2 = (21)^{j+k} \hat{} d \hat{} p_4$ we then have $d \hat{} p_2 = (21)^k \hat{} d \hat{} p_4$. But this is impossible because $2 \sqsubseteq p_2$. Hence $p_1 = p_3$. But then $p_3 = p_4$ also follows.

Conversely, suppose $p = p_1 \hat{} d \hat{} p_2$ where $p_1$ and $p_2$ are 2,1-pseudotallies and d a single digit. Then $1 \sqsupseteq p_1$ and $2 \sqsubseteq p_2$, and whether d = 1 or d = 2 the string p cannot contain 3 consecutive 1's or 2's. If d = 1, then $1 \hat{} 1 \subseteq_p p$ and, since neither $p_1$ and $p_2$ contain any duplicates, that is the only duplicate occurrence in s. Likewise if d = 2. Hence $\mathfrak{S}^* \vDash \text{Block}(\underline{p})$ holds. □

Let PowerBlock(w,z) abbreviate the $\mathcal{L}^*$-formula

$$\text{Block}(w) \,\&\, \text{Pseudotally}_{21}(z) \,\&\, \bigvee_{d \in D} (z\hat{}\underline{d} \sqsubseteq w \,\&\, \underline{d}_{\vee}(z\hat{}21) \sqsupseteq w) \,\&$$

$$\&\, \neg z\hat{}21 \sqsubseteq w \,\&\, \neg 21_{\vee}(z\hat{}21) \sqsupseteq w.$$

Then we have:

6.11(a) For any $p, p_1, p_2 \in \Sigma^*$, $d \in D$,

$\mathfrak{S}^* \vDash \forall z\, (\text{PowerBlock}(\underline{p},z) \,\&\, \text{Pseudotally}_{21}(\underline{p}_1) \,\&\, \text{Pseudotally}_{21}(\underline{p}_2) \,\&$

$\& \underline{p} = (\underline{p_1}\hat{\ }\underline{d})\hat{\ }\underline{p_2} \rightarrow \underline{p_1} = z \ \& \ \underline{p_2} = z\hat{\ }21).$

(b) $\mathfrak{S}^* \vDash \forall w_1,w_2,z$ (PowerBlock($w_1$,z) & PowerBlock($w_2$,z) &

$\& \bigvee_{d \in D} (z\hat{\ }\underline{d} \sqsubseteq w_1 \ \& \ z\hat{\ }\underline{d} \sqsubseteq w_2) \rightarrow w_1 = w_2).$

Proof: Assume that $\mathfrak{S}^* \vDash$ PowerBlock($\underline{p}$,z) and $p = p_1\hat{\ }d\hat{\ }p_2$ where p, z, $p_1$, $p_2 \in \Sigma^*$, and $p_1$,$p_2$ are 2,1-pseudotallies with $d \in D$. We have that $z\hat{\ }d_0 \sqsubseteq p$ for some $d_0 \in D$. Since $p_1$ is a 2,1-pseudotally, we have from $p_1 \sqsubseteq p$ that $p_1 \sqsubseteq z$ by choice of z. Suppose for a reductio that $p_1 \neq z$. Then $p_1\hat{\ }21 \sqsubseteq z$. But then $p_1\hat{\ }21\hat{\ }z_1 = z$ for some $z_1 \in \Sigma^*$. So we have

$$(p_1\hat{\ }21\hat{\ }z_1)\hat{\ }d_0\hat{\ }w_1 = p = p_1\hat{\ }d\hat{\ }p_2$$

for some $w_1 \in \Sigma^*$, whence $21\hat{\ }z_1\hat{\ }d_0\hat{\ }w_1 = d\hat{\ }p_2$, contradicting the hypothesis that $p_2$ is a 2,1-pseudotally. Hence $p_1 = z$. On the other hand, we have that $d_1\hat{\ }(z\hat{\ }21) \sqsupseteq p$ for some $d_1 \in D$ and also $p_2 \sqsupseteq z\hat{\ }21$ by choice of $z\hat{\ }21$, since $p_2$ is a 2,1-pseudotally. Suppose, for a reductio, that $p_2 \neq z\hat{\ }21$. Then $21\hat{\ }p_2 \sqsupseteq z\hat{\ }21$. But then $z_2\hat{\ }(21\hat{\ }p_2) = p = z\hat{\ }21$ for some $z_2 \in \Sigma^*$, whence $p_1\hat{\ }d\hat{\ }p_2 = p = w_2\hat{\ }(d_1\hat{\ }z\hat{\ }21) = w_2\hat{\ }d_1\hat{\ }(z_2\hat{\ }21\hat{\ }p_2)$ for some $w_2 \in \Sigma^*$. We then obtain $p_1\hat{\ }d = w_2\hat{\ }d_1\hat{\ }z_2\hat{\ }21$, which contradicts the hypothesis that $p_1$ is a 2,1-pseudotally. Hence also $p_2 = z\hat{\ }21$. This proves (a). For (b), assume that $\mathfrak{S}^* \vDash$ PowerBlock($w_1$,z) & PowerBlock($w_2$,z) where $w_1$, $w_2$, $z \in \Sigma^*$ and $z\hat{\ }d \sqsubseteq w_1$ and $z\hat{\ }d \sqsubseteq w_2$ for some $d \in D$. Then $\mathfrak{S}^* \vDash$ Block($\underline{w_1}$) & Block($\underline{w_2}$), so by 6.10 there are unique 2,1-pseudotallies $p_{11}$, $p_{12}$, $p_{21}$, $p_{22}$ and single digits $d_1$, $d_2$ such that $w_1 = p_{11}\hat{\ }d_1\hat{\ }p_{12}$ and $w_2 = p_{21}\hat{\ }d_2\hat{\ }p_{22}$. By (a) from $\mathfrak{S}^* \vDash$ PowerBlock($\underline{w_1}$,z) & PowerBlock($\underline{w_2}$,z) it follows that $p_{11} = z = p_{21}$ and $p_{12} = z\hat{\ }21 = p_{22}$. Hence from $z\hat{\ }d \sqsubseteq w_1 = p_{11}\hat{\ }d_1\hat{\ }p_{12}$ it follows that $d = d_1$, and likewise from $z\hat{\ }d \sqsubseteq w_2 = p_{21}\hat{\ }d_2\hat{\ }p_{22}$ we have $d = d_2$. But then

$w_1 = p_{11}\hat{}\, d_1\hat{}\; p_{12} = z\hat{}\, d\hat{}\, z\hat{}\, 21 = p_{21}\hat{}\, d_2\hat{}\; p_{22} = w_2$ as claimed. □

Let BlockSeq(x) abbreviate the formula

$$21 \sqsubseteq x \;\&\; 21 \sqsupseteq x \;\&\; \neg 1^3 \subseteq_p x \;\&\; \neg 2^3 \subseteq_p x \;\&\; \bigvee_{d \in D} \underline{d}\hat{}\underline{d} \subseteq_p x \;\&$$

$$\&\; \bigwedge_{d,d^* \in D} \neg(\underline{d}\hat{}\underline{d})\hat{}(\underline{d^*}\hat{}\underline{d^*}) \subseteq_p x,$$

reading "x is a block sequence". Then we have

6.12 For any $s \in \Sigma^*$,

$\mathfrak{S}^* \vDash \mathrm{BlockSeq}(\underline{s}) \Leftrightarrow$ for some $m \geq 2$, $s = p_1\hat{}\, d_1\hat{}\, \ldots\hat{}\, d_{m-1}\hat{}\, p_m$ for some (unique) locally maximal 2,1-pseudotallies $p_1,\ldots,p_m$ and digits $d_1,\ldots,d_{m-1}$.

<u>Proof:</u> Suppose $\mathfrak{S}^* \vDash \mathrm{BlockSeq}(\underline{s})$. Let $(21)^j$ be the longest proper initial segment of s that is a 2,1-pseudotally, which must exist because s contains a duplicate part $d\hat{}\, d$. Then $(21)^j\hat{}\, d \sqsubseteq s$. Let $p_1 = (21)^j$ and $d_1 = d$. Assuming that we have found l.m.p.'s $p_1,\ldots,p_k$, and digits $d_1,\ldots,d_k$ such that $s_k = \frown_{1\leq i\leq k} (p_i\hat{}\, d_i) \sqsubseteq s$, we show how to find $p_{k+1}$. Since $21 \sqsupseteq s$ by hypothesis and $1 \sqsupseteq p_k$, we have that $s_k \neq s$.

Case 1. $d_k = 1$.

Then $11 \sqsupseteq s_k$, and we must have $s_k\hat{}\, 2 \sqsubseteq s$ since $\neg 1^3 \subseteq_p s$. Now, $s_k\hat{}\, 2 \neq s$ because $21 \sqsupseteq s$. Hence $s_k\hat{}\, 21 \sqsubseteq s$ or $s_k\hat{}\, 22 \sqsubseteq s$. But the latter is ruled out because $11 \sqsupseteq s_k$ and $s_k$ has no consecutive duplicates. Hence $s_k\hat{}\, 21 \sqsubseteq s$. Let $s_k\hat{}\, s^* = s$. Let $p_{k+1}$ be the longest initial segment $(21)^q$ of $s^*$ such that $s_k\hat{}\, (21)^q \sqsubseteq s$.

Case 2. $d_k = 2$.

Then $12 \sqsupseteq s_k$. Since $21 \sqsupseteq s$, we must have $s_k \sqsubseteq s$ and $s_k \neq s$. Hence $s_k\hat{}\,1 \sqsubseteq s$ or $s_k\hat{}\,2 \sqsubseteq s$. If $s_k\hat{}\,1 \sqsubseteq s$, then $p_k\hat{}\,d_k\hat{}\,1 = p_k\hat{}\,21$, contradicting the hypothesis that the 2,1-pseudotally $p_k$ is locally maximal in s. So we must have $s_k\hat{}\,2 \sqsubseteq s$. Since $21 \sqsupseteq s$, then $s_k\hat{}\,2 \neq s$, so $s_k\hat{}\,21 \sqsubseteq s$ or $s_k\hat{}\,22 \sqsubseteq s$. But the latter is ruled out because $d_k = 2$ and $\neg 2^3 \subseteq_p s$. Hence we must have $s_k\hat{}\,21 \sqsubseteq s$, and we may again let $p_{k+1}$ be the longest initial subsegment $(21)^q$ of s* such that $s_k\hat{}\,(21)^q \sqsubseteq s$.

If $p_{k+1} = s^*$, we are done; otherwise, let $d_{k+1} = d$ where $p_{k+1}\hat{}\,d \sqsubseteq s^*$. Then the uniqueness of $p_1,\ldots,p_m$ and digits $d_1,\ldots,d_{m-1}$ follows by construction.

Conversely, suppose $s = p_1\hat{}\,d_1\hat{}\,\ldots\hat{}\,d_{m-1}\hat{}\,p_m$ for some l.m.p.'s $p_i$ in s and digits $d_j$, $1 \leq i \leq m$, $1 \leq j < m$, for $m \geq 2$. Clearly, there cannot be any 3 consecutive 1's or 2's in s. There is at least one duplicate in s, but all the duplicates in s are of the form $d\hat{}\,d_j$ where $d \sqsupseteq p_j$ or $d_j\hat{}\,d$ where $d \sqsubseteq p_{i+1}$ depending on whether $d_j = 1$ or $d_j = 2$, respectively. Since $p_j$ and $p_{j+1}$ are 2,1-pseudotallies, there cannot be consecutive duplicates in s. Hence $\mathfrak{S}^* \vDash \text{BlockSeq}(\underline{s})$. □

Note that if $\mathfrak{S}^* \vDash \text{BlockSeq}(\underline{s})$ where $s = s^*\hat{}\,d\hat{}\,p$ for some digit d and l.m.p. p in s, then also $\mathfrak{S}^* \vDash \text{BlockSeq}(\underline{s^*})$ unless s* is a 2,1-pseudotally.

Let $s = p_1\hat{}\,d_1\hat{}\,\ldots\hat{}\,d_{m-1}\hat{}\,p_m$ be a block sequence. For any digit d, we say that <u>$s\hat{}\,d$ is a blockchain for the string $d_1\hat{}\,\ldots\hat{}\,d_{m-1}\hat{}\,d$</u> if and only if $p_1 = 21$ and for each i, $1 \leq i < m$, $p_{i+1} = p_i\hat{}\,21$.

Let BlockChain(x,y) abbreviate the conjunction of (bc1)-(bc5) where

(bc1) $21 \sqsubseteq x$ & $\exists x_1 \sqsubseteq x \bigvee_{d \in D} ((x = x_1\hat{}\underline{d}$ &

$\& ((x_1 = 21 \; \& \; y = \underline{d}) \vee (\mathrm{BlockSeq}(x_1) \; \& \; \underline{d} \sqsupseteq y \; \& \; d \neq y))),$

(bc2) $\bigvee_{d \in D} 21\hat{}\underline{d} \sqsubseteq x \; \& \; \neg 21\hat{}21 \sqsubseteq x \; \& \; \underline{e} \sqsubseteq y \; \& \; \underline{e} \neq y \; \&$

$\& \bigwedge_{d \in D} ((21\hat{}\underline{d} \sqsubseteq x \leftrightarrow \underline{d} \sqsubseteq y) \; \& \; (\underline{d} \sqsupseteq x \leftrightarrow \underline{d} \sqsupseteq y))$

(bc3) $\forall y_1 \sqsubseteq y \bigwedge_{d,d^* \in D} (\underline{d} \sqsupseteq y_1 \; \& \; y_1\hat{}\underline{d^*} \sqsubseteq y \; \& \; y_1\hat{}\underline{d^*} \neq y \rightarrow$

$\rightarrow \; \exists v \sqsubseteq x \; \exists u \sqsubseteq v \; \exists w_1, w_2 \subseteq_p x \; \exists z \sqsubseteq w_1 ( w_1 \sqsupseteq v \; \& \; \mathrm{PowerBlock}(w_1,z) \; \&$

$\& \; \mathrm{PowerBlock}(w_2,z\hat{}21) \; \& \; z \, \underline{\in}_u \, x \; \& \; z\hat{}\underline{d} \sqsubseteq w_1 \; \& \; z\hat{}21 \sqsupseteq w_1 \; \&$

$\& \; z\hat{}21 \, \underline{\in}_v \, x \; \& \; v\hat{}\underline{d^*} \sqsubseteq w_2 \; \& \; (z\hat{}21)\hat{}\underline{d^*} \sqsubseteq w_2))),$

(bc4) $\forall v \sqsubseteq x \; \forall u \sqsubseteq v \; \forall z \sqsupseteq u \; (z \, \underline{\in}_u \, x \; \& \; z\hat{}21 \, \underline{\in}_v \, x \; \rightarrow$

$\rightarrow \bigwedge_{d,d^* \in D} (u\hat{}\underline{d} \sqsubseteq x \; \& \; v\hat{}\underline{d^*} \sqsubseteq x \rightarrow \exists y_2 \sqsubseteq y \; \exists y_1 \sqsubseteq y_2 \; ( \underline{d} \sqsupseteq y_1 \; \& \; y_2 = y_1\hat{}\underline{d^*} )))$

(bc5) $\forall v \sqsubseteq x \; \forall z \sqsupseteq v \; (z \, \underline{\in}_v \, x \; \rightarrow \; \bigvee_{d \in D} v\hat{}\underline{d} \sqsubseteq x \; \& \; ((v = 21 \; \& \; z = v) \vee$

$\vee \; \exists u \sqsubseteq v \; \exists w \sqsupseteq v \; \exists z_1 \sqsubseteq w \; (\mathrm{PowerBlock}(w,z_1) \; \& \; z_1 \, \underline{\in}_u \, x \; \& \; z = z_1\hat{}21))).$

6.13 Let $p, s \in \Sigma^*$ and $d \in D$. Then for some $n \geq 1$,

$\mathfrak{S}^* \vDash \mathrm{BlockChain}(\underline{p}\hat{}\underline{d},\underline{s}) \rightarrow \forall z \subseteq_p \underline{p}\hat{}\underline{d} \; (\exists u \sqsubseteq \underline{p}\hat{}\underline{d} \;\; z \, \underline{\in}_u \, \underline{p}\hat{}\underline{d} \; \leftrightarrow \bigvee_{1 \leq j \leq n} z = (21)^j).$

<u>Proof</u>: Assume that $\mathfrak{S}^* \vDash \mathrm{BlockChain}(\underline{p}\hat{}\underline{d},s)$ for some $p, s \in \Sigma^*$ and $d \in D$.

From (bc1) we have that $\mathfrak{S}^* \vDash \underline{p} = 21 \vee \mathrm{BlockSeq}(\underline{p})$. If $p = 21$, the claim is immediate for $n = 1$. So we may assume that $\mathfrak{S}^* \vDash \mathrm{BlockSeq}(\underline{p})$. Then by 6.12 there are unique l.m.p.'s $p_1,\ldots,p_m$ and digits $d_1,\ldots,d_{m-1}$ such that $p = \frown_{1 \leq i < m} (p_i \hat{}\, d_i) \hat{}\, p_m$ for some $m \geq 2$. We argue by induction on $m\text{-}1 \geq 1$. If $m\text{-}1 = 1$, then $p = 21 \hat{}\, d_1 \hat{}\, p_2$, and $p$ is a 2,1-pseudotally block $w$ such that $21 \hat{}\, d_1 \sqsubseteq w$ and $d_1 \hat{}\, p_2 \sqsupseteq w$. Then from $\mathfrak{S}^* \vDash z \, \underline{\in}_u \, p\hat{}\underline{d}$ we have by (bc5) that $z = 21$ or $z = p_2 = 21\hat{}21$, as claimed. So we may let $n = 2$. Assume now that the claim holds for $m\text{-}1 = k$ and consider $p = \frown_{1\leq i<k+1} (p_i \hat{}\, d_i)\hat{}\, p_{k+1}$. Then

$p_k \hat{}\, d_k \hat{}\, p_{k+1}$ is a 2,1-pseudotally block w such that $p_k \hat{}\, d_k \sqsubseteq w$ and $p_{k+1} \sqsupseteq w$. But then from $\mathfrak{S}^* \vDash p_{k+1} \in_u p\text{\textasciicircum}\underline{d}$ we have by (bc5) that $p_{k+1} = p_k\text{\textasciicircum}21$, and the desired claim follows from the induction hypothesis. □

6.14 Suppose $\mathfrak{S}^* \vDash \text{BlockChain}(\underline{p},\underline{s})$, where $p, s \in \Sigma^*$ and $s = \frown_{1\leq i<n} s_i$, for $s_i \in D$, $1 \leq i \leq n$. Then $p = \frown_{1\leq i<n} ((21)^i \hat{}\, d_i)$ for some $d_i \in D$, $1 \leq i \leq n$, such that, for each k, $1 \leq k \leq n$,

(a) $\mathfrak{S}^* \vDash \bigwedge_{1 \leq i \leq k} \underline{d}_i = \underline{s}_i$,

(b) $\mathfrak{S}^* \vDash \bigwedge_{d \in D}(\exists u \sqsubseteq \underline{p}\ ((21)^k \in_u \underline{p}\ \&\ u\text{\textasciicircum}\underline{d} \sqsubseteq \underline{p}) \leftrightarrow \underline{d} = \underline{s}_k)$.

<u>Proof</u>: Assume $\mathfrak{S}^* \vDash \text{BlockChain}(\underline{p},\underline{s})$, where $p, s \in \Sigma^*$ and $s = \frown_{1\leq i<n} s_i$, for $s_i \in D$, $1 \leq i \leq n$. We argue by induction on the number m of 2,1-pseudotallies in p. By (bc1) we have that $m \geq 1$. For $m = 1$ we must have

$$\mathfrak{S}^* \vDash \exists x_1 \bigvee_{d \in D} (\underline{p} = x_1\text{\textasciicircum}\underline{d}\ \&\ x_1 = 21\ \&\ \underline{s} = \underline{d}),$$

for otherwise if $\mathfrak{S}^* \vDash \exists x_1 \bigvee_{d \in D} (\underline{p} = x_1\text{\textasciicircum}\underline{d}\ \&\ \text{BlockSeq}(x_1))$ then $m \geq 2$ by 6.12. So $p = 21\hat{}\, d$ and $s = d$ for $d \in D$, and thus $s = s_1 = d_1$ and both (a) and (b) hold. Assume now that $m = k+1$ for $k \geq 1$. Then from (bc1) in the principal hypothesis $\mathfrak{S}^* \vDash \text{BlockChain}(\underline{p},\underline{s})$ we have that

$$\mathfrak{S}^* \vDash \exists x_1 \bigvee_{d \in D} (\underline{p} = x_1\text{\textasciicircum}\underline{d}\ \&\ \text{BlockSeq}(x_1)).$$

Let $p^* \in \Sigma^*$ be the block sequence such that $p = p^* \hat{}\, d$. By 6.12 there are unique l.m.p.'s $p_1,\ldots,p_{k+1}$ in $p^*$ and digits $d_1,\ldots,d_k$ such that

$$p = p^* \hat{}\, d = (p_1 \hat{}\, d_1 \hat{}\, \ldots \hat{}\, d_k \hat{}\, p_{k+1}) \hat{}\, d.$$

From 6.13 we have that $p_j = (21)^j$ for each j, $1 \leq j \leq k+1$. Now suppose that $\mathfrak{S}^* \vDash \exists u \sqsubseteq \underline{p}\ ((21)^{k+1} \in_u \underline{p}\ \&\ u\hat{}\underline{d^*} \sqsubseteq \underline{p})$ where $d^* \in D$. By (bc5) we have that

$\mathfrak{S}^* \vDash \exists v \sqsubseteq u\ \exists w \sqsupseteq u\ \exists z_1 \sqsubseteq w$ (PowerBlock$(w,z_1)$ & $z_1 \in_v x$ &

$\&\ \bigvee_{d \in D} (z_1\hat{}\underline{d} \sqsubseteq w\ \&\ \underline{d}_\vee(z_1\hat{}21) \sqsupseteq w\ \&\ (21)^{k+1} = z_1\hat{}21)$.

Then $z_1 = (21)^k$ and $w = (21)^k \frown d \frown (21)^{k+1}$ for some $d \in D$, and by choice of $p^*$, we have that $d = d_k$. Since by the induction hypothesis (a) and (b) hold for k, we have that $d_k = s_k$ and $\mathfrak{S}^* \vDash \bigwedge_{d \in D}((21)^k \sqsupseteq v\ \&\ v\hat{}\underline{d} \sqsubseteq \underline{p} \leftrightarrow \underline{d} = \underline{s}_k)$.

Then from $\mathfrak{S}^* \vDash (21)^k \in_v \underline{p}\ \&\ (21)^{k+1} \in_u \underline{p}\ \&\ v\hat{}\underline{s}_k \sqsubseteq \underline{p}\ \&\ u\hat{}\underline{d^*} \sqsubseteq \underline{p}$ we have by (bc4) that $\mathfrak{S}^* \vDash \exists y_2 \sqsubseteq \underline{s}\ \exists y_1 \sqsubseteq y_2\ (\underline{d}_k \sqsupseteq y_1\ \&\ y_2 = y_1\hat{}\underline{d^*})$. But then $y_1 = \frown_{1 \leq i \leq k} s_i$ and $y_2 = (\frown_{1 \leq i \leq k} s_i) \frown d^* = \frown_{1 \leq i \leq k+1} s_i$ that is, $d^* = s_{k+1}$, as needed.

Conversely, suppose $d = s_{k+1}$. We need to show that

$\mathfrak{S}^* \vDash \exists u \sqsubseteq \underline{p}\ ((21)^{k+1} \in_u \underline{p}\ \&\ u\hat{}\underline{d} \sqsubseteq \underline{p})$. Suppose $k+1 < n$. Then from $s_k \sqsupseteq \frown_{1 \leq i \leq k} s_i$ and $(\frown_{1 \leq i \leq k} s_i) \frown s_{k+1} \sqsubseteq s$ and $(\frown_{1 \leq i \leq k} s_i) \frown s_{k+1} \neq s$ we have by (bc3) in the principal hypothesis that

$\mathfrak{S}^* \vDash \exists v \sqsubseteq \underline{p}\ \exists u \sqsupseteq v\ \exists w_1,w_2 \subseteq_p \underline{p}\ \exists z \sqsubseteq w_1$(PowerBlock$(w_1,z)$ &

& Block$(w_2, z\hat{}21)$ & $w_1 \sqsupseteq v\ \&\ z \in_u \underline{p}\ \&\ z\hat{}\underline{s}_k \sqsubseteq w_1\ \&\ z\hat{}21 \in_v \underline{p}\ \&\ z\hat{}21 \sqsupseteq w_1$ &

& $v\hat{}\underline{s}_{k+1} \sqsubseteq x\ \&\ (z\hat{}21)\hat{}\underline{s}_{k+1} \sqsubseteq w_2)$).

From the induction hypothesis (a) and (b) for k and the uniqueness of $p_i$'s and $d_i$'s in $p^*$ it follows that $d = s_{k+1}$ and $z = (21)^k$ and $z\hat{}21 = (21)^{k+1}$. Hence we obtain $\mathfrak{S}^* \vDash \exists v \sqsubseteq \underline{p}\ ((21)^{k+1} \in_v \underline{p}\ \&\ v\hat{}\underline{s_{k+1}} \sqsubseteq \underline{p})$, as needed.

Finally, suppose that $k+1 = n$. Then $(21)^{k+1}$ is the last 2,1-pseudotally in p and the desired claim follows from (bc2). □

6.15 Let $p, s \in \Sigma^*$. For each $n \geq 2$,

$\mathfrak{S}^* \vDash \forall u \sqsubseteq \underline{p} \bigwedge_{d \in D} (\text{BlockChain}(\underline{p},\underline{s}) \ \& \ (21)^n \in_u \underline{p} \ \& \ (21)^n \ {}^\wedge\underline{d} \sqsupseteq \underline{p} \ \rightarrow$

$$\rightarrow \bigwedge_{d^* \in D} \text{BlockChain}((\underline{p}{}^\wedge(21)^{n+1}){}^\wedge\underline{d^*}, \underline{s}{}^\wedge\underline{d^*})).$$

Proof: Assume that $\mathfrak{S}^* \vDash \text{BlockChain}(\underline{p},\underline{s}) \ \& \ (21)^n \in_u \underline{p} \ \& \ (21)^n \ {}^\wedge\underline{d} \sqsupseteq \underline{p}$ where $u \sqsubseteq p$. We verify that (bc1)-(bc5) hold for $(\underline{p}{}^\frown(21)^{n+1}){}^\frown\underline{d^*}$ and $\underline{s}{}^\wedge\underline{d^*}$ for any $d^* \in D$. This is immediate for (bc2). If $p = p^*{}^\frown d$, we have from hypothesis that $\mathfrak{S}^* \vDash \text{BlockSeq}(\underline{p^*})$, and $\mathfrak{S}^* \vDash \text{BlockSeq}((\underline{p^*}{}^\wedge\underline{d}){}^\wedge (21)^{n+1})$ follows from 6.12, whence also (bc1) holds. For (bc3), assume $y_1 \sqsubseteq s{}^\frown d^*$ where $d_0 \sqsupseteq y_1$ and $y_1{}^\frown d_1 \sqsubseteq s{}^\frown d^*$, for $d_0,d_1 \in D$. We may assume that $d_1$ is the last digit d of s, for otherwise the claim follows from the hypothesis for p and s. But then from 6.12 and 6.10 it follows that $(21)^{n-1}{}^\frown d_0{}^\frown(21)^n$ and $(21)^n{}^\frown d_1{}^\frown(21)^{n+1}$ are the desired 2,1-pseudotally blocks $w_1$, $w_2$ in $p{}^\frown(21)^{n+1}{}^\frown d^*$. For (bc5), assume $\mathfrak{S}^* \vDash z \in_u \underline{q}$ where $q = p{}^\frown(21)^{n+1}{}^\frown d^*$ and $v \sqsubseteq q$. We may assume that $\neg v \sqsubseteq p$. Otherwise, from $\mathfrak{S}^* \vDash z \in_v \underline{q}$ we have that

$$\mathfrak{S}^* \vDash \text{PsTally}_{21}(z) \ \& \ z \sqsupseteq v \ \& \ \neg 21{}^\wedge z \sqsupseteq v \ \& \ \neg v{}^\wedge 21 \sqsubseteq \underline{q}.$$

Since $p \sqsubseteq q$, we have in fact that

$$\mathfrak{S}^* \vDash \text{PsTally}_{21}(z) \ \& \ v \sqsubseteq \underline{p} \ \& \ z \sqsupseteq v \ \& \ \neg 21{}^\wedge z \sqsupseteq v \ \& \ \neg v{}^\wedge 21 \sqsubseteq \underline{p},$$

that is, $\mathfrak{S}^* \vDash z \in_v \underline{p}$. Then the desired claim follows from (bc5) in the principal hypothesis. Now, if $\mathfrak{S}^* \vDash v \sqsubseteq \underline{q} \ \& \ \neg v \sqsubseteq \underline{p} \ \& \ z \in_v \underline{q}$, then we must have $v = p{}^\frown(21)^{n+1}$ because $\mathfrak{S}^* \vDash \text{PsTally}_{21}(z) \ \& \ z \sqsupseteq v \ \& \ \neg v{}^\wedge 21 \sqsubseteq \underline{q}$. And since $\neg 21{}^\wedge z \sqsupseteq v$, it follows that $z = (21)^{n+1}$. But then, since by the principal hypothesis $\mathfrak{S}^* \vDash (21)^n \in_u \underline{p} \ \& \ (21)^n \ {}^\wedge\underline{d} \sqsupseteq \underline{p}$ for some $d \in D$, we have that $w = (21)^n{}^\frown d{}^\frown(21)^{n+1}$ is the desired 2,1-pseudotally block such that $w \sqsupseteq v$. Hence (bc5) also holds for q and $s{}^\frown d^*$. □

6.16 For any $p, s \in \Sigma^*$,

$\mathfrak{S}^* \vDash \text{BlockChain}(\underline{p},\underline{s}) \iff$ p is a blockchain for s.

<u>Proof</u>: The 'if' part follows from 6.14. Conversely, assume that $p = p_1 \hat{}\, d_1 \hat{}\, \ldots \hat{}\, p_m \hat{}\, d_m$ and p is a blockchain for $s = d_1 \hat{}\, \ldots \hat{}\, d_m$. Then a straightforward induction on $m \geq 1$ using 6.15 shows that $\mathfrak{S}^* \vDash \text{BlockChain}(\underline{p},\underline{s})$. □

Let $\text{Tlmp}(z,x) \equiv: \exists y \sqsubseteq x\, (z \in_y x \,\&\, \bigvee_{d \in D} y\hat{}\underline{d} \sqsupseteq x)$,

reading "z is (the) terminal l.m.p. in x". Let Conc(x,y,z) abbreviate the formula

$(x = \underline{e} \,\&\, z = y)\ \text{v}\ (y = \underline{e} \,\&\, z = x)\ \text{v}\ (x \neq \underline{e} \,\&\, y \neq \underline{e} \,\&\, z \neq \underline{e} \,\&$

$\&\, \exists u,v,w\, (\text{BlockChain}(u,x) \,\&\, \text{BlockChain}(v,y) \,\&\, \text{BlockChain}(w,z) \,\&$

$\&\, \exists z_1 \subseteq_p u\, \exists z_2 \subseteq_p v\, \exists z_3 \subseteq_p w\, (\text{Tlmp}(z_1,u) \,\&\, \text{Tlmp}(z_2,v) \,\&\, \text{Tlmp}(z_3,w) \,\&$

$\&\, \text{PsTC}(z_1,z_2,z_3) \,\&\, \forall z_0 \sqsubseteq z_3\, \forall w_1 \sqsubseteq w \bigwedge_{d \in D} (z_0 \in_{w(1)} w \,\&\, w_1\hat{}\underline{d} \sqsubseteq w \rightarrow$

$\rightarrow (z_0 \sqsubseteq z_1 \,\&\, z_0 \in_{w(1)} u \,\&\, w_1\hat{}\underline{d} \sqsubseteq u)\ \text{v}$

$\text{v}\ (z_1\hat{}21 \sqsubseteq z_0 \,\&\, \forall z_4 \sqsubseteq z_2\, \forall v_1 \sqsubseteq v\, (\text{PsTC}(z_1,z_4,z_0) \,\&\, z_4 \in_{v(1)} v \rightarrow z_4\hat{}\underline{d} \sqsubseteq v)))) \,\&$

$\&\, \forall u_0 \sqsubseteq z \bigwedge_{d \in D} (u_0\hat{}\underline{d} \sqsubseteq z \leftrightarrow (u_0 \sqsubset x \,\&\, \exists u_1 \sqsubset x\, \exists v_1 \sqsubseteq z_1 (v_1 \in_{u(1)} u \,\&\, v_1\hat{}\underline{d} \sqsubseteq u))\ \text{v}$

$\text{v}\ (\neg u_0 \sqsubset x \,\&\, x \sqsubseteq u_0 \,\&\, \exists u_0 \sqsubset v\, \exists v_2 \sqsubseteq z_2 (v_2 \in_{u(2)} v \,\&\, v_2\hat{}\underline{d} \sqsubseteq v)))\ \&\ x \neq z)$

We then have:

6.17 For each $r, s, t \in \Sigma^*$, $\mathfrak{S}^* \vDash \text{Conc}(\underline{r},\underline{s},\underline{t}) \iff r \hat{}\, s = t$.

<u>Proof</u>: We may assume that r, s, t, are all distinct from e. Let u, v, w be the (unique) blockchains for r, s, t, resp., with terminal l.m.p.'s $u_1 = (21)^m$, $v_1 = (21)^k$, $w_1 = (21)^n$, resp., for m, k, n ≥ 1, such that $\mathfrak{S}^* \vDash PsTC((21)^m,(21)^k,(21)^n)$. So $(21)^m \hat{} (21)^k = (21)^n$. But then $n = m+k$, and for each j, $1 \leq j \leq k$, $(21)^m \hat{} (21)^j = (21)^{m+j}$. Now, we have that

(i) $u = \frown_{1\leq j\leq m} ((21)^j \hat{} d_{1(i)})$, $v = \frown_{1\leq j\leq k} ((21)^j \hat{} d_{2(i)})$ and $w = \frown_{1\leq j\leq n} ((21)^j \hat{} d_{3(i)})$,

where by 6.14 $d_{i(j)}$'s are single digits such that $r = \frown_{1\leq j\leq m} d_{1(j)}$, $s = \frown_{1\leq j\leq k} d_{2(j)}$ and $t = \frown_{1\leq j\leq n} d_{3(j)}$.

Assume that $\mathfrak{S}^* \vDash Conc(\underline{r},\underline{s},\underline{t})$. From the hypothesis we also have that $\bigwedge_{1\leq j\leq m} d_{1(j)} = d_{3(j)}$ and $\bigwedge_{1\leq j\leq k} d_{2(j)} = d_{3(m+j)}$. But then

$$\text{(ii)}\quad r \hat{} s = (\frown_{1\leq j\leq m} d_{1(j)}) \hat{} (\frown_{1\leq j\leq k} d_{2(j)}) = (\frown_{1\leq j\leq m} d_{1(j)}) \hat{} (\frown_{1\leq j\leq k} d_{3(m+j)}) =$$

$$= \frown_{1\leq j\leq n} d_{3(j)} = t,$$

as needed. Conversely, suppose that $r \hat{} s = t$ holds, with blockchains u, v and w satisfying (ii) and with terminal l.m.p.'s $(21)^m$, $(21)^k$ and $(21)^n$, where $m+k = n$. Then we have that $\mathfrak{S}^* \vDash PsTC((21)^m,(21)^k,(21)^n)$ as well as

$$\mathfrak{S}^* \vDash \forall z_0 \sqsubseteq (21)^n\ (PsTally_{21}(z_0) \leftrightarrow (PsTally_{21}(z_0)\ \&\ z_0 \sqsubseteq (21)^m)\ v$$

$$v\ \exists z \sqsubseteq (21)^k\ PsTC((21)^m,z,z_0)).$$

Finally, assume that $\mathfrak{S}^* \vDash \underline{t}_0 \sqsubset \underline{t}$ where $r \hat{} s = t$, and let $d \in D$. Then

$$\mathfrak{S}^* \vDash \underline{t}_0\hat{}\underline{d} \sqsubseteq \underline{t} \leftrightarrow (\underline{t}_0 \sqsubset \underline{r}\ \&\ \underline{t}_0\hat{}\underline{d} \sqsubseteq \underline{r})\ v\ (\underline{r} \sqsubseteq \underline{t}_0\ \&\ \underline{t}_0\hat{}\underline{d} \sqsubseteq \underline{r\hat{}s}).$$

But $\mathfrak{S}^* \vDash \underline{t}_0\hat{}\underline{d} \sqsubseteq \underline{r\hat{}s} \leftrightarrow \exists y_1 \sqsubset \underline{s}\ y_1\hat{}\underline{d} \sqsubseteq \underline{s}$ if $r \sqsubseteq t_0$, whereas

$\mathfrak{S}^* \vDash \underline{t}_0 \sqsubset \underline{r}\ \&\ \underline{t}_0\hat{}\underline{d} \sqsubseteq \underline{r} \leftrightarrow \underline{d} = \underline{d}_j$ for some j, $1 \leq j \leq m$ where $r = \frown_{1\leq i\leq m} d_i$, $d_i \in D$, and $\mathfrak{S}^* \vDash \underline{r} \sqsubseteq \underline{t}_0\ \&\ \underline{t}_0\hat{}\underline{d} \sqsubseteq \underline{r\hat{}s} \leftrightarrow \underline{d} = \underline{d}_j$ for some j, $1 \leq j \leq k$ where $s = \frown_{1\leq i\leq k} d_i$, $d_i \in D$. From 6.14 it follows that

$$\mathfrak{S}^* \vDash \underline{t}_0 \sqsubset \underline{r}\ \&\ \underline{t}_0\hat{}\underline{d} \sqsubseteq \underline{r\hat{}s} \leftrightarrow \exists u_1 \sqsubseteq \underline{u}\ ((21)^j \in_{u(1)} \underline{u}\ \&\ u_1\hat{}\underline{d} \sqsubseteq \underline{u})$$

and $\mathfrak{S}^* \vDash \underline{r} \sqsubseteq \underline{t}_0 \;\&\; \underline{t}_0\hat{}\,\underline{d} \sqsubseteq \underline{r\hat{}s} \;\leftrightarrow\; \exists v_1 \sqsubseteq \underline{v}\,((21)^j \in_{v(1)} \underline{v} \;\&\; v_1\hat{}\,\underline{d} \sqsubseteq \underline{v})$.

But then $\mathfrak{S}^* \vDash Conc(\underline{r},\underline{s},\underline{t})$ follows straightforwardly from 6.13, 6.14 and the definitions. □

6.18 For each $r, s \in \Sigma^*$, $VW^* \vdash Conc(\underline{r},\underline{s},\underline{r\hat{}s}) \;\&\; \forall z\,(Conc(\underline{r},\underline{s},z) \rightarrow \underline{r\hat{}s} \sqsubseteq z)$.

Proof: That $VW^* \vdash Conc(\underline{r},\underline{s},\underline{r\hat{}s})$ follows from 6.17 and 6.8. Assume now that

$M \vDash Conc(\underline{r},\underline{s},z)$ for $z \in M$. Suppose that $r = e$ or $s = e$. Then from $M \vDash Conc(\underline{r},\underline{s},z)$ we have that $M \vDash z = \underline{s}$ or $M \vDash z = \underline{r}$, with $M \vDash \underline{r\hat{}s} \sqsubseteq \underline{s}$ and $M \vDash \underline{r\hat{}s} \sqsubseteq \underline{r}$, resp., by 6.4(a). Either way, $M \vDash \underline{r\hat{}s} \sqsubseteq z$ follows immediately. So we may assume that $r \neq e$ and $s \neq e$. From the hypothesis we have that

$$M \vDash \forall u_0 \sqsubset z \bigwedge_{d\in D} \underline{u}_0\hat{}\,\underline{d} \sqsubseteq z \;\leftrightarrow\; (u_0 \sqsubset \underline{r} \;\&\; \exists u_1 \sqsubseteq \underline{r}\, \exists v_1 \sqsubseteq \underline{z}_1(v_1 \in_{u(1)} \underline{v} \;\&$$

$$\& \; v_1\hat{}\,\underline{d} \sqsubseteq \underline{u})) \vee (\neg u_0 \sqsubset \underline{r} \;\&\; \underline{r} \sqsubseteq u_0 \;\&\; \exists u_2 \sqsubset \underline{v}\, \exists v_2 \sqsubseteq \underline{z}_2(v_2 \in_{u(2)} \underline{v} \;\&\; v_2\hat{}\,\underline{d} \sqsubseteq \underline{v}))))$$

where $u,v,z_1,z_2 \in \Sigma^*$ such that

$$M \vDash BlockChain(\underline{u},\underline{r}) \;\&\; BlockChain(\underline{v},\underline{s}) \;\&\; BlockChain(w,z) \;\&\; Tlmp(\underline{z}_1,\underline{u}) \;\&$$

$$\& \; Tlmp(\underline{z}_2,\underline{v}) \;\&\; \underline{r} \neq z)$$

for some $w \in M$.

We now argue by induction on $lh(t)$ that $M \vDash \underline{t} \sqsubseteq z$ for each $t \sqsubseteq r\hat{}s$. Suppose $t = e$. Then $M \vDash e \sqsubseteq z$ and also $M \vDash e \sqsubset z$ follows from $M \vDash BlockChain(w,z)$ by (bc2).

Consider now $t\hat{}\,d_k \sqsubseteq r\hat{}s$ where $M \vDash \underline{t} \sqsubset z$ and $d_k \in D$. Suppose (1) that $M \vDash \underline{t} \sqsubset \underline{r}$, that is, $r = \frown_{1\leq i\leq m} d_i$, $d_i \in D$, and $1\leq k\leq m$. Then $z_1 = (21)^m$, and M ⊨ $\underline{t} \sqsubseteq \underline{r}$ by 6.7. Also by 6.7 we have from M ⊨ BlockChain($\underline{u}$,$\underline{r}$) that 𝔖* ⊨ BlockChain($\underline{u}$,$\underline{r}$), whence by 6.14 it follows that

$$\mathfrak{S}^* \vDash \exists u_1 \sqsubset \underline{u}\ ((21)^k \in_{u(1)} \underline{u}\ \&\ u_1{}^\wedge \underline{d_k} \sqsubseteq \underline{r})$$

where k≤m. By 6.7 and the hypothesis it then follows that M ⊨ $\underline{t}{}^\frown \underline{d_k} \sqsubseteq z$, as needed.

Suppose (2) that M ⊨ ¬$\underline{t} \sqsubset \underline{r}$, hence 𝔖* ⊨ ¬$\underline{t} \sqsubset \underline{r}$ by 6.7. Then 𝔖* ⊨ $\underline{r} \sqsubseteq \underline{t}$ follows from hypothesis $t \sqsubseteq r{}^\frown s$. We have from (1) that M ⊨ $\underline{t} \sqsubseteq z$ if t = r, and M ⊨ $\underline{r} \neq z$ by hypothesis. Hence M ⊨ $\underline{t} \sqsubset z$. We consider $t{}^\frown d_{k+1} \sqsubseteq r{}^\frown s$ where

$s = {}^\frown{}_{1\leq j\leq n}\, d_j$, $d_j \in D$, and 1≤k≤n. Then $z_2 = (21)^n$. Again by 6.7 from M ⊨ BlockChain($\underline{v}$,$\underline{s}$) we have that 𝔖* ⊨ BlockChain($\underline{v}$,$\underline{s}$) whence

$$\mathfrak{S}^* \vDash \exists u_2 \sqsubset \underline{v}\ ((21)^k \in_{u(2)} \underline{v}\ \&\ u_2{}^\wedge \underline{d_k} \sqsubseteq \underline{s})$$

where k≤n. By 6.7 and the hypothesis it then follows that M ⊨ $\underline{t}{}^\frown \underline{d_k} \sqsubseteq z$, as needed.

From (1) and (2) we have established that M ⊨ $\bigwedge_{t \sqsubseteq r{}^\frown s} \underline{t} \sqsubseteq z$. But then M ⊨ $\underline{r{}^\frown s} \sqsubseteq z$, as claimed. □

Let Conc*(x,y,z) abbreviate

$$\exists w\ (\mathrm{Conc}(x,y,w)\ \&\ \forall v\ (\mathrm{Conc}(x,y,v) \rightarrow w \sqsubseteq v)\ \&\ z = w)\ \ \mathrm{v}$$

$$\mathrm{v}\ \ (\neg\exists w(C(x,y,w)\ \forall v\ (\mathrm{Conc}(x,y,v) \rightarrow w \sqsubseteq v))\ \&\ z = \underline{e}).$$

We now proceed to define a direct interpretation $K$ of $WD_R$ in VW*.

For each constant $\underline{s}$, for s ∈ Σ*, of $WD_R$, let $\underline{s}^K$ be the $\mathcal{L}$*-canonical term $\underline{s}$ for s.

Let $x \preccurlyeq^K y$ ≡: $x \sqsubseteq y$, and let $C^K(x,y,z)$ ≡: Conc*(x,y,z).

We then have:

Theorem 6.19 VW* interprets $WD_R$.

Proof: We first show that VW* ⊢ $C^K(\underline{s}^K,\underline{t}^K,\underline{(s^\frown t)}^K)$ for each s, t ∈ Σ*. By 6.17 we have that $\mathfrak{S}^* \models Conc(\underline{s},\underline{t},\underline{(s^\frown t)})$. Taking into account 4.9, we have that $Conc(\underline{s},\underline{t},\underline{(s^\frown t)})$ is equivalent over $\mathfrak{S}^*$ to a $\Sigma$-sentence of $\mathcal{L}^*$. Then VW* ⊢ $Conc(\underline{s},\underline{t},\underline{(s^\frown t)})$ follows by 6.8. From 6.18 we then have that

VW* ⊢ $C^K(\underline{s}^K,\underline{t}^K,\underline{(s^\frown t)}^K)$. That VW* ⊢ $(WD_R2)^K$ is immediate by (VW1), and VW* ⊢ $(WD_R3)^K$ holds by 6.5. That VW* ⊢ $(WD_R4)^K$ and VW* ⊢ $(WD_R5)^K$ follows from the definition of Conc(x,y,z) by (VW10). □

## §7. Unary concatenation interpreted away

We reformulate VW* as a relational theory $VW^*_{Rel}$ expressed in the vocabulary $\mathcal{L}^*_{Rel} = \{\underline{c}_0, \underline{c}_1, \underline{c}_2,\ldots,R_1,R_2,R_3,R_4,\sqsubseteq,\sqsupseteq\}$ interpreted in the relational structure $<\Sigma^*,s_0,s_1,s_2,\ldots,\mathbf{R}_1,\mathbf{R}_2,\mathbf{R}_3,\mathbf{R}_4,\sqsubseteq,\sqsupseteq>$, with constants $\underline{c}_i$ standing for the distinct strings $s_i$ in Σ*, the relational symbols $R_1,R_2,R_3,R_4$ expressing the graphs $\mathbf{R}_1,\mathbf{R}_2,\mathbf{R}_3,\mathbf{R}_4$ of the unary concatenation operations $x^\frown\underline{1}$, $x^\frown\underline{2}$, $\underline{1}^\frown x$ and $\underline{2}^\frown x$, resp., and ⊑ and ⊒ interpreted as the initial and end segment relations on Σ*, resp.. We write $\underline{s}_i$ for the constants $\underline{c}_i$, omitting the indices for convenience, and $A^{\wedge}_1$ , $A^{\wedge}_2$ , $A^{\vee}_1$ , $A^{\vee}_2$ for $R_1,R_2,R_3,R_4$, resp.. Let D = {1, 2}. The axioms of $VW^*_{Rel}$ are all instances of the schemas

| | | |
|---|---|---|
| ($VW^*_{Rel}1$) | $\neg(\underline{s} = \underline{t})$ | for any distinct s, t ∈ Σ*, |
| ($VW^*_{Rel}2$) | $A^{\wedge}_d(\underline{s},\underline{s^\frown d})$ | for each s ∈ Σ*, d ∈ D |

(VW*$_{Rel}$3) $A^{\vee}_{d}(\underline{s},\underline{d\hat{\ }s})$ for each $s \in \Sigma^*$, $d \in D$

(VW*$_{Rel}$4) $\forall y\,(A^{\wedge}_{d}(\underline{t},y) \rightarrow \forall z\,(z \sqsubseteq y \leftrightarrow z \sqsubseteq \underline{t} \vee z = y)$ for each $t \in \Sigma^*$, $d \in D$

(VW*$_{Rel}$5) $\forall y\,(A^{\vee}_{d}(\underline{t},y) \rightarrow \forall z\,(z \sqsupseteq y \leftrightarrow z \sqsupseteq \underline{t} \vee z = y)$ for each $t \in \Sigma^*$, $d \in D$

(VW*$_{Rel}$6) $\forall x\, \exists y\, A^{\wedge}_{d}(x,y)$

(VW*$_{Rel}$7) $\forall x\, \exists y\, A^{\vee}_{d}(x,y)$

(VW*$_{Rel}$8) $\forall x,y,z\,(A^{\wedge}_{d}(x,y)\ \&\ A^{\wedge}_{d}(x,z) \rightarrow y = z)$

(VW*$_{Rel}$9) $\forall x,y,z\,(A^{\vee}_{d}(x,y)\ \&\ A^{\vee}_{d}(x,z) \rightarrow y = z)$

plus (VW3), (VW4), (VW3*), (VW4*) and (VW10).

We now proceed to formally interpret VW*$_{Rel}$ in the theory VW$_0$ formulated in the reduced vocabulary $\mathcal{L}_0 = \{\underline{c}_0, \underline{c}_1, \underline{c}_2, \ldots, \sqsubseteq, \sqsupseteq\}$ interpreted in the structure $\mathfrak{S}_0 = <\Sigma^*, s_0, s_1, s_2, \ldots, \sqsubseteq, \sqsupseteq>$.

For $d \in D$, let $B_d(x,y)$ abbreviate the $\mathcal{L}_{\sqsubseteq,\sqsupseteq,c}$ -formula

$$x \sqsubseteq y\ \&\ x \neq y\ \&\ \underline{d} \sqsupseteq y\ \&\ \forall u \sqsubseteq y\,(u \sqsubseteq x \vee u = y)\ \&\ \forall u \sqsubseteq x\,(u \sqsubseteq y\ \&\ u \neq y)\ \&$$
$$\&\ y \sqsubseteq y\ \&\ \textstyle\bigwedge_{d^* \in D,\, d^* \neq d} \neg \underline{d^*} \sqsupseteq y,$$

and let $B_{-d}(x,y)$ abbreviate

$$x \sqsupseteq y\ \&\ x \neq y\ \&\ \underline{d} \sqsubseteq y\ \&\ \forall u \sqsupseteq y\,(u \sqsupseteq x \vee u = y)\ \&\ \forall u \sqsupseteq x\,(u \sqsupseteq y\ \&\ u \neq y)\ \&$$
$$\&\ y \sqsubseteq y\ \&\ \textstyle\bigwedge_{d^* \in D,\, d^* \neq d} \neg \underline{d^*} \sqsubseteq y.$$

We have that:

7.1(a) For any $s, t \in \Sigma^*$ and $d \in D$, $\mathfrak{S}_0 \vDash B_d(\underline{s},\underline{t}) \Leftrightarrow t = s\hat{\ }d$.

(b) For any $s, t \in \Sigma^*$ and $d = -1,-2$, $\mathfrak{S}_0 \vDash B_d(\underline{s},\underline{t}) \Leftrightarrow t = d\hat{\ }s$.

7.2 Let $s, t \in \Sigma^*$. Then: (a) if $\mathfrak{S}_0 \vDash \underline{s} \sqsubseteq \underline{t}$, then $VW_0 \vdash \underline{s} \sqsubseteq \underline{t}$,

(b) if $\mathfrak{S}_0 \vDash \underline{s} \sqsupseteq \underline{t}$, then $VW_0 \vdash \underline{s} \sqsupseteq \underline{t}$,

(c) if $\mathfrak{S}_0 \nvDash \underline{s} \sqsubseteq \underline{t}$, then $VW_0 \vdash \neg\underline{s} \sqsubseteq \underline{t}$,

(d) if $\mathfrak{S}_0 \nvDash \underline{s} \sqsupseteq \underline{t}$, then $VW_0 \vdash \neg\underline{s} \sqsupseteq \underline{t}$.

Proof: (a) and (b) follows immediately from ($VW_0$2) and ($VW_0$3). For (c), assume $\mathfrak{S}_0 \nvDash \underline{s} \sqsubseteq \underline{t}$. Then $\mathfrak{S}_0 \vDash \neg(\underline{s} = \underline{u})$ for any $u \in \Sigma^*$ such that $u \sqsubseteq t$, so $\mathfrak{S}_0 \vDash \bigwedge_{u \sqsubseteq t} \neg(\underline{s} = \underline{u})$. But then $VW_0 \vdash \neg\underline{s} \sqsubseteq \underline{t}$ follows by ($VW_0$2). Analogously for (d). □

We can now reformulate 1.12 for $\mathcal{L}_0$ and $VW_0$, with analogous proof.

We define the class of bounded $\mathcal{L}_0$-formulae inductively as follows: (i) any atomic formula of $\mathcal{L}_0$ is a bounded $\mathcal{L}_0$-formula; (ii) if $\varphi_1$, $\varphi_2$ are bounded $\mathcal{L}_0$-formulae, so are $\neg\varphi$, $\varphi_1 \,\&\, \varphi_2$, $\varphi_1 \vee \varphi_2$; (iii) if $\varphi$ is a bounded $\mathcal{L}_0$-formula and t is a constant $\underline{c}_i$ or a variable $x_i$ distinct from x, then $\forall x \sqsubseteq t\, \varphi$, $\exists x \sqsubseteq t\, \varphi$, $\forall x \sqsupseteq t\, \varphi$ and $\exists x \sqsupseteq t\, \varphi$ are bounded $\mathcal{L}_0$-formulae.

7.3 Let $\varphi$ be a bounded $\mathcal{L}_0$-sentence. Then: (a) if $\mathfrak{S}_0 \vDash \varphi$, then $VW_0 \vdash \varphi$,

(b) if $\mathfrak{S}_0 \nvDash \varphi$, then $VW_0 \vdash \neg\varphi$.

This is proved analogously to 1.17 making use of 7.2.

Let M be any model of $VW_0$.

7.4 Let $s \in \Sigma^*$. Then:

(a) For d = 1,2, $VW_0 \vdash B_d(\underline{s},\underline{s\hat{}d})$ & $\forall z\,(B_d(\underline{s},z) \rightarrow z = \underline{s\hat{}d})$.

(b) For d = -1,-2, $VW_0 \vdash B_d(\underline{s},\underline{d\hat{}s})$ & $\forall z\,(B_d(\underline{s},z) \rightarrow z = \underline{d\hat{}s})$.

Proof: For (a), note that $B_d(\underline{s},\underline{s\hat{}d})$ is a bounded $\mathcal{L}_0$-sentence such that $\mathfrak{S}_0 \vDash B_d(\underline{s},\underline{s\hat{}d})$ by 7.1(a). Hence $VW_0 \vdash B_d(\underline{s},\underline{s\hat{}d})$ follows by 7.3(a). Assume now $M \vDash B_d(\underline{s},z)$ where $z \in M$. By 7.2(b), $M \vDash \underline{d} \sqsupseteq \underline{s\hat{}d}$, and $M \vDash \underline{d} \sqsupseteq z$ from hypothesis $M \vDash B_d(\underline{s},z)$. Assume $M \vDash y \sqsubseteq z \;\&\; y \neq z$. Then $M \vDash y \sqsubseteq \underline{s}$. Hence from $M \vDash B_d(\underline{s},\underline{s\hat{}d})$ it follows that $M \vDash y \sqsubseteq \underline{s\hat{}d} \;\&\; y \neq \underline{s\hat{}d}$. Hence $M \vDash \forall y\,(y \sqsubset z \rightarrow y \sqsubset \underline{s\hat{}d})$. Conversely, assume $M \vDash y \sqsubseteq \underline{s\hat{}d} \;\&\; y \neq \underline{s\hat{}d}$. Then $M \vDash y = \underline{r}$ for some $r \sqsubseteq s$ by ($VW_0$2) Then $M \vDash \underline{r} \sqsubseteq \underline{s}$ by 7.2(a). Hence from the hypothesis $M \vDash B_d(\underline{s},z)$ we have that $M \vDash \underline{r} \sqsubseteq z \;\&\; \underline{r} \neq z$. Therefore also $M \vDash \forall y\,(y \sqsubset \underline{s\hat{}d} \rightarrow y \sqsubset z)$. But then $M \vDash z = \underline{s\hat{}d}$ follows by ($VW_0$4). The proof of (b) is analogous to (a), making use of 7.1(b), 7.2(a), ($VW_0$3), 7.2(b) and ($VW_0$5). □

For d = 1,2,-1,-2, let

$$C_d(x,y) \;\equiv:\; \exists! z\,(B_d(x,z) \;\&\; z = y) \;\;v\;\; (\neg\exists! z\, B_d(x,z) \;\&\; y = x).$$

We can then prove:

Theorem 7.5 $VW^*_{Rel} \leq_I VW_0$.

Proof: We define a direct interpretation $L$ of $VW^*_{Rel}$ in $VW_0$ as follows. The constants $\underline{c}_0, \underline{c}_1, \underline{c}_2,\ldots$ of $\mathcal{L}^*_{Rel}$ are mapped to the constants $\underline{c}_0, \underline{c}_1, \underline{c}_2,\ldots$ of $\mathcal{L}_0$, the relational symbols $A^{\wedge}_{d}$, $d \in D$, are interpreted by the formulae $C_d(x,y)$, d = 1,2, whereas $A^{\vee}_{d}$, $d \in D$, are interpreted by $C_d(x,y)$, d = -1,-2. Finally, we let $x \sqsubseteq^{L} y \;\equiv:\; x \sqsubseteq y$ and $x \sqsupseteq^{L} y \;\equiv:\; x \sqsupseteq y$. We now proceed to verify in $VW_0$ the $\mathcal{L}$-translations of the axioms of $VW^*_{Rel}$. This is immediate for ($VW^*_{Rel}1$). For ($VW^*_{Rel}2$) and ($VW^*_{Rel}3$) this follows from 7.4(a) and 7.4(b). For ($VW^*_{Rel}4$), let $t \in \Sigma^*$ and $d \in D$. Suppose $M \vDash C_d(\underline{t},x)$ for $x \in M$.

Case 1. $M \vDash \exists! y\, B_d(\underline{t},y)$.

Then $M \vDash x = a$ for some unique $a \in M$ such that $M \vDash B_d(\underline{t},a)$, that is,

$$M \vDash \underline{t} \sqsubseteq a \;\&\; \underline{t} \neq a \;\&\; \underline{d} \sqsupseteq a \;\&\; \forall u \sqsubseteq a\, (u \sqsubseteq \underline{t} \vee u = a) \;\&\; \forall u \sqsubseteq \underline{t}\, (u \sqsubseteq a \;\&\; u \neq a) \;\&\; \;\&\; a \sqsubseteq a \;\&\; \textstyle\bigwedge_{d^* \in D} (\underline{d^*} \neq \underline{d} \rightarrow \neg \underline{d^*} \sqsupseteq a).$$

Suppose now that $M \vDash z \sqsubseteq x$. Then $M \vDash z \sqsubseteq a$, and we have that $M \vDash z \sqsubseteq \underline{t} \vee z = a$, that is, $M \vDash \forall z\, (z \sqsubseteq x \rightarrow z \sqsubseteq \underline{t} \vee z = x)$. Conversely, suppose $M \vDash z \sqsubseteq \underline{t}$. Then $M \vDash z \sqsubseteq a$. And if $M \vDash z = x$, then from $M \vDash a \sqsubseteq a \;\&\; x = a$, we have $M \vDash z \sqsubseteq x$. Hence also $M \vDash \forall z\, (z \sqsubseteq \underline{t} \vee z = x \rightarrow z \sqsubseteq x)$.

Case 2. $M \vDash \neg\exists! y\, B_d(\underline{t},y)$.

Then $M \vDash x = \underline{t}$. Suppose $M \vDash z \sqsubseteq \underline{t}$. Then, a fortiori, $M \vDash z \sqsubseteq \underline{t} \vee z = \underline{t}$. Conversely, suppose $M \vDash z = \underline{t}$. By 7.2(a), we have that $M \vDash \underline{t} \sqsubseteq \underline{t}$, hence $M \vDash z \sqsubseteq \underline{t}$. Therefore, also $M \vDash \forall z\, ( z \sqsubseteq \underline{t} \vee z = \underline{t} \rightarrow z \sqsubseteq \underline{t})$, that is, $M \vDash \forall z\, ( z \sqsubseteq \underline{t} \vee z = x \rightarrow z \sqsubseteq x)$.

So we have established that

$$M \vDash \forall x\, (C_d(\underline{t},x) \rightarrow \forall z\, (z \sqsubseteq x \leftrightarrow z \sqsubseteq \underline{t} \vee z = x)).$$

But then $VW_0 \vdash [(VW^*_{Rel}4)]^{\mathcal{L}}$. Analogously for ($VW^*_{Rel}5$). The proofs for ($VW^*_{Rel}6$), ($VW^*_{Rel}7$), ($VW^*_{Rel}8$) and ($VW^*_{Rel}9$) follow immediately from the definition of $C_d(x,y)$, and (VW3), (VW4), (VW3*) and (VW4*) follow from ($VW_02$) and ($VW_03$), resp.. This suffices to establish that $\mathcal{L}$ interprets $VW^*_{Rel}$ in $VW_0$. □

## §8. Very weak dyadic arithmetic with inverse

We now consider the theory VWDI formulated in the vocabulary $\mathcal{L}_{\sqsubseteq,f,c} = \{ \underline{c}_0, \underline{c}_1, \underline{c}_2,\ldots,f,\sqsubseteq \}$ with infinitely many individual constants $\underline{c}_0$, $\underline{c}_1$, $\underline{c}_2,\ldots$, one relational symbol $\sqsubseteq$ and a single unary function symbol $f$. We let M be any model of VWDI.

8.1 VWDI $\vdash \forall x,y\ (f(x) = f(y) \rightarrow x = y)$.

Proof: Assume $M \vDash f(x) = f(y)$ where $x,y \in M$. Then $M \vDash f(f(x)) = f(f(y))$, whence $M \vDash x = y$ follows from (VWDI4). □

Then we have:

Theorem 8.2 $VW_0 \leq_I$ VWDI.

Proof: We define a direct interpretation $P$ of $VW_0$ in VWDI as follows: each $\mathcal{L}_0$ individual constant $\underline{c}_i$, $i \in N$, is interpreted by the corresponding $\mathcal{L}_{\sqsubseteq,f,c}$-constant $\underline{c}_i$, assuming, as stated earlier, a fixed enumeration $s_0, s_1, s_2,\ldots$ of the strings in $\Sigma^*$, with $s_0 = e$, $s_1 = 1$, $s_2 = 2$. We let

$$x \sqsubseteq^P y \equiv: x \sqsubseteq y \text{ and } x \sqsupseteq^P y \equiv: f(x) \sqsubseteq f(y).$$

We first verify that VWDI $\vdash [(VW_0 3)]^P$, that is, that for each $t \in \Sigma^*$,

$$\text{VWDI} \vdash \forall x\ (f(x) \sqsubseteq f(\underline{t}) \leftrightarrow \bigvee_{s \sqsupseteq t} x = \underline{s}).$$

Here we argue analogously to 1.11 making use of (VWDI3), (VWDI2) and (VWDI4).

Note that $[(VW_04)]^P$ is (VWDI5). To complete the proof, it will suffice to show that VWDI $\vdash [(VW_05)]^P$. Let $x \in M$. From (VWDI5), substituting $f(x)$ for x and $\underline{f(t)}$ for $\underline{t}$ we have

$$(1) \quad M \vDash \forall y\,(y \sqsubset f(x) \leftrightarrow y \sqsubset \underline{f(t)}) \to \bigwedge_{d \in D} (\underline{d} \sqsubseteq f(f(x))\ \&\ \underline{d} \sqsubseteq f(\underline{f(t)}) \to \to f(x) = \underline{f(t)}) .$$

Assume (2) $M \vDash \forall y\,(f(y) \sqsubset f(x) \leftrightarrow f(y) \sqsubset f(\underline{t}))$ along with (3) $M \vDash \underline{d^*} \sqsubseteq x\ \&\ \underline{d^*} \sqsubseteq \underline{t}$ where $d^* \in D$. From (2) we have in particular that $M \vDash \forall y\,(f(f(y)) \sqsubset f(x) \leftrightarrow f(f(y)) \sqsubset f(\underline{t}))$. From (VWDI4) we obtain $M \vDash \forall y\,(y \sqsubset f(x) \leftrightarrow y \sqsubset f(\underline{t}))$, whence from (1) and (VWDI2) we derive

$$M \vDash \bigwedge_{d \in D} (\underline{d} \sqsubseteq f(f(x))\ \&\ \underline{d} \sqsubseteq f(f(\underline{t})) \to f(x) = f(\underline{t})).$$

From (VWDI4) again we obtain $M \vDash \bigwedge_{d \in D} (\underline{d} \sqsubseteq x\ \&\ \underline{d} \sqsubseteq \underline{t} \to f(x) = f(\underline{t}))$. Hence from (3) we have $M \vDash f(x) = f(\underline{t})$, and then $M \vDash x = \underline{t}$ follows by 8.1. Therefore, we derived

$$\text{VWDI} \vdash \forall x\,(\forall y\,(f(y) \sqsubset f(x) \leftrightarrow f(y) \sqsubset f(\underline{t})) \to \bigwedge_{d \in D} (\underline{d} \sqsubseteq x\ \&\ \underline{d} \sqsubseteq \underline{t} \to x = \underline{t})),$$

that is, $\text{VWDI} \vdash \forall x\,(\forall y\,(y \sqsupset x \leftrightarrow y \sqsupset \underline{t}) \to \bigwedge_{d \in D} (\underline{d} \sqsubseteq x\ \&\ \underline{d} \sqsubseteq \underline{t} \to x = \underline{t}))$. □

## §9. Putting It All Together

Theorem 9.1 $\quad VW^+ \equiv_I VW^* \equiv_I VW_0 \equiv_I VWDI \equiv_I WD \equiv_I WT \equiv_I R$

Proof: Note that every finite subset of axioms of VWDI has a finite model, hence by a theorem of Visser (see [9]) we have that VWDI $\leq_I$ R. For the same reason $VW^+ \leq_I R$. On the other hand, from Theorem 3.7 we have that WT $\leq_I$ VW$^+$, and from Theorem 5.6 that WD $\leq_I$ VW$^+$, whereas from Theorem 6.19 also WD $\leq_I$ VW*. From Theorem 7.5 it follows that VW* $\leq_I$ VW$_0$, while VW$_0$ $\leq_I$ VWDI by Theorem 8.2. On other hand, R $\leq_I$ WD by [6] and R $\leq_I$ WT by [3]. □

We now consider the structure $\mathfrak{S}_1 = < \Sigma^*, 1, 2, \sqsubseteq, \sqsupseteq >$ interpreting the vocabulary $\mathcal{L}_1 = \{ \underline{1}, \underline{2}, \sqsubseteq, \sqsupseteq \}$ in the usual way and its theory $T_1 = Th(\Sigma^*, \mathcal{L}_1)$.

Theorem 9.2 Binary concatenation is definable in T$_1$.

Proof: Consider the $\mathcal{L}^*$ formula Conc(x,y) from §6. Using the formulae $A_d(x,y)$, d = 1,2,-1,-2, from §7, we can eliminate all occurrences of the operation symbols x^1, x^2, $1_\vee x$ and $2_\vee x$. Finally, occurrences of the constant $\underline{e}$ can be eliminated using the predicate $E(x) \equiv \forall z\, (z \sqsubseteq x \rightarrow z = x)$. The resulting formula Conc*(x,y) is an $\mathcal{L}_1$ formula, which, on account of 6.17, defines the binary concatenation operation x ̂ y in $\mathfrak{S}_1$. □

Let F(x,y) abbreviate the $\mathcal{L}^*$ formula

$(x = \underline{e}\ \&\ y = \underline{e})$ v $\exists u,v\ \exists z \subseteq_p u$ (BlockChain(u,x) & BlockChain(v,y) &

& Tlmp(z,u) & Tlmp(z,v) & $\forall w_1 \sqsubseteq u\ \forall w_2 \sqsubseteq v\ \forall z_1,z_2 \sqsubseteq z\ (z_1 \underline{\in}_{w(1)} u\ \&\ z_2 \underline{\in}_{w(2)} v$ &

& PsTC($z_1$,$z_2$,z^21) $\rightarrow \bigwedge_{d \in D} (w_1$^$\underline{d} \sqsubseteq u \leftrightarrow w_2$^$\underline{d} \sqsubseteq v$)).

9.3 For each $s \in \Sigma^*$, $\mathfrak{S}^* \vDash F(\underline{s}, \underline{s^{-1}})$.

Proof: We argue by induction on lh(s). The case where s = e is immediate. Assuming the claim holds for s, we consider $s^\frown d$. Note that $(s^\frown d)^{-1} = d^\frown s^{-1}$. From the induction hypothesis there are $u^*, v^*, z^* \in \Sigma^*$, such that

$$\mathfrak{S}^* \vDash \text{BlockChain}(\underline{u^*},\underline{s}) \ \& \ \text{BlockChain}(\underline{v^*},\underline{s^{-1}}),$$

where $\mathfrak{S}^* \vDash \text{Tlmp}(\underline{z^*},\underline{u^*}) \ \& \ \text{Tlmp}(\underline{z^*},\underline{v^*})$. We also have that

$\mathfrak{S}^* \vDash \forall w_1 \sqsubseteq \underline{u^*} \ \forall w_2 \sqsubseteq \underline{v^*} \ \forall z_1,z_2 \sqsubseteq \underline{z^*} \ (\text{PsTC}(z_1,z_2,\underline{z^*}\hat{\ }21) \ \&$

$\& \ z_1 \in_{w(1)} \underline{u^*} \ \& \ z_2 \in_{w(2)} \underline{v^*} \ \rightarrow \ \bigwedge_{d \in D} (z_1\hat{\ }\underline{d} \sqsubseteq \underline{u^*} \leftrightarrow z_2\hat{\ }\underline{d} \sqsubseteq \underline{v^*})).$

By 6.14 we have that $u^* = {}^\frown_{1\leq i\leq k} (21)^{i\frown} d_{1(i)})$ and $v^* = {}^\frown_{1\leq i\leq k} (21)^{i\frown} d_{2(i)})$ for some unique $d_{1(i)}$ and $d_{2(i)}$, $1 \leq i \leq k$, such that $s = {}^\frown_{1\leq i\leq k} d_{1(i)}$ and $s^{-1} = {}^\frown_{1\leq i\leq k} d_{2(i)}$ and $z^* = (21)^k$. From the induction hypothesis we have that $d_{1(j)} = d_{2(k+1-j)}$ for each j, $1\leq j \leq k$. Now, let u, v be the blockchains for $s^\frown d$ and $(s^\frown d)^{-1}$, resp.. Then

$$u = u^{*\frown} (21)^{k+1\frown} d \ \text{ and } \ v = (21)^\frown d^\frown ({}^\frown_{1\leq i\leq k} ((21)^{i+1\frown} d_{2(i)},)),$$

where $d_0$, the last digit of $s^{-1}$ (= $d_{2(k)}$), is the first digit of s. We need to show that u, v and $z = (21)^{k+1}$ have the required properties.

Suppose $z_1,z_2 \sqsubseteq (21)^{k+1}$ where

$\mathfrak{S}^* \vDash \underline{z_1} \in_{w(1)} \underline{u} \ \& \ \underline{z_2} \in_{w(2)} \underline{v} \ \& \ \text{PsTC}(\underline{z_1},\underline{z_2},(21)^{k+2})$ with $w_1 \sqsubseteq u$ and $w_2 \sqsubseteq v$.

Then $z_1{}^\frown z_2 = (21)^{k+2}$, and by 6.14, we have that $z_1 = (21)^i$ and $z_2 = (21)^j$ for some i, $j \leq k+1$ such that $i+j = k+2$ where $z_1 \sqsupseteq w_1$ and $z_2 \sqsupseteq w_2$. If i = 1, then j = k+1, and we have that $(21)^\frown d_0 \sqsubseteq u$ and $(21)^{j\frown} d_{2(k)} = (21)^{k+1\frown} d_0 \sqsupseteq v$, as needed. If i = k+1, then j = 1, and we have that $(21)^{k+1\frown} d \sqsupseteq u$ and $(21)^\frown d \sqsubseteq v$, as needed.

If $1 < i < k+1$, then for any $d^* \in D$, and $j^*$ such that $i+j^* = k+1$,

$(21)^i \,\hat{}\, d^* \subseteq_p u \;\Leftrightarrow\; (21)^i \,\hat{}\, d_{1(i)} \subseteq_p u^* \;\Leftrightarrow\; (21)^{j^*} \,\hat{}\, d_{1(i)} \subseteq_p v^* \Leftrightarrow$

$\Leftrightarrow\; (21)^{j^*} \,\hat{}\, d_{2(k+1-i)} \subseteq_p v^* \Leftrightarrow\; (21)^{j^*} \,\hat{}\, d_{2(j^*)} \subseteq_p v^* \Leftrightarrow\; (21)^{j^*+1} \,\hat{}\, d_{2(j^*)} \subseteq_p v.$

But $j^* = k+1-i = k+1-(k+2-j) = j-1$. Hence we obtain

$$w_1 \,\hat{}\, d^* \sqsubseteq u \;\Leftrightarrow\; w_2 \,\hat{}\, d^* \sqsubseteq v,$$

as needed. □

We further have:

9.4 (a) $\mathfrak{S}^* \vDash \forall x,y,z\ (F(x,y)\ \&\ F(x,z) \rightarrow y = z)$.

(b) $\mathfrak{S}^* \vDash \forall x,y\ (F(x,y) \rightarrow F(y,x))$.

(c ) $\mathfrak{S}^* \vDash \forall x,y,u,v\ (F(x,u)\ \&\ F(y,v)\ \&\ u = v \rightarrow x = y)$.

<u>Proof</u>: (a) Assume $\mathfrak{S}^* \vDash F(\underline{s},\underline{u})\ \&\ F(\underline{s},\underline{v})$ for s, u, v $\in \Sigma^*$. Then there are blockchains

$\frown_{1\le i\le k} ((21)^i \,\hat{}\, d_{1(i)})$ and $\frown_{1\le i\le k} ((21)^i \,\hat{}\, d_{2(i)})$ and $\frown_{1\le i\le k} ((21)^i \,\hat{}\, d_{3(i)})$

for s, u, v, resp., where $s = \frown_{1\le i\le k}\ d_{1(i)}$, $u = \frown_{1\le i\le k}\ d_{2(i)}$ and $v = \frown_{1\le i\le k}\ d_{3(i)}$, and, furthermore, for each i, $1 \le i \le k$, $d_{2(k+1-i)} = d_{3(k+1-i)}$, and hence, for each j, where $k \ge j \ge 1$ and $j = k+1-i$, $d_{2(j)} = d_{3(j)}$. But then

$$u = \frown_{1\le j\le k} d_{2(i)} = \frown_{1\le j\le k}\ d_{3(i)} = v,$$

as claimed. For (b), assume $\mathfrak{S} \vDash F(\underline{s},\underline{t})$ for s, t $\in \Sigma^*$. Then there are blockchains $\frown_{1\le i\le k} ((21)^i \,\hat{}\, d_{1(i)})$ and $\frown_{1\le i\le k} ((21)^i \,\hat{}\, d_{2(i)})$ for s, t, resp., where $s = \frown_{1\le i\le k}\ d_{1(i)}$, $t = \frown_{1\le i\le k}\ d_{2(i)}$ and for each i, $1 \le i \le k$, $d_{1(i)} = d_{2(k+1-i)}$. For each j such that $i+j = k+1$ we have $j = k+1-i$ and $i = k+1-j$. Hence $d_{2(i)} = d_{1(k+1-j)}$ for each j, $1 \le j \le k+1$. Since

$\mathfrak{S}^* \vDash PsTC((21)^i,(21)^j,(21)^k) \;\Leftrightarrow\; \mathfrak{S}^* \vDash PsTC((21)^j,(21)^i,(21)^k)$,

it follows from definitions that $\mathfrak{S}^* \vDash F(\underline{t},\underline{s})$.

Finally, for (c ), assume $\mathfrak{S}^* \vDash F(\underline{s},\underline{u})$ & $F(\underline{t},\underline{v})$ & $\underline{u} = \underline{v}$ for s, t, u, v $\in \Sigma^*$. By 8.3 we have that

(1) $\mathfrak{S}^* \vDash F(\underline{s},\underline{s^{-1}})$ & $F(\underline{t},\underline{t^{-1}})$.

From (a) and the hypothesis it follows that $\mathfrak{S}^* \vDash \underline{u}=\underline{s^{-1}}$ & $\underline{v}=\underline{t^{-1}}$. Hence from hypothesis $\mathfrak{S}^* \vDash \underline{u} = \underline{v}$ we have that $\mathfrak{S}^* \vDash \underline{s^{-1}} = \underline{t^{-1}}$. But from (1) by (b) it follows that $\mathfrak{S}^* \vDash F(\underline{s^{-1}},\underline{s})$ & $F(\underline{t^{-1}},\underline{t})$, hence also $\mathfrak{S}^* \vDash F(\underline{s^{-1}},\underline{s})$ & $F(\underline{s^{-1}},\underline{t})$. But then $\mathfrak{S}^* \vDash \underline{s} = \underline{t}$ follows from (a). □

Let $\mathcal{L}_c = \{*,a,b\}$ where * is a binary operation symbol for the concatenation operation and a, b, individual constants. We then have:

<u>Theorem 9.5</u> $Th(\Sigma^*,\mathcal{L}_c) \equiv_I Th(\Sigma^*,\mathcal{L}^*) \equiv_I Th(\Sigma^*,\mathcal{L}_1) \equiv_I Th(\Sigma^*,\mathcal{L}_{\sqsubseteq,f,c})$.

<u>Proof</u>: To show that $Th(\Sigma^*,\mathcal{L}_{\sqsubseteq,f,c}) \leq_I Th(\Sigma^*,\mathcal{L}_1)$ we use the formula F(x,y) and then eliminate all occurrences of $\underline{e}$ and the operation symbols of $\mathcal{L}^*$ as described in the proof of 9.2. The resulting formulae interpret the non-logical vocabulary of $\mathcal{L}_{\sqsubseteq,f,c}$ in $\mathcal{L}_1$ on account of 9.3. Since $\sqsupseteq$ is definable in terms of $\sqsubseteq$ and $f$ in $\mathcal{L}_{\sqsubseteq,f,c}$, the converse $Th(\Sigma^*,\mathcal{L}_1) \leq_I Th(\Sigma^*,\mathcal{L}_{\sqsubseteq,f,c})$ also holds. That $Th(\Sigma^*,\mathcal{L}_1) \leq_I Th(\Sigma^*,\mathcal{L}_c)$ follows from the fact '$x \sqsubseteq y$' and '$x \sqsupseteq y$' are defined in terms of * by '$\exists z\, x*z = y$' and '$\exists z\, z*x = y$', resp., and $\underline{1}$ and $\underline{2}$ are interpreted by $\underline{a}$ and $\underline{b}$, resp.. That $Th(\Sigma^*,\mathcal{L}_c) \leq_I Th(\Sigma^*,\mathcal{L}^*)$ follows from 6.17, and $Th(\Sigma^*,\mathcal{L}^*) \leq_I Th(\Sigma^*,\mathcal{L}_1)$ follows modulo elimination of the occurrences of $\underline{e}$ and the operation symbols of $\mathcal{L}^*$ as described above. □

Remark: Karlov [5] has proved that $T_1 = \mathrm{Th}(\Sigma^*, \mathcal{L}_1)$, where $\mathcal{L}_1 = \{ \underline{1}, \underline{2}, \sqsubseteq, \sqsupseteq \}$, is mutually interpretable with elementary arithmetic $\mathrm{Th}(\mathbb{N}, +, \cdot)$, and that so is $\mathrm{Th}(\Sigma^*, \mathcal{L}_{\sqsubseteq, f})$. Hence we conclude that the same holds of the other theories listed in Theorem 9.5. Neither one of $\sqsupseteq$, $\sqsubseteq$ can be omitted from $\mathcal{L}_1$, nor $f$ from $\mathcal{L}_{\sqsubseteq, f, \mathrm{c}}$, because the resulting theories are known to be decidable.